\documentclass[12pt]{iopart}
\usepackage{iopams}

\newcommand{\Z}{\mathbb{Z}}
\newcommand{\R}{\mathbb{R}}

\newcommand{\mrm}[1]{\mathrm{#1}}
\newcommand{\mcl}[1]{\mathcal{#1}}

\usepackage{enumitem}
\usepackage{float}
\usepackage{graphicx}
\usepackage{subcaption}

\usepackage{tikz}
\usetikzlibrary{positioning}
\usetikzlibrary{arrows}
\usetikzlibrary{shapes.geometric}
\usetikzlibrary{arrows}
\usetikzlibrary{decorations.markings}
\usetikzlibrary{arrows.meta}

\usepackage{xcolor}
\usepackage{cite}
\usepackage{xifthen}
\usepackage{multirow}

\definecolor{xi_orange}{HTML}{FF7F00}
\definecolor{xi_green}{HTML}{66CC00}
\definecolor{xi_blue}{HTML}{0080FF}
\definecolor{xi_pink}{HTML}{FF0080}
\definecolor{xi_dark_green}{HTML}{00B359}
\definecolor{xi_purple}{HTML}{7F00FF}

\definecolor{mytextgreen}{HTML}{28A428}
\definecolor{mytextblue}{HTML}{3A53F8}
\definecolor{mytextred}{HTML}{F83A53}

\usepackage{mathtools}

\usepackage{amsthm}
\theoremstyle{definition}
\newtheorem{definition}{Definition}[section]

\theoremstyle{plain}

\newcommand{\poinin}[2][]{\ifthenelse{\isempty{#1}}{\mathbf{H}_{\mrm{#2}}^{\mathrm{in}}}{\mathbf{H}_{\mrm{#2}}^{\mathrm{in},\mrm{#1}}}}
\newcommand{\poinout}[2][]{\ifthenelse{\isempty{#1}}{\mathbf{H}_{\mrm{#2}}^{\mathrm{out}}}{\mathbf{H}_{\mrm{#2}}^{\mathrm{out},\mrm{#1}}}}

\newcommand{\cyc}[1]{\mathcal{C}_{#1}}
\newcommand{\net}[1]{\mathcal{N}_{#1}}

\newcommand{\invar}[1]{\mathrm{in},\mrm{#1}}
\newcommand{\outvar}[1]{\mathrm{out},\mrm{#1}}

\newcommand{\order}[1]{\mathcal{O}\left(#1\right)}

\usepackage{stackengine}
\stackMath
\newcommand\reallywidehat[1]{%
\savestack{\tmpbox}{\stretchto{%
  \scaleto{%
    \scalerel*[\widthof{\ensuremath{#1}}]{\kern-.6pt\bigwedge\kern-.6pt}%
    {\rule[-\textheight/2]{1ex}{\textheight}}
  }{\textheight}%
  }{0.5ex}}%
  \stackon[1pt]{#1}{\tmpbox}%
}

\def\Rminusscale{0.7}
\newcommand{\negativeR}[1]{\R_{\scalebox{\Rminusscale}[1.0]{--}}^{#1}}

\let\originalleft\left
\let\originalright\right
\renewcommand{\left}{\mathopen{}\mathclose\bgroup\originalleft}
\renewcommand{\right}{\aftergroup\egroup\originalright}

\newcommand{\rmP}{\mrm{P}}
\newcommand{\rmD}{\mrm{D}}
\newcommand{\rms}{\mrm{s}}
\newcommand{\rmc}{\mrm{c}}

\newcommand{\rmA}{\mrm{A}}
\newcommand{\rmB}{\mrm{B}}

\newcommand{\rmX}{\mrm{X}}
\newcommand{\rmY}{\mrm{Y}}
\newcommand{\rmXY}{\mrm{X}\mrm{Y}}
\newcommand{\rmC}{\mrm{C}}
\newcommand{\rmks}{\mrm{KS}}

\usepackage{amsmath}

\begin{document}

\title[]{Continuity of projected maps of heteroclinic networks in \(\mathbb{R}^{4}\)}

\author{David C Groothuizen Dijkema\(^{1}\), Claire M Postlethwaite\(^{2}\), Alastair M Rucklidge\(^{3}\)}

\address{\(^{1}\) School of Mathematical Sciences, University College Cork --- National University of Ireland, Cork T12 XF62, Ireland}
\address{\(^{2}\) Department of Mathematics, University of Auckland, Auckland 1010, New Zealand}
\address{\(^{3}\) School of Mathematics, University of Leeds, Leeds LS2 9JT, United Kingdom}
\ead{DGroothuizenDijkema@ucc.ie}
\vspace{10pt}
\begin{indented}
\item[]June 2026
\end{indented}

\begin{abstract}
  Stability of robust heteroclinic cycles and networks is typically studied by constructing return maps and analysing their associated transition matrices. This analysis can be simplified with the network's \textit{projected map}, derived by projecting the linear action of the transition matrix onto a simplex. This projection produces a piecewise-smooth map, which is one-dimensional for heteroclinic networks in \(\R^{4}\). We consider two such networks, the Kirk--Silber network and the \(\Delta\)-clique network. For the Kirk--Silber network, the projected map is discontinuous on its switching manifold, while it is continuous for the \(\Delta\)-clique network. In this paper, we address the dynamical phenomena that produce a discontinuity in the case of the Kirk--Silber network, and explain the value of the projected map at the switching manifold. We construct a \textit{completed} return map near both networks, which captures the behaviour of all trajectories that begin near the network but may move away from it temporarily. We show that the discontinuity in the projected map of the Kirk--Silber network emerges due to two phenomena: first, there exists a discontinuity in a component of the completed return map in the limit as trajectories approach the Kirk--Silber network, as a result of the presence of a separatrix near the network, and, second, the procedure that defines the projected map. For the \(\Delta\)-clique network, there is no such separatrix, and so no such discontinuity emerges. This analysis is a necessary step towards understanding the more complicated dynamics observed near heteroclinic networks of five or more equilibria.
\end{abstract}

\vspace{2pc}
\noindent 37C29, 34C37 

\submitto{\NL}

\section{Introduction}\label{sec:intro}

Heteroclinic cycles and networks are flow-invariant structures in dynamical systems composed of hyperbolic saddle equilibria and heteroclinic orbits between them. In systems with appropriate invariant subspaces, these structures can be robust. The behaviour of trajectories near an attracting cycle or network is typically intermittent: trajectories spend a long period of time in a small neighbourhood of one equilibrium before rapidly switching to the next.

Heteroclinic cycles and networks are typically analysed using return maps to cross-sections defined in a small neighbourhood of each equilibrium. These return maps are sometimes not defined on the entire cross-section: we may have to exclude cusp-shaped sets of initial conditions that give rise to trajectories that do not remain sufficiently close to the network for certain linearisations to be valid approximations of the dynamics near the network. These return maps can be transformed into logarithmic coordinates, and, for certain classes of heteroclinic cycles---including those that are of type \(Z\) in the terminology of Podvigina \cite{podvigina_2012} and those that are \textit{quasi-simple} in the terminology of Garrido da Silva and Castro \cite{garrido_da_silva_castro_2019}---the action of the leading order approximation of the map becomes linear, described by a \textit{transition matrix}. Podvigina \cite{podvigina_2012} provides sufficient conditions for a heteroclinic cycle to be \textit{fragmentarily asymptotically stable} based on the eigenvalues and eigenvectors of its transition matrices.

\begin{figure}
  \centering
  \newlength{\thickness}
  \setlength{\thickness}{0.5mm}

  \newlength{\xdisp}
  \setlength{\xdisp}{1.3cm}
  \newlength{\ydisp}
  \setlength{\ydisp}{1cm}
  \newlength{\ydispp}
  \setlength{\ydispp}{2.2cm}

  \begin{subfigure}[t]{0.45\linewidth}
    \centering
    \begin{tikzpicture}
      \node [circle,minimum size=2.5mm,inner sep=0pt,draw=xi_pink,fill=xi_pink,    label={[font=\normalsize, label distance=-0.5mm]above:\(A\)}]                                        (one)   {};
      \node [circle,minimum size=2.5mm,inner sep=0pt,draw=xi_green,fill=xi_green,  label={[font=\normalsize, label distance=-0.5mm]below:\(B\)},below = \ydispp of one]                 (two)   {};
      \node [circle,minimum size=2.5mm,inner sep=0pt,draw=xi_orange,fill=xi_orange,    label={[font=\normalsize, label distance=-1.0mm]left: \(Y\)},below left  = \ydisp and \xdisp of one] (three) {};
      \node [circle,minimum size=2.5mm,inner sep=0pt,draw=xi_blue,fill=xi_blue,label={[font=\normalsize, label distance=-0.5mm]right:\(X\)},below right = \ydisp and \xdisp of one] (four)  {};

      \draw [-{Straight Barb},line width=\thickness,shorten >=5pt,shorten <=5pt] (one)   to [] node[below] {} (two);
      \draw [-{Straight Barb},line width=\thickness,shorten >=5pt,shorten <=5pt] (two)   to [] node[below] {} (three);
      \draw [-{Straight Barb},line width=\thickness,shorten >=5pt,shorten <=5pt] (three) to [] node[below] {} (one);
      \draw [-{Straight Barb},line width=\thickness,shorten >=5pt,shorten <=5pt] (two)  to [] node[below] {} (four);
      \draw [-{Straight Barb},line width=\thickness,shorten >=5pt,shorten <=5pt] (four) to [] node[below] {} (one);
    \end{tikzpicture}
    \caption{The Kirk--Silber network.}
    \label{fig:R4_networks:kirk_silber_net}
  \end{subfigure}%
  \hfill%
  \begin{subfigure}[t]{0.45\linewidth}
    \centering
    \begin{tikzpicture}
      \node [circle,minimum size=2.5mm,inner sep=0pt,draw=xi_orange,fill=xi_orange,label={[font=\normalsize, label distance=-0.5mm]above:\(Y\)}]                                         (four)  {};
      \node [circle,minimum size=2.5mm,inner sep=0pt,draw=xi_green,fill=xi_green,  label={[font=\normalsize, label distance=-0.5mm]below:\(B\)},below = \ydispp of four]                 (two)   {};
      \node [circle,minimum size=2.5mm,inner sep=0pt,draw=xi_pink,fill=xi_pink,    label={[font=\normalsize, label distance=-1.0mm]left: \(A\)},below left  = \ydisp and \xdisp of four] (one)   {};
      \node [circle,minimum size=2.5mm,inner sep=0pt,draw=xi_blue,fill=xi_blue,    label={[font=\normalsize, label distance=-0.5mm]right:\(X\)},below right = \ydisp and \xdisp of four] (three) {};

      \draw [-{Straight Barb},line width=\thickness,shorten >=5pt,shorten <=5pt] (one)   to [] node[below] {} (two);
      \draw [-{Straight Barb},line width=\thickness,shorten >=5pt,shorten <=5pt] (two)   to [] node[below] {} (three);
      \draw [-{Straight Barb},line width=\thickness,shorten >=5pt,shorten <=5pt] (two)   to [] node[below] {} (four);
      \draw [-{Straight Barb},line width=\thickness,shorten >=5pt,shorten <=5pt] (three) to [] node[below] {} (four);
      \draw [-{Straight Barb},line width=\thickness,shorten >=5pt,shorten <=5pt] (four)  to [] node[below] {} (one);
    \end{tikzpicture}
    \caption{The \(\Delta\)-clique network.}
    \label{fig:R4_networks:delta_clique_net}
  \end{subfigure}%
  \caption{Diagrammatic representations of the two heteroclinic networks between four equilibria in \(\R^{4}\) that we study in this paper. Coloured vertices represent hyperbolic saddle equilibria, and directed edges represent robust heteroclinic orbits.}
  \label{fig:R4_networks}
\end{figure}
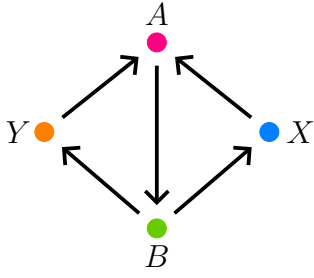
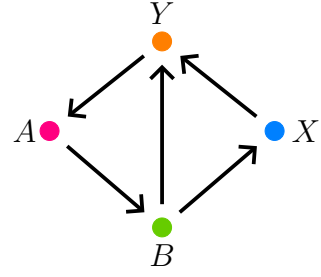

In \cite{article}, we developed an extension to this methodology, called the \textit{projected map}, which is a piecewise-smooth discrete map defined by the induced action of certain transition matrices on a simplex. The projected map can be used to determine when a cycle is fragmentarily asymptotically stable, and to more easily analyse the behaviour of trajectories near a network. For networks in \(\R^{4}\), this map is one-dimensional. In \cite{article}, we used the projected map to classify the dynamics of trajectories near three heteroclinic networks in \(\R^{4}\) that consist of four equilibria. We consider two of these heteroclinic networks in this paper: the Kirk--Silber network and the \(\Delta\)-clique network. (See \Fref{fig:R4_networks}.)

The \textit{switching manifold} of a piecewise map is the boundary between the domains of definition of its components. Although the projected map of a heteroclinic network is not defined on its switching manifold, we can consider the limiting value of the map at this point. In the case of the \(\Delta\)-clique network, we observed in \cite{article} that the left-hand and right-hand limits are equal, and so we say the projected map is continuous. However, for the Kirk--Silber network, these two limits are generically not equal, and so the projected map is discontinuous. While in \cite{article} we used the projected map to classify the dynamics of trajectories near three heteroclinic networks in \(\R^{4}\), in this paper, our focus is to further our understanding of the properties of the projected map by understanding the dynamical properties of the network which result in a continuous or discontinuous projected map.

We focus on the continuity of the projected map in this paper because the continuity or discontinuity of a general piecewise-smooth map at its switching manifold can greatly affect the resulting dynamics seen in the system \cite{zhusubaliyev_mosekilde_2003,di_bernado_budd_champneys_koqalczyk_2008,simpson_2010,avrutin_gardini_sushko_tramontana_2019}. For example, Arnold tongues in piecewise-linear \textit{continuous} maps often form chains with codimension-\(2\) points of zero width, and the tongues are likened to a ``string of sausages'' \cite{yang_hao_1987,simpson_2017}. However, in a piecewise-smooth \textit{discontinuous} map, the rotation numbers of periodic orbits form a devil's staircase in parameter space \cite{di_bernado_budd_champneys_koqalczyk_2008}. Therefore, in the study of heteroclinic networks, we would like to understand what properties of the original dynamical system result in a continuous or discontinuous projected map. This paper is an important first step in studying the (dis)continuity of a projected map in general. We expect that this analysis will allow an understanding of the continuity of projected maps for heteroclinic networks in higher dimensions, such as the Rock-Paper-Scissors-Lizard-Spock network \cite{postlethwaite_rucklidge_2022}. This understanding may then assist in the proof of the existence of the complicated structures and patterns observed in the bifurcation sets of regular and irregular cycling near such networks.

The return map used to construct the projected map is not defined for the set of trajectories that correspond to the switching manifold. Therefore, we cannot investigate the dynamical origin of this discontinuity from this return map. In this paper, we instead use the methodology in \cite{kirk_lane_postlethwaite_rucklidge_silber_2010} and \cite{kirk_postlethwaite_rucklidge_2012} to construct return maps near the Kirk--Silber and \(\Delta\)-clique networks that include all trajectories near the network. We refer to such a return map as a \textit{completed} return map. As it does not exclude any trajectories near the network, we can study the dynamics of trajectories that correspond to the switching manifold. We are thus able to determine the features of the dynamical system that cause the (dis)continuity of the projected map. Specifically, we are able to show that the presence of a separatrix near the Kirk--Silber network creates, in the limit as a trajectory approaches the network, a discontinuity in the completed return map. In contrast, in the case of the \(\Delta\)-clique network, no such discontinuity emerges. We are then able to demonstrate how this continuity or discontinuity manifests in the value of the projected map at the switching manifold, as the limiting value of the projected map at the switching manifold is found in the limit as trajectories approach the network.

The remainder of this paper proceeds as follows. In \sref{sec:problem}, we define both the Kirk--Silber and \(\Delta\)-clique networks, summarise how the projected map of each network is constructed, and demonstrate their continuity or discontinuity. In \sref{sec:problem:ssec:limit}, we highlight a particular part of the construction of the projected map important to understanding the emergence of the discontinuity. We proceed in \sref{sec:ret_map} to define cross-sections near both networks and construct local and global maps between these sections that account for all trajectories. Although we cannot construct these completed maps explicitly, we can derive certain properties that we can then use to analyse the behaviour of trajectories near the networks. Then, in \sref{sec:anal:ssec:proj}, we define projected maps from these completed return maps. In \sref{sec:anal:ssec:calc}, we calculate an approximation of the completed return map for certain trajectories sufficiently close to the network. Lastly, in \sref{sec:anal:ssec:expl}, we explain how the discontinuity of the projected map of the Kirk--Silber network emerges, and why the projected map of the \(\Delta\)-clique network is continuous.

\Sref{sec:disc} concludes with a particular view towards the consequences of this work in the study of heteroclinic networks.

\section{Problem statement}\label{sec:problem}

We begin by defining the two heteroclinic networks we study in this paper, reviewing the projected map of each network, and discussing the continuity of these projected maps.

\subsection{The heteroclinic networks}

In this paper, we are interested in dynamical systems defined by a system of ordinary differential equations
\begin{equation}\label{eqn:dyn_sys}
  \dot{x}=f(x),
\end{equation}
where \(x=\left(x_{1},x_{2},x_{3},y_{3}\right)\in\R^{4}\) and \(f\colon\R^{4}\to\R^{4}\) is a sufficiently smooth vector field. For a point \(x_{0}\in\R^{4}\), let \(\phi\left(x_{0},t\right)\) be the flow generated by \eqref{eqn:dyn_sys} through \(x_{0}\). Let \(\xi_{j}\) and \(\xi_{k}\) be two equilibria of the system. Then, a solution \(\phi\left(x,t\right)\) is a \textit{heteroclinic orbit} from \(\xi_{j}\) to \(\xi_{k}\) if \(\phi\left(x,t\right)\to\xi_{j}\) as \(t\to-\infty\) and \(\phi\left(x,t\right)\to\xi_{k}\) as \(t\to\infty\). We also write \(\xi_{j}\to\xi_{k}\) for a heteroclinic orbit. We note that, by this definition as a specific solution to the ODEs \eqref{eqn:dyn_sys}, a heteroclinic orbit is a specific trajectory in phase-space. A heteroclinic orbit is thus also a one-dimensional submanifold of \(W^{\mathrm{u}}\left({\xi_{j}}\right)\cap W^{\mathrm{s}}\left(\xi_{k}\right)\neq\emptyset\). For convenience, we say that \(\xi_{j}\to\xi_{k}\) is an \textit{outgoing} orbit of \(\xi_{j}\) and \textit{incoming} orbit of \(\xi_{k}\).

\begin{definition}\label{def:het_cyc}
  A \textit{heteroclinic cycle} \(\cyc{}\subset\R^{4}\) is the union of a set \(\left\{\xi_{1},\dots,\xi_{m}\right\}\) of equilibria of \eqref{eqn:dyn_sys}, and a set \(\left\{\phi_{1}(t),\dots,\phi_{m}(t)\right\}\) of heteroclinic orbits, such that \(\phi_{j}(t)\) is a heteroclinic orbit from \(\xi_{j}\) to \(\xi_{j+1}\), where all subscripts are taken modulo \(m\).
\end{definition}
\begin{definition}\label{def:het_net}
  Let \(\cyc{1},\ldots,\cyc{s}\), for \(s\geq 2\), be a collection of heteroclinic cycles. Then the union \(\net{}=\bigcup_{j=1}^{s}\cyc{j}\) is a \textit{heteroclinic network} if, for every pair of equilibria, there exists a sequence of heteroclinic orbits connecting them; that is, if \(\xi_{\rmA},\xi_{\rmB}\in\net{}\), then there exists a set of equilibria \(\left\{\xi_{\ell_{1}},\dots,\xi_{\ell_{p}}\right\}\subseteq\net{}\) and heteroclinic orbits \(\left\{\phi_{\ell_{1}\ell_{2}}(t),\dots,\phi_{\ell_{p}\ell_{1}}(t)\right\}\subseteq\net{}\), such that \(\xi_{\ell_{1}}\equiv\xi_{\rmA}\), \(\xi_{\ell_{p}}\equiv\xi_{\rmB}\), and \(\phi_{\ell_{i}\ell_{j}}(t)\) is a heteroclinic orbit from \(\xi_{\ell_{i}}\) to \(\xi_{\ell_{j}}\).
\end{definition}

We consider the system \eqref{eqn:dyn_sys} under two different cases, each of which contains a different heteroclinic network. These two networks are composed of four different hyperbolic saddle equilibria, \(A\), \(B\), \(X\), and \(Y\). We assume that each of these equilibria lie at unit distance from the origin on each of the four coordinate axes in the positive orthant. We construct three different cycles from these equilibria:
\begin{enumerate}
  \item The \(\cyc{X}\) cycle, \(A\to B\to X\to A\);
  \item The \(\cyc{Y}\) cycle, \(A\to B\to Y\to A\); and
  \item The \(\cyc{XY}\) cycle, \(A\to B\to X\to Y\to A\).
\end{enumerate}
The cycles \(\cyc{X}\) and \(\cyc{Y}\) are both of type \(B_{3}^{-}\) in the classification of cycles in \(\R^{4}\) given in \cite{krupa_melbourne_2004}, and the cycle \(\cyc{XY}\) is of type \(C_{4}^{-}\).
\begin{enumerate}
  \item The Kirk--Silber network \cite{kirk_silber_1994}, \(\net{KS}=\cyc{X}\cup\mathcal{C}_{\rmY}\); and
  \item The \(\Delta\)-clique network, \(\net{\Delta}=\cyc{XY}\cup\mathcal{C}_{\rmY}\).
\end{enumerate}
We present diagrammatic representations of these networks in \Fref{fig:R4_networks}. The directed triangle formed by the orbits \(B\to Y\) and \(B\to X\to Y\) was named a \(\Delta\)-clique in \cite{ashwin_castro_lohse_2020}, and as such we name it the \(\Delta\)-clique network.

The following vector field \eqref{eqn:dyn_sys} is produced with the simplex realisation of Ashwin and Postlethwaite \cite{ashwin_postlethwaite_2013}, and is an example of a system of ODEs that contains these networks:
\begin{equation}\label{eqn:ODEs}
  \begin{aligned}
    \dot{x}_{1}&=x_{1}\left(1 - \left\lVert x\right\rVert_{2}^{2} + \alpha_{21}x_{2}^{2} + \alpha_{31}x_{3}^{2} + \alpha_{41}y_{3}^{2}\right),\\
    \dot{x}_{2}&=x_{2}\left(1 - \left\lVert x\right\rVert_{2}^{2} + \alpha_{12}x_{1}^{2} + \alpha_{32}x_{3}^{2} + \alpha_{42}y_{3}^{2}\right),\\
    \dot{x}_{3}&=x_{3}\left(1 - \left\lVert x\right\rVert_{2}^{2} + \alpha_{13}x_{1}^{2} + \alpha_{23}x_{2}^{2} + \alpha_{43}y_{3}^{2}\right),\\
    \dot{y}_{3}&=y_{3}\left(1 - \left\lVert x\right\rVert_{2}^{2} + \alpha_{14}x_{1}^{2} + \alpha_{24}x_{2}^{2} + \alpha_{34}x_{3}^{2}\right),
  \end{aligned}
\end{equation}
where each \(\alpha_{jk}\in\R\) is a parameter. We impose the genericity conditions that all \(\alpha_{jk}\neq0\). It is straightforward to verify that these equations contain four hyperbolic saddle equilibria at unit distance from the origin in the positive orthant, one on each of the four coordinates axes. Specific choices of the sign of the \(\alpha_{jk}\) terms ensure that certain equilibria are sinks or saddles for the flow restricted to corresponding flow-invariant coordinate planes, and thus the existence of robust heteroclinic orbits between two given equilibria. As such, they can be chosen in such a way as to give either of the two networks we consider. It is also straightforward to verify that the above ODEs are equivariant with respect to the group \(\Z_{2}^{4}\), generated by reflections across each of the coordinate hyperplanes. This equivariance ensures the flow-invariance of the coordinate axes, coordinate planes, and coordinate hyperplanes, and thus the robustness of the necessary heteroclinic orbits to sufficiently small symmetric perturbations.

For both networks we consider, the dimension of the unstable manifold of the equilibrium \(B\) is \(2\), and the stable manifolds of \(A\), \(X\), and \(Y\) have dimension \(3\). Moreover, for example, the intersection \(W^{\mathrm{u}}\left(\rmB\right)\cap W^{\mathrm{s}}\left(\rmY\right)\) is a submanifold of dimension \(2\). An infinite number of heteroclinic orbits \(B\to Y\) therefore exist, which also holds for \(B\to X\) in the case of the Kirk--Silber network. In all cases where \(\dim\left(W^{\mathrm{u}}\left(\xi_{j}\right)\cap W^{\mathrm{s}}\left(\xi_{k}\right)\right)\geq 2\), only the single heteroclinic orbit that lies in the coordinate plane containing \(\xi_{j}\) and \(\xi_{k}\) is considered to be part of the cycle and network. We exclude all connections that do not properly lie in a two-dimensional coordinate plane, such as the continuum of orbits \(B\to Y\) in the \(x_{1}=0\) subspace. From Definition~\ref{def:het_cyc}, all heteroclinic cycles we consider are one-dimensional flow-invariant structures. In Definition~\ref{def:het_net}, we allow only a finite number of component heteroclinic cycles, and so the two networks we study are also one-dimensional flow-invariant structures. This definition is in comparison to \cite{kirk_lane_postlethwaite_rucklidge_silber_2010,kirk_postlethwaite_rucklidge_2012}, for example, where a network may contain possibly an infinite number of heteroclinic cycles, such as occurs when one equilibrium has a two-dimensional unstable manifold, and a continuum of heteroclinic orbits exists between two equilibria. Although such a surface of connections exists in the two networks we consider, we only consider one of these connections to be part of the network, as only then can we define the projected map, which is the focus of our study here.

\begin{table}[t!]
  \centering
  \caption{Eigenvalues of \(\rmD f\) evaluated at the four equilibria of the two networks.}\label{tbl:eigs}
  \begin{minipage}{.5\linewidth}
    \centering
    \caption*{(a) Kirk--Silber Network}
    \begin{tabular}{||c|c c c c||}
      \hline
      Eigen-    & \multicolumn{4}{|c||}{Equilibrium} \\
      space     & \(A\)       & \(B\)      & \(X\)       & \(Y\) \\ [0.5ex]
      \hline
      \hline
      \(x_{1}\) & \(-r_{\rmA}\)  & \(-c_{\rmB}\) & \(e_{\rmX\rmA}\)  & \(e_{\rmY\rmA}\)  \\ [0.5ex]
      \(x_{2}\) & \(e_{\rmA}\)   & \(-r_{\rmB}\) & \(-c_{\rmX\rmB}\) & \(-c_{\rmY\rmB}\) \\ [0.5ex]
      \(x_{3}\) & \(-c_{\rmA\rmX}\) & \(e_{\rmB\rmX}\) & \(-r_{\rmX}\)  & \(-c_{\rmY\rmX}\) \\ [0.5ex]
      \(y_{3}\) & \(-c_{\rmA\rmY}\) & \(e_{\rmB\rmY}\) & \(-c_{\rmXY}\) & \(-r_{\rmY}\)  \\ [0.5ex]
      \hline
    \end{tabular}
  \end{minipage}%
  \begin{minipage}{.5\linewidth}
    \caption*{(b) \(\Delta\)-clique Network}
    \begin{tabular}{||c|c c c c||}
      \hline
      Eigen-   & \multicolumn{4}{|c||}{Equilibrium} \\
      space     & \(A\)       & \(B\)      & \(X\)       & \(Y\) \\ [0.5ex]
      \hline
      \hline
      \(x_{1}\) & \(-r_{\rmA}\)  & \(-c_{\rmB}\) & \(-c_{\rmX\rmA}\) & \(e_{\rmY\rmA}\)  \\ [0.5ex]
      \(x_{2}\) & \(e_{\rmA}\)   & \(-r_{\rmB}\) & \(-c_{\rmX\rmB}\) & \(-c_{\rmY\rmB}\) \\ [0.5ex]
      \(x_{3}\) & \(-c_{\rmA\rmX}\) & \(e_{\rmB\rmX}\) & \(-r_{\rmX}\)  & \(-c_{\rmY\rmX}\) \\ [0.5ex]
      \(y_{3}\) & \(-c_{\rmA\rmY}\) & \(e_{\rmB\rmY}\) & \(e_{\rmXY}\)  & \(-r_{\rmY}\)  \\ [0.5ex]
      \hline
    \end{tabular}
  \end{minipage}
\end{table}

We list the eigenvalues of each equilibrium, along with the corresponding coordinate axis that is its eigenspace, in \Tref{tbl:eigs}. Note that each letter in \Tref{tbl:eigs} represents a positive value, and so a negative eigenvalue is indicated with a minus sign. The signs of these eigenvalues ensure the existence of the necessary heteroclinic orbits needed to form each of the two networks we study. The eigenvalues \(-r_{j}\) are known as \textit{radial} eigenvalues in the classification of Krupa and Melbourne \cite{krupa_melbourne_2004}, and the corresponding eigenspace is the radial direction.

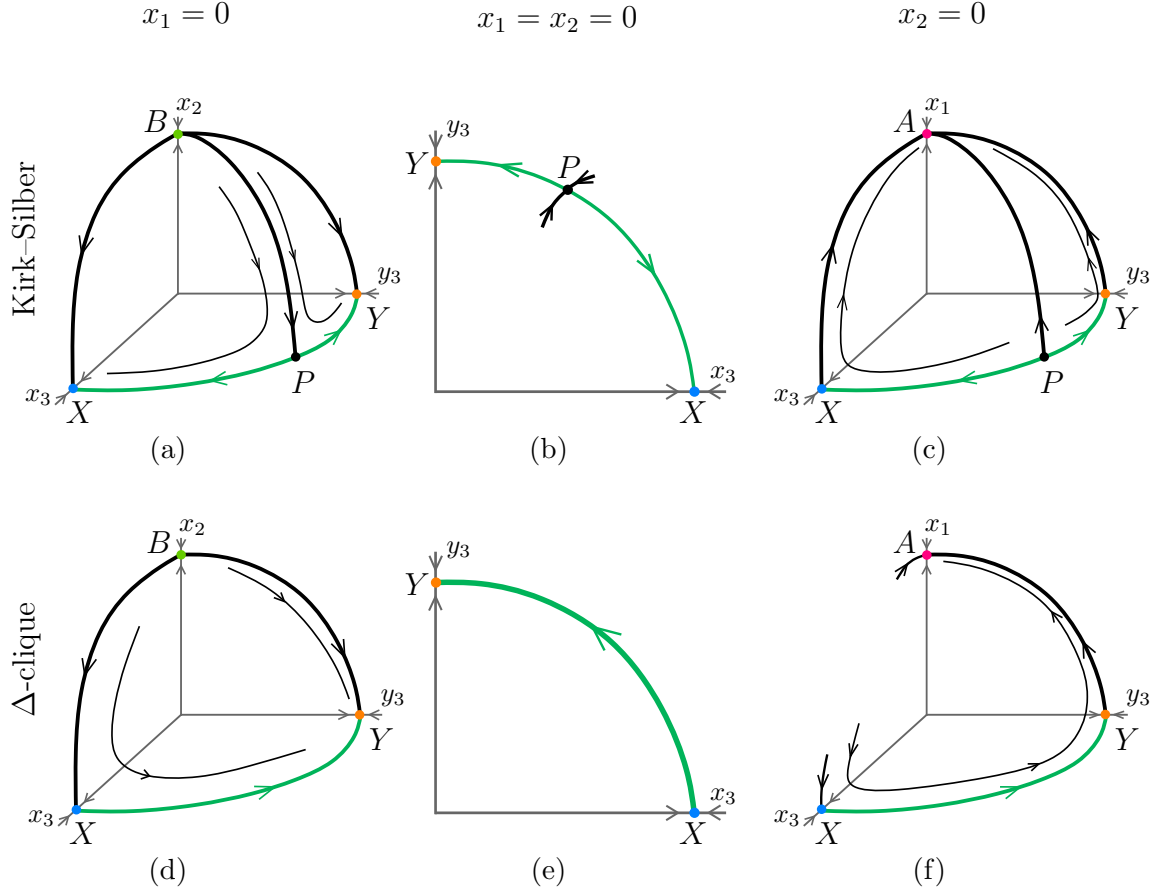
\begin{figure}
  \begin{subfigure}[t]{0.28\linewidth}
    \centering
    \begin{tikzpicture}
      \node at (0,0) () {\includegraphics[width=\linewidth]{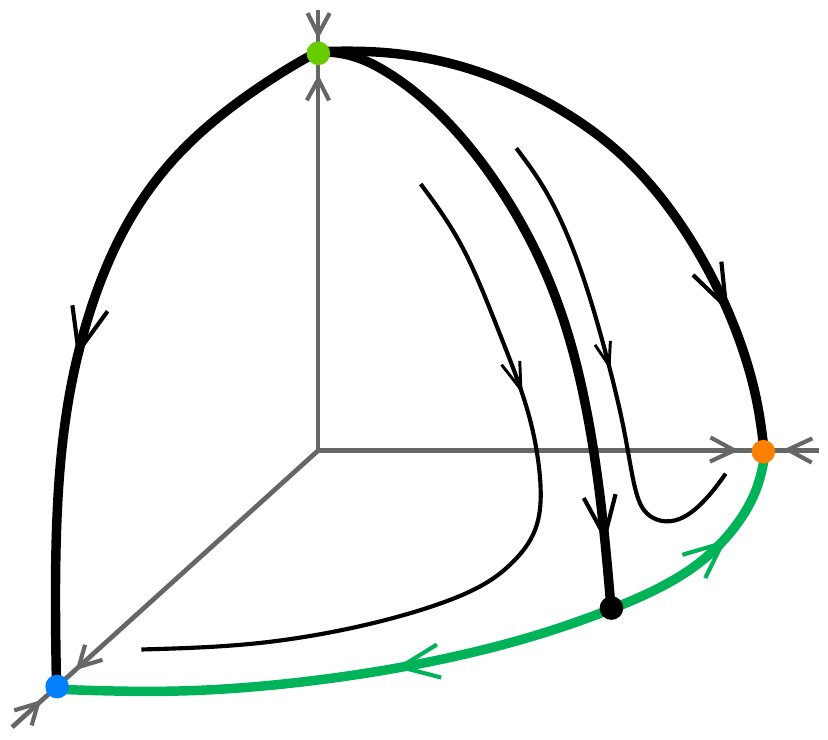}};

      \node at (-0.4,3.25)  () {\(x_{1}=0\)};
      \node[rotate=90] at (-2.55,0.3)  () {Kirk--Silber};

      \node at (-0.8,1.85)  () {\(B\)};
      \node at (1.15,-1.6)  () {\(P\)};
      \node at (-1.8,-2)    () {\(X\)};
      \node at (2.15,-0.75) () {\(Y\)};

      \node at (-0.35,2.05) () {\footnotesize\(x_{2}\)};
      \node at (-2.35,-1.8) () {\footnotesize\(x_{3}\)};
      \node at (2.3,-0.2)  () {\footnotesize\(y_{3}\)};
    \end{tikzpicture}
    \vspace{-8mm}
    \caption{}
    \vspace{5mm}
    \label{fig:dynamics:KS:x2_x3_y3}
  \end{subfigure}%
  \hspace{8mm}%
  \begin{subfigure}[t]{0.26\linewidth}
    \centering
    \begin{tikzpicture}
      \node at (0,0) () {\includegraphics[width=\linewidth]{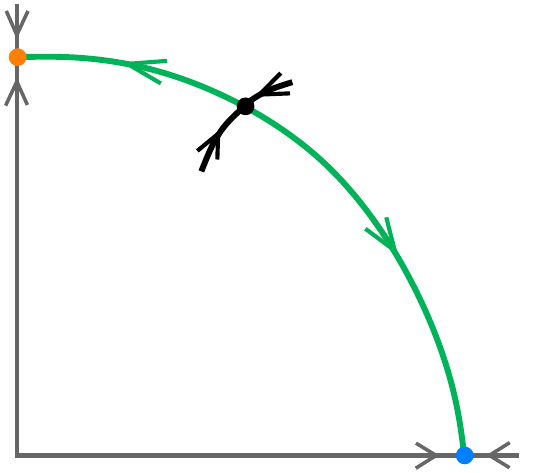}};

      \node at (-0.3,3.25)  () {\(x_{1}=x_{2}=0\)};

      \node at (1.5,-1.95)  () {\(X\)};
      \node at (-2.15,1.35) () {\(Y\)};
      \node at (-0.15,1.3)  () {\(P\)};

      \node at (-1.6,1.8) () {\footnotesize\(y_{3}\)};
      \node at (1.9,-1.45) () {\footnotesize\(x_{3}\)};
    \end{tikzpicture}
    \vspace{-8mm}
    \caption{}
    \vspace{5mm}
    \label{fig:dynamics:KS:C}
  \end{subfigure}%
  \hspace{8mm}%
  \begin{subfigure}[t]{0.28\linewidth}
    \centering
    \begin{tikzpicture}
      \node at (0,0) () {\includegraphics[width=\linewidth]{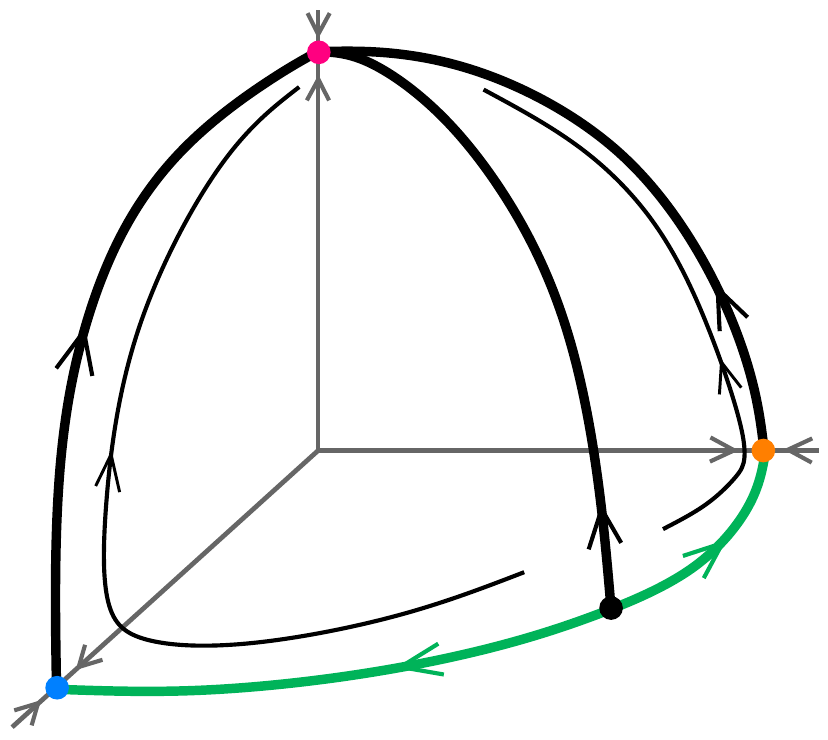}};

      \node at (-0.3,3.25)  () {\(x_{2}=0\)};

      \node at (-0.8,1.85)  () {\(A\)};
      \node at (1.15,-1.6)  () {\(P\)};
      \node at (-1.8,-2)    () {\(X\)};
      \node at (2.15,-0.75) () {\(Y\)};

      \node at (-0.35,2.05) () {\footnotesize\(x_{1}\)};
      \node at (-2.35,-1.8) () {\footnotesize\(x_{3}\)};
      \node at (2.3,-0.2)  () {\footnotesize\(y_{3}\)};
    \end{tikzpicture}
    \vspace{-8mm}
    \caption{}
    \vspace{5mm}
    \label{fig:dynamics:KS:x1_x3_y3}
  \end{subfigure}

  \begin{subfigure}[t]{0.28\linewidth}
    \centering
    \begin{tikzpicture}
      \node at (0,0) () {\includegraphics[width=\linewidth]{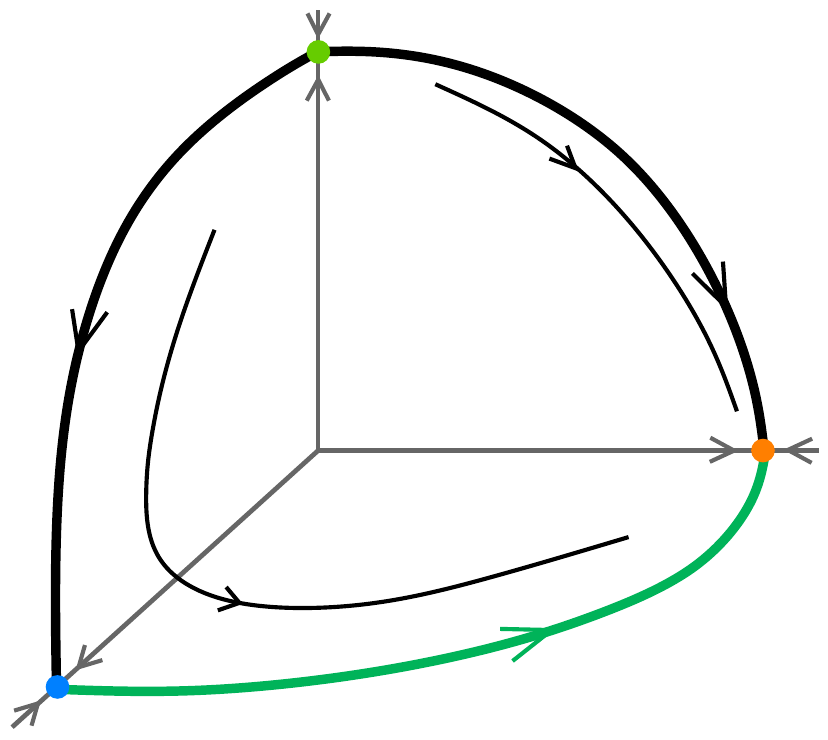}};

      \node[rotate=90] at (-2.55,0.3)  () {\(\Delta\)-clique};

      \node at (-0.8,1.85)  () {\(B\)};
      \node at (-1.8,-2)    () {\(X\)};
      \node at (2.15,-0.75) () {\(Y\)};

      \node at (-0.35,2.05) () {\footnotesize\(x_{2}\)};
      \node at (-2.35,-1.8) () {\footnotesize\(x_{3}\)};
      \node at (2.3,-0.2)  () {\footnotesize\(y_{3}\)};
    \end{tikzpicture}
    \vspace{-8mm}
    \caption{}
    \label{fig:dynamics:D:x2_x3_y3}
  \end{subfigure}%
  \hspace{8mm}%
  \begin{subfigure}[t]{0.26\linewidth}
    \centering
    \begin{tikzpicture}
      \node at (0,0) () {\includegraphics[width=\linewidth]{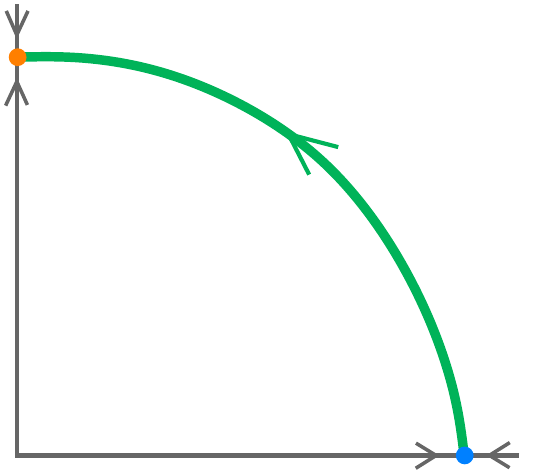}};

      \node at (1.5,-1.95)  () {\(X\)};
      \node at (-2.15,1.35) () {\(Y\)};

      \node at (-1.6,1.8) () {\footnotesize\(y_{3}\)};
      \node at (1.9,-1.45) () {\footnotesize\(x_{3}\)};
    \end{tikzpicture}
    \vspace{-8mm}
    \caption{}
    \label{fig:dynamics:D:C}
  \end{subfigure}%
  \hspace{8mm}%
  \begin{subfigure}[t]{0.28\linewidth}
    \centering
    \begin{tikzpicture}
      \node at (0,0) () {\includegraphics[width=\linewidth]{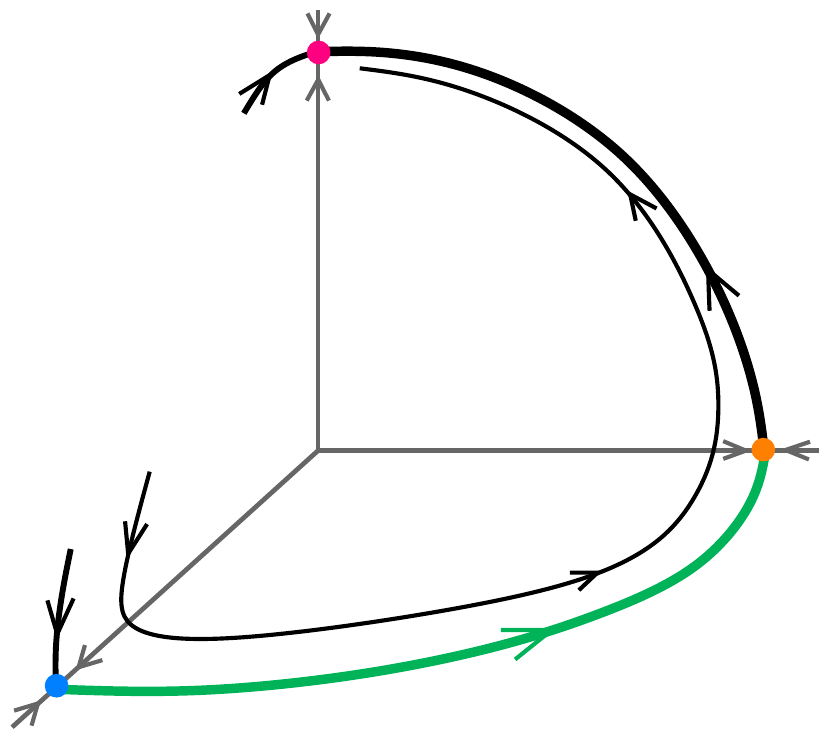}};

      \node at (-0.8,1.85)  () {\(A\)};
      \node at (-1.8,-2)    () {\(X\)};
      \node at (2.15,-0.75) () {\(Y\)};

      \node at (-0.35,2.05) () {\footnotesize\(x_{1}\)};
      \node at (-2.35,-1.8) () {\footnotesize\(x_{3}\)};
      \node at (2.3,-0.2)  () {\footnotesize\(y_{3}\)};
    \end{tikzpicture}
    \vspace{-8mm}
    \caption{}
    \label{fig:dynamics:D:x1_x3_y3}
  \end{subfigure}

  \caption{Dynamics within various subspaces of the two dynamical systems being considered. Equilibria are given by solid dots. The equilibria that are components of the network are coloured as in \Fref{fig:R4_networks}, while the equilibrium \(P\) of the Kirk--Silber network is black. Heteroclinic orbits are thick lines. The orbits that are part of the flow-invariant curve \(C\) are coloured green, and all other orbits are coloured black. Example trajectories are thin black lines, and the flow-invariant coordinate axes are thin grey lines. (a) Kirk--Silber and (d) \(\Delta\)-clique: The \(x_{1}=0\) subspace, showing the two dimensional unstable manifold of \(B\) and the curve \(C\), and how these structures interact. (b) Kirk--Silber and (e) \(\Delta\)-clique: The \(\left(x_{3},y_{3}\right)\)-plane and details of the invariant curve \(C\). (c) Kirk--Silber and (f) \(\Delta\)-clique: The \(x_{2}=0\) subspace, showing two dimensional stable manifold of \(A\) and how it interacts with \(C\).}
  \label{fig:dynamics}
\end{figure}

We present diagrammatic representations of the dynamics in relevant subspaces of \(\R^{4}\) in \Fref{fig:dynamics}.

In the case of the Kirk--Silber network, we assume that there is an additional hyperbolic saddle equilibrium \(P\) in the \(\left(x_{3},y_{3}\right)\)-plane, and that there are robust heteroclinic orbits \(P\to X\) and \(P\to Y\) in this plane, and \(B\to P\) in the \(x_{1}=0\) subspace. See Figures~\ref{fig:dynamics}(\subref{fig:dynamics:KS:x2_x3_y3}) and \ref{fig:dynamics}(\subref{fig:dynamics:KS:C}). The equilibrium \(P\) and these additional orbits are themselves not a component of the Kirk--Silber network, but their existence is needed for the analysis presented in \sref{sec:ret_map}.

For both networks, we consider the dynamics near a flow-invariant curve that exists in the \(\left(x_{3},y_{3}\right)\)-plane. In both cases, we label this curve \(C\), and we colour this curve green in \Fref{fig:dynamics}. For the Kirk--Silber network, the curve is the union of the three equilibria \(X\), \(Y\), and \(P\), and the two orbits \(P\to X\) and \(P\to Y\). For the \(\Delta\)-clique network, the curve \(C\) is the union of the two equilibria \(X\) and \(Y\), and the orbit \(X\to Y\). We assume that the curve \(C\) can be parametrised by an angle \(\theta_{3}\), defined by \(\tan\theta_{3}=y_{3}/x_{3}\).

For the \(\Delta\)-clique network, we make no particular assumptions about the dynamics in the \(\left(x_{1},x_{3}\right)\)-plane, other than the sign of the eigenvalue \(-c_{\rmA\rmX}\) at \(A\) in the \(x_{3}\)-direction, and the sign of the eigenvalue \(-c_{\rmX\rmA}\) at \(X\) in the \(x_{1}\)-direction. For this reason, we have shown in \Fref{fig:dynamics}(\subref{fig:dynamics:D:x1_x3_y3}) only a portion of the dynamics in this \(\left(x_{1},x_{3}\right)\)-plane, locally near \(A\) and \(X\).

\subsection{Return maps and transition matrices}\label{sec:problem:ssec:ret_maps}

The usual method to analyse the dynamics near heteroclinic networks is to construct return maps to cross-sections defined transverse to the flow near the network. There are many examples of the process in this literature: see any of \cite{melbourne_1991,krupa_melbourne_1995a,kirk_silber_1994,chossat_krupa_melbourne_scheel_1997,krupa_melbourne_2004,postlethwaite_dawes_2006,postlethwaite_dawes_2010,postlethwaite_2010,podvigina_ashwin_2011,podvigina_2012,castro_lohse_2016b,garrido_da_silva_castro_2019,postlethwaite_rucklidge_2022,article}.

Typically, return maps are defined on cross-sections defined in a small neighbourhood of each equilibrium of the network, and transverse to either the incoming or outgoing heteroclinic orbit. If the equilibrium is a part of more than one cycle of the network, the return map is constructed in a piecewise manner, with components representing each cycle.

When studying, for example, a periodic orbit, only a single cross-section transverse to the flow of the orbit needs to be constructed, and the return map defined on this section is qualitatively the same, regardless of where the cross-section is defined. Analysing a heteroclinic cycle, however, presents additional complications. First, there may be an open set of points on one cross-section that gives rise to trajectories that do not return to that cross-section, while, for a different cross-section, all trajectories may return to that cross-section. Second, trajectories may be asymptotic to the intersection of the heteroclinic cycle and a cross-section at one equilibrium, but not at other equilibria. See \cite{kirk_silber_1994,postlethwaite_rucklidge_2022} for examples. As such, a return map needs to be constructed and analysed for each equilibrium in the cycle.

For the Kirk--Silber and \(\Delta\)-clique networks, trajectories at \(B\) sufficiently close to the network leave either in the direction of \(X\) or \(Y\), and so cycle around one of the two component cycles of the network. As such, we refer to the equilibrium \(B\) as a \textit{splitting equilibrium} of the network. For these two networks, \(B\) is the only splitting equilibrium.

The return map to \(B\) is defined on
\begin{equation*}
  \poinin[A]{B}=\left\{(x_{1},x_{2},x_{3},y_{3})\in\R^{4}\mid x_{1}=h,|x_{2}-1|<h,0<x_{3},y_{3}<h\right\},
\end{equation*}
which is defined in a small neighbourhood of \(B\), transverse to the heteroclinic orbit \(A\to B\). The constant \(h\ll 1\) is chosen to be sufficiently small that the cross-section is contained within a neighbourhood where the flow is approximately linear. In the definition of \(\poinin[A]{B}\), \(x_{1}\) is fixed, and, as a consequence of the invariant sphere theorem \cite{field_1996}, we can ignore the radial coordinate, which near \(B\) is \(x_{2}\). As such, the relevant coordinates of \(\poinin[A]{B}\) are \(\left(x_{3},y_{3}\right)\).

For each network, we construct a return map \(\Phi\colon\poinin[A]{B}\to\poinin[A]{B}\) by linearising the flow near each equilibrium and along heteroclinic orbits. Trajectories asymptotic to the network are represented in the return map by orbits converging to the origin in \(\poinin[A]{B}\). However, not all trajectories starting on \(\poinin[A]{B}\) remain sufficiently close to the heteroclinic network for these linearisations to remain a valid approximation, and so the return map \(\Phi\) is not, in fact, defined on all of \(\poinin[A]{B}\).

We define the curve
\begin{equation*}
  \Sigma_{\mrm{s}}=\left\{\left(x_{3},y_{3}\right)\in\poinin[A]{B}\mid x_{3}^{e_{\rmB\rmY}}=y_{3}^{e_{\rmB\rmX}}\right\}.
\end{equation*}
This curve, which we refer to as the \textit{switching curve}, partitions \(\poinin[A]{B}\) into two sets of initial conditions, for which either \(x_{3}\) or \(y_{3}\) becomes \(\order{1}\) first under the local flow near \(B\). Trajectories starting at initial conditions sufficiently close to this curve do not leave a small neighbourhood of \(B\) in the direction of \(X\) or \(Y\); that is, after leaving a small neighbourhood of \(B\), a trajectory enters a small neighbourhood of \(C\), but does not enter this neighbourhood within distance \(h\) from either \(X\) or \(Y\). The trajectory may, however, leave a small neighbourhood of \(C\) near \(X\) or \(Y\).
In the case of the Kirk--Silber network, such a trajectory may enter a small neighbourhood of the curve \(C\) near the equilibrium \(P\), or somewhere along the heteroclinic orbits \(P\to X\) or \(P\to Y\), but not close to \(X\) or \(Y\).  In the case of the \(\Delta\)-clique network, a trajectory starting at a point in \(\Gamma_{\rmc}\) enters a small neighbourhood of \(C\) along the heteroclinic orbit \(X\to Y\), but not close to \(X\) or \(Y\).

These trajectories do not remain sufficiently close to the network for the linearisations that are used to construct the return map to remain valid. As such, we must exclude from the domain of definition of \(\Phi\) those points that are too close to the curve \(\Sigma_{\rms}\). The set of points that we exclude is:
\begin{equation*}
  \Gamma_{\rmc}=\left\{\left(x_{1},x_{2},x_{3},y_{3}\right)\in\poinin[A]{B}\mid x_{1}=h,|x_{2}-1|<h,\left(1-\epsilon\right)x_{3}^{\frac{e_{\rmB\rmY}}{e_{\rmB\rmX}}}<y_{3}<\left(\frac{x_{3}}{1-\epsilon}\right)^{\frac{e_{\rmB\rmY}}{e_{\rmB\rmX}}}\right\}
\end{equation*}
where \(\epsilon\) is some positive real number. We refer to \(\Gamma_{\rmc}\) as the \textit{excluded cusp}, as it is the set of initial conditions excluded from the domain of \(\Phi\). Trajectories that begin in \(\Gamma_{c}\) do not remain sufficiently close to the network for the linearisations used to construct \(\Phi\) to remain valid. In particular, after a trajectory starting at a point in \(\Gamma_{c}\) leaves a small neighbourhood of the equilibrium \(B\), it does not strike cross-sections defined in a small neighbourhood of \(X\) or \(Y\) that are transverse to the heteroclinic orbit \(B\to X\) or \(B\to Y\), respectively.

We next define the two sets
\begin{equation}\label{eqn:gamma_sets}
  \begin{aligned}
    \Gamma_{\rmX}=&\left\{\left(x_{1},x_{2},x_{3},y_{3}\right)\in\poinin[A]{B}\mid x_{1}=h,|x_{2}-1|<h,\left(1-\epsilon\right)x_{3}^{\frac{e_{\rmB\rmY}}{e_{\rmB\rmX}}}\geq y_{3}\right\},\ \textrm{ and}\\
    \Gamma_{\rmY}=&\left\{\left(x_{1},x_{2},x_{3},y_{3}\right)\in\poinin[A]{B}\mid x_{1}=h,|x_{2}-1|<h,\left(1-\epsilon\right)y_{3}^{\frac{e_{\rmB\rmX}}{e_{\rmB\rmY}}}\geq x_{3}\right\}.\\
  \end{aligned}
\end{equation}
Trajectories that begin in \(\Gamma_{\rmX}\) leave a small neighbourhood of \(B\) in the direction of the equilibrium \(X\), and enter a small neighbourhood of \(C\) near \(X\), striking a cross-section that is transverse to the orbit \(B\to X\), and similarly for \(\Gamma_{\rmY}\).

\begin{figure}
  \centering
  \includegraphics[width=0.4\linewidth]{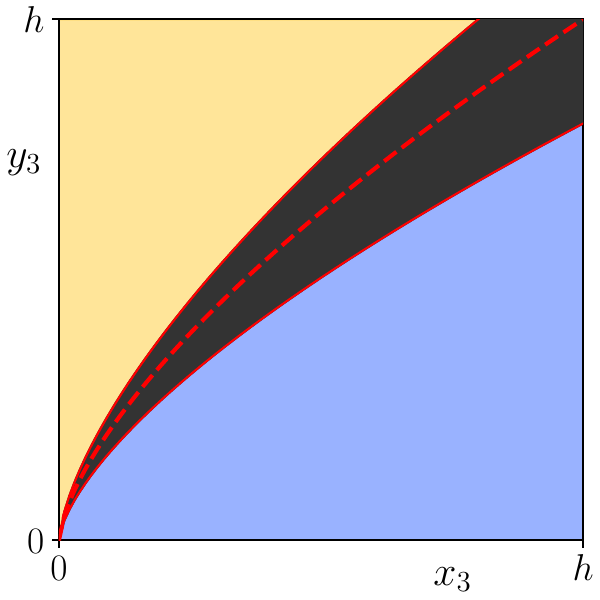}
  \caption{A representation of the cross-section \(\poinin[A]{B}\). Points in the set \(\Gamma_{\rmX}\) are shaded blue, and those in the set \(\Gamma_{\rmY}\) are shaded orange. The excluded cusp \(\Gamma_{\rmc}\) is shaded black. The curve defined by \(x_{3}^{e_{\rmB\rmY}}=y_{3}^{e_{\rmB\rmX}}\) is a dashed red line, while the solid red lines are the boundaries between \(\Gamma_{\rmX}\), \(\Gamma_{\rmY}\), and \(\Gamma_{\rmc}\).}
  \label{fig:section}
\end{figure}

The three subsets \(\Gamma_{\rmX}\), \(\Gamma_{\rmY}\), and \(\Gamma_{\rmc}\) form a partition of \(\poinin[A]{B}\), and \Fref{fig:section} gives a representation of these three sets.

Following methodology similar to that in \cite{kirk_silber_1994}, we can derive a return map for each network, \(\Phi_{\mathrm{KS}},\Phi_{\Delta}\colon\Gamma_{\rmX}\cup\Gamma_{\rmY}\to\poinin[A]{B}\), which is defined in a piecewise manner, and is given to leading order by
\begin{equation}\label{eqn:ks_ret_map}
  \Phi_{\mathrm{KS}}\left(x_{3},y_{3}\right)=\begin{cases}
    \Phi_{\rmX}\left(x_{3},y_{3}\right)=\left(C_{1}x_{3}^{\delta_{\rmX}},C_{2}x_{3}^{\rho_{\rmX}}y_{3}\right)& \textrm{if}\ \left(x_{3},y_{3}\right)\in\Gamma_{\rmX},\\
    \Phi_{\rmY}\left(x_{3},y_{3}\right)=\left(D_{1}x_{3}y_{3}^{\rho_{\rmY}},D_{2}y_{3}^{\delta_{\rmY}}\right)& \textrm{if}\ \left(x_{3},y_{3}\right)\in\Gamma_{\rmY},
  \end{cases}
\end{equation}
and
\begin{equation}\label{eqn:delta_ret_map}
  \Phi_{\Delta}\left(x_{3},y_{3}\right)=\begin{cases}
    \Phi_{\rmXY}\left(x_{3},y_{3}\right)=\left(E_{1}x_{3}^{\alpha_{1}}y_{3}^{\alpha_{2}},E_{2}x_{3}^{\alpha_{3}}y_{3}^{\alpha_{4}}\right)& \textrm{if}\ \left(x_{3},y_{3}\right)\in\Gamma_{\rmX},\\
    \Phi_{\rmY}\left(x_{3},y_{3}\right)=\left(F_{1}x_{3}y_{3}^{\rho_{\rmY}},F_{2}y_{3}^{\delta_{\rmY}}\right)& \textrm{if}\ \left(x_{3},y_{3}\right)\in\Gamma_{\rmY},
  \end{cases}
\end{equation}
where the terms \(C_{j},D_{j},E_{j},F_{j}\) are \(\order{1}\) positive constants, and, for the maps \(\Phi_{\rmX}\) and \(\Phi_{\rmY}\):
\begin{align*}
  \delta_{\rmX}
  &=\frac{c_{\rmA\rmX}c_{\rmB}c_{\rmX\rmB}}{e_{\rmA}e_{\rmB\rmX}e_{\rmX\rmA}},
  \qquad
  \rho_{\rmX}
  =-\frac{e_{\rmB\rmY}}{e_{\rmB\rmX}}+\frac{c_{\rmB}c_{\rmXY}}{e_{\rmB\rmX}e_{\rmX\rmA}}+\frac{c_{\rmA\rmY}c_{\rmB}c_{\rmX\rmB}}{e_{\rmA}e_{\rmB\rmX}e_{\rmX\rmA}},
  \\[0.4em]
  \delta_{\rmY}
  &=\frac{c_{\rmA\rmY}c_{\rmB}c_{\rmY\rmB}}{e_{\rmA}e_{\rmB\rmY}e_{\rmY\rmA}},
  \qquad
  \rho_{\rmY}=
  -\frac{e_{\rmB\rmX}}{e_{\rmB\rmY}}+\frac{c_{\rmB}c_{\rmY\rmX}}{e_{\rmB\rmY}e_{\rmY\rmA}}+\frac{c_{\rmA\rmX}c_{\rmB}c_{\rmY\rmB}}{e_{\rmA}e_{\rmB\rmY}e_{\rmY\rmA}};
\end{align*}
and, for the map \(\Phi_{\rmXY}\):
\begin{align*}
  \alpha_{1}
  &=\frac{c_{\rmB}c_{\rmY\rmX}}{e_{\rmB\rmX}e_{\rmY\rmA}} + \frac{c_{\rmB}c_{\rmY\rmB}c_{\rmA\rmX}}{e_{\rmB\rmX}e_{\rmY\rmA}e_{\rmB}} - \frac{c_{\rmX\rmB}e_{\rmB\rmY}c_{\rmA\rmX}}{e_{\rmXY}e_{\rmB\rmX}e_{\rmB}} - \frac{c_{\rmY\rmX}e_{\rmB\rmY}c_{\rmX\rmA}}{e_{\rmY\rmA}e_{\rmB\rmX}e_{\rmXY}} - \frac{c_{\rmY\rmB}e_{\rmB\rmY}c_{\rmA\rmX}c_{\rmX\rmA}}{e_{\rmY\rmA}e_{\rmB\rmX}e_{\rmB}e_{\rmXY}},
  \\[0.4em]
  \alpha_{2}
  &=\frac{c_{\rmX\rmB}c_{\rmA\rmX}}{e_{\rmXY}e_{\rmB}} + \frac{c_{\rmY\rmX}c_{\rmX\rmA}}{e_{\rmY\rmA}e_{\rmXY}} + \frac{c_{\rmY\rmB}c_{\rmA\rmX}c_{\rmX\rmA}}{e_{\rmY\rmA}e_{\rmB}e_{\rmXY}},
  \\[0.4em]
  \alpha_{3}
  &=\frac{c_{\rmA\rmY}c_{\rmB}c_{\rmY\rmB}}{e_{\rmB}e_{\rmB\rmX}e_{\rmY\rmA}} - \frac{c_{\rmA\rmY}c_{\rmX\rmB}e_{\rmB\rmY}}{e_{\rmB}e_{\rmXY}e_{\rmB\rmX}} - \frac{c_{\rmA\rmY}c_{\rmY\rmB}e_{\rmB\rmY}c_{\rmX\rmA}}{e_{\rmB}e_{\rmY\rmA}e_{\rmB\rmX}e_{\rmXY}},
  \\[0.4em]
  \alpha_{4}
  &=\frac{c_{\rmA\rmY}c_{\rmX\rmB}}{e_{\rmB}e_{\rmXY}} + \frac{c_{\rmA\rmY}c_{\rmY\rmB}c_{\rmX\rmA}}{e_{\rmB}e_{\rmY\rmA}e_{\rmXY}}.
\end{align*}
The map \(\Phi_{\rmY}\) is the same for both networks, as the cycle \(\cyc{Y}\) is common to both. The map \(\Phi_{\rmXY}\) is so labelled as it corresponds to trajectories that visit \(X\) and then \(Y\) before returning to \(B\), whereas \(\Phi_{\rmX}\) and \(\Phi_{\rmY}\) describe trajectories that visit only one of \(X\) or \(Y\) before returning to \(B\).

Podvigina \cite{podvigina_2012} and Garrido da Silva and Castro \cite{garrido_da_silva_castro_2019} have shown that the general form of return maps near heteroclinic cycles to a \(p\)-dimensional cross-section is
\begin{equation}\label{eqn:gen_ret_map}
  \Phi\left(x_{1},\dots,x_{p}\right)=\left(C_{1}x_{1}^{\alpha_{11}}\dots x_{p}^{\alpha_{1p}},\dots,C_{p}x_{1}^{\alpha_{p1}}\dots x_{p}^{\alpha_{pp}}\right),
\end{equation}
where the \(C_{j}\) are \(\order{1}\) positive constants. From return maps of this form, we can compute \textit{transition matrices}, defined from the exponents of the components of the return map:
\begin{equation*}
  M\left(\Phi\right)=\begin{pmatrix}
    \alpha_{11} & \hdots & \alpha_{1p} \\
    \vdots & \ddots & \vdots \\
    \alpha_{p1} & \hdots & \alpha_{pp}
  \end{pmatrix}.
\end{equation*}
These transition matrices act on the coordinates \(X_{j}\coloneqq\log x_{j}\). Since \(x_{j}<h\ll 1\), the logarithmic coordinates \(X_{j}\ll 0\). For this reason, when transforming \(\Phi\) to logarithmic coordinates, the logarithm of the \(\order{1}\) constants \(C_{j}\) can be ignored as they do not appear in asymptotically significant expressions. A trajectory asymptotic to the network is represented in logarithmic coordinates by both components of the resulting vector \(\left(X_{3},Y_{3}\right)\) diverging to \(-\infty\). The return maps in \eqref{eqn:ks_ret_map} and \eqref{eqn:delta_ret_map} are represented by a piecewise-linear map \(M\), whose components are the transition matrices \(M_{\rmX}\coloneqq M(\Phi_{\rmX})\), \(M_{\rmY}\coloneqq M(\Phi_{\rmY})\), or \(M_{\rmXY}\coloneqq M(\Phi_{\rmXY})\).

In \cite[Theorem~5]{podvigina_2012}, Podvigina provides necessary and sufficient conditions for the stability of a heteroclinic cycle, based on properties of the eigenvalues and eigenvectors of the transition matrices of that cycle. These properties generically need to be checked for the transition matrix of each return map of the cycle.

\subsection{The projected map}

In \cite{article}, we presented an extension of this methodology, called the \textit{projected map}. The main idea behind this approach is to simplify the analysis of dynamics near heteroclinic networks by identifying all trajectories that have qualitatively the same behaviour. We summarise the results of \cite{article} in the remainder of this section.

The projected map is constructed by first taking a slice across the negative quadrant of the \(\left(X_{3},Y_{3}\right)\)-plane, defined by
\begin{equation}\label{eqn:S}
  S=\left\{\left(X_{3},Y_{3}\right)\in\negativeR{2}\mid X_{3}+Y_{3}=-1\right\}.
\end{equation}
We then identify each linear subspace of \(\R^{2}\)---restricted to the negative quadrant---by where it intersects with the set \(S\). Since \(S\) is one-dimensional, we use the \(X_{3}\) coordinate of a point on \(S\) as its coordinate.  Therefore, we identify \(S\) with the open interval \(\left(-1,0\right)\). The projection of any point \(\left(X_{3},Y_{3}\right)\in\negativeR{2}\) onto \(S\) is given by
\begin{equation}\label{eqn:proj}
  \Pi\left(X_{3},Y_{3}\right)=\frac{-X_{3}}{X_{3}+Y_{3}}.
\end{equation}

The curve \(\Sigma_{\rms}\), around which the excluded cusp \(\Gamma_{\rmc}\) is defined, is represented by a point \(\vartheta_{\rms}\in S\):
\begin{equation*}
  \vartheta_{\rms}=\dfrac{-1}{1+\frac{e_{\rmB\rmY}}{e_{\rmB\rmX}}}.
\end{equation*}
The set \(\Gamma_{\rmX}\) is represented by the open interval \(\Theta_{\rmX}=\left(\vartheta_{\rms},0\right)\subseteq S\), and \(\Gamma_{\rmY}\) by the open interval \(\Theta_{\rmY}=\left(-1,\vartheta_{\rms}\right)\subseteq S\), and so \(S=\Theta_{\rmX}\cup\left\{\vartheta_{\rms}\right\}\cup\Theta_{\rmY}\).

\begin{figure}
  \centering
  \begin{subfigure}[t]{0.5\linewidth}
    \centering
    \includegraphics[width=\linewidth]{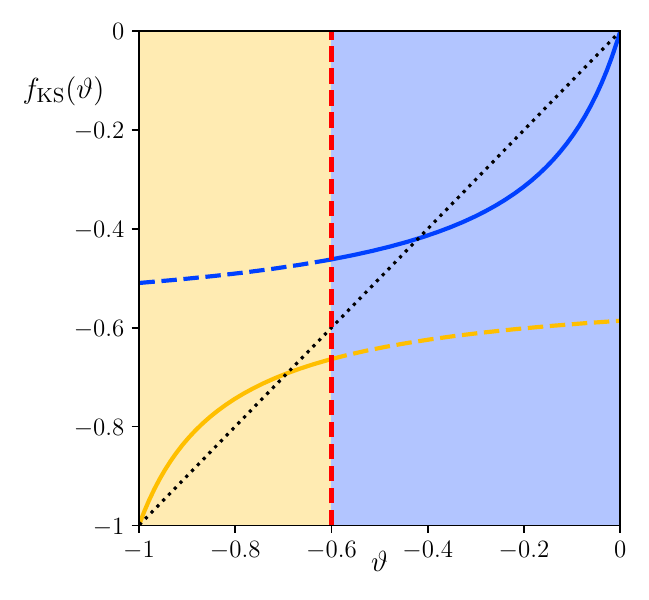}
    \caption{\(f_{\mathrm{KS}}\)}
    \label{fig:proj_maps:kirk_silber_net}
  \end{subfigure}%
  \hfill%
  \begin{subfigure}[t]{0.5\linewidth}
    \centering
    \includegraphics[width=\linewidth]{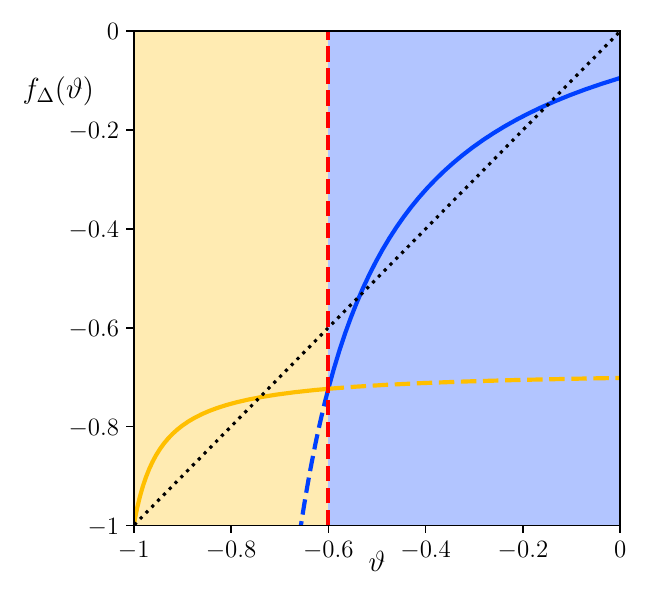}
    \caption{\(f_{\Delta}\)}
    \label{fig:proj_maps:delta_clique_net}
  \end{subfigure}%
  \caption{The projected map of each network, with parameter values chosen so that both cycles of each network are stable. The yellow shaded region is \(\Theta_{\rmY}\), and the blue shaded region is \(\Theta_{\rmX}\). Each component function is shown as solid within its domain of definition, and dashed without. The dashed red line corresponds to \(\vartheta_{\rms}\), and the line \(\vartheta=f\left(\vartheta\right)\) is the dotted black line. The projected map of the Kirk--Silber network can be seen to be discontinuous on \(\vartheta_{\rms}\), while that of the \(\Delta\)-clique network is continuous. For both networks, fixed points of the map, where the graph intersects the diagonal, correspond to the respective heteroclinic cycle, and it is necessary for the fixed point to lie in the correct subdomain of \(f\) for the corresponding cycle to be fragmentarily asymptotically stable. Parameter values for these Figures~can be found in \ref{app:parms}.}
  \label{fig:proj_maps}
\end{figure}

A projected map can then be defined for each network, \(f_{\mathrm{KS}}\colon\Theta_{\rmX}\cup\Theta_{\rmY}\to S\) and \(f_{\Delta}\colon\Theta_{\rmX}\cup\Theta_{\rmY}\to S\). These maps are defined as the induced action of the transition matrices on \(S\), and are piecewise-smooth maps. The point \(\vartheta_{\rms}\), which divides its two components, is known as the \textit{switching manifold}.

The projected map of the Kirk--Silber network is
\begin{equation}\label{eqn:ks_proj_map}
  f_{\mathrm{KS}}\left(\vartheta\right)=\begin{cases}
    f_{\rmX}\left(\vartheta\right)=\dfrac{-\delta_{\rmX}\vartheta}{\left(\delta_{\rmX}+\rho_{\rmX}-1\right)\vartheta-1} & \textrm{if}\ \vartheta\in\Theta_{\rmX},\\[1em]
    f_{\rmY}\left(\vartheta\right)=\dfrac{\left(1-\rho_{\rmY}\right)\vartheta-\rho_{\rmY}}{\left(\delta_{\rmY}+\rho_{\rmY}-1\right)\vartheta+\delta_{\rmY}+\rho_{\rmY}} & \textrm{if}\ \vartheta\in\Theta_{\rmY},
  \end{cases}
\end{equation}
and the projected map of the \(\Delta\)-clique network is
\begin{equation}\label{eqn:delta_proj_map}
  f_{\Delta}\left(\vartheta\right)=\begin{cases}
    f_{\rmXY}\left(\vartheta\right)=\dfrac{\left(\alpha_{2}-\alpha_{1}\right)\vartheta+\alpha_{2}}{\left(\alpha_{1}-\alpha_{2}+\alpha_{3}-\alpha_{4}\right)\vartheta} & \textrm{if}\ \vartheta\in\Theta_{\rmX},\\[1em]
    f_{\rmY}\left(\vartheta\right)=\dfrac{\left(1-\rho_{\rmY}\right)\vartheta-\rho_{\rmY}}{\left(\delta_{\rmY}+\rho_{\rmY}-1\right)\vartheta+\delta_{\rmY}+\rho_{\rmY}} & \textrm{if}\ \vartheta\in\Theta_{\rmY}.
  \end{cases}
\end{equation}
\Fref{fig:proj_maps} shows examples of these projected maps.

The projected map simplifies the analysis of the dynamics near these two networks to the study of the dynamics of a piecewise-smooth, one-dimensional discrete map. Heteroclinic cycles that are fragmentarily asymptotically stable in the dynamical system are represented in the projected map by asymptotically stable fixed points that lie in the correct subdomain of the map; that is, for example, it is necessary for the fixed point of the map \(f_{\rmX}\) to lie in \(\Theta_{\rmX}\). In the theory of piecewise-smooth dynamical systems, such a fixed point is called \textit{admissible} \cite{simpson_2016}. A trajectory asymptotic to a cycle is represented in the projected map by an orbit asymptotic to these fixed points. Trajectories that make a full excursion around one and then make a full excursion around the other cycle---that is, trajectories that switch between the two cycles---are represented by a point in one of \(\Theta_{\rmX}\) or \(\Theta_{\rmY}\) that is mapped into the other domain by the projected map.

We prove in \cite{article} that, for the Kirk--Silber and \(\Delta\)-clique networks, all stability conditions of component cycles are represented by conditions on the stability and admissibility of the fixed points of \(f\). Therefore, we can study the dynamics near these networks solely with the projected map.

\subsection{Continuity of projected maps}\label{sec:problem:ssec:cont}

In \cite{article}, we observed that the projected map of the Kirk--Silber network is generically discontinuous on its switching manifold, whereas the projected map of the \(\Delta\)-clique network is continuous on its switching manifold. We formalise this difference as follows.

In the case of the Kirk--Silber network, we find
\begin{equation}\label{eqn:ks_lim_right}
  \lim_{\vartheta\searrow\vartheta_{\rms}}f_{\mathrm{KS}}(\vartheta)=\lim_{\vartheta\searrow\vartheta_{\rms}}f_{\rmX}(\vartheta)=\frac{-c_{\rmA\rmX}c_{\rmX\rmB}}{(c_{\rmA\rmX}+c_{\rmA\rmY})c_{\rmX\rmB}+e_{\rmA}c_{\rmXY}},
\end{equation}
and
\begin{equation}\label{eqn:ks_lim_left}
  \lim_{\vartheta\nearrow\vartheta_{\rms}}f_{\mathrm{KS}}(\vartheta)=\lim_{\vartheta\nearrow\vartheta_{\rms}}f_{\rmY}(\vartheta)=\frac{-\left(c_{\rmA\rmX}c_{\rmY\rmB}+e_{\rmA}c_{\rmY\rmX}\right)}{(c_{\rmA\rmX}+c_{\rmA\rmY})c_{\rmY\rmB}+e_{\rmA}c_{\rmY\rmX}}.
\end{equation}
In the case of the \(\Delta\)-clique network, we find
\begin{equation}\label{eqn:delta_lim}
  \lim_{\vartheta\to\vartheta_{\rms}}f_{\Delta}(\vartheta)=\lim_{\vartheta\searrow\vartheta_{\rms}}f_{\rmXY}(\vartheta)=\lim_{\vartheta\nearrow\vartheta_{\rms}}f_{\rmY}(\vartheta)=\frac{-\left(c_{\rmA\rmX}c_{\rmY\rmB}+e_{\rmA}c_{\rmY\rmX}\right)}{(c_{\rmA\rmX}+c_{\rmA\rmY})c_{\rmY\rmB}+e_{\rmA}c_{\rmY\rmX}}.
\end{equation}
Therefore, we consider the projected map of the Kirk--Silber network to be discontinuous, and the projected map of the \(\Delta\)-clique network to be continuous, which can be seen in \Fref{fig:proj_maps}. The discontinuity in the case of the Kirk--Silber network is perhaps unexpected, since the projected map is ultimately derived from a smooth flow through a process defined only with continuous functions.

We emphasise that the projected maps are not defined at \(\vartheta_{\rms}\), due to the mechanism by which the map was constructed, excluding \(\Gamma_{\rmc}\) from the domain of \(\Phi\). However, we prove in \cite{article} that only the point \(\vartheta_{\rms}\) needs to be excluded from the domain of definition of the projected map, provided the trajectories are sufficiently close to the network.

We highlight that the limit in \eqref{eqn:delta_lim} is the same as in \eqref{eqn:ks_lim_left}. Both of these limits are calculated from the part of the projected map that describes the dynamics of trajectories that leave a neighbourhood of the curve \(C\) near the equilibrium \(Y\), with \(\theta_{3}\) close to \(\frac{\pi}{2}\). This similarity provides us with an initial indication of the dynamical properties of the system that produce the discontinuity in the case of the Kirk--Silber network. In particular, all trajectories sufficiently close to the \(\Delta\)-clique network leave a neighbourhood of the invariant curve \(C\) near the equilibrium \(Y\), even if they have first visited \(X\) and so are in the domain of \(f_{\rmXY}\). In contrast, trajectories sufficiently close to the Kirk--Silber network can leave a small neighbourhood of \(C\) near either \(X\) or \(Y\). The left-hand limit of the Kirk--Silber network---which corresponds to trajectories that leave near \(Y\)---is equal to the limit of the \(\Delta\)-clique network. However, the right-hand limit differs, and corresponds to trajectories that leave near \(X\). We can see in \Fref{fig:dynamics}(\subref{fig:dynamics:KS:x2_x3_y3}) that there is a separatrix that divides trajectories near the Kirk--Silber network between those that pass close to \(X\) or \(Y\). In \Fref{fig:dynamics}(\subref{fig:dynamics:D:x2_x3_y3}), we see for the \(\Delta\)-clique network that there is a surface of orbits from \(B\) to \(Y\), with some passing close to \(X\).

In Sections \ref{sec:ret_map} and \ref{sec:anal}, we explain the dynamical origin of the continuity or discontinuity of these projected maps.

\subsection{Evaluation of the limit \(\vartheta\to\vartheta_{\rms}\)}\label{sec:problem:ssec:limit}

An important aspect of our explanation in the remainder of this paper is that the discontinuity emerges, in part, because of how the projected map is defined. We now recall the relevant parts of \cite{article}.

The boundary of the excluded cusp, \(\Gamma_{c}\) is defined by the two curves
\begin{equation*}
  \Sigma_{\rms}^{+}=\left\{\left(x_{3},y_{3}\right)\in\poinin[A]{B}\mid y_{3}=\left(1-\epsilon\right)x_{3}^{\frac{e_{\rmB\rmY}}{e_{\rmB\rmX}}}\right\}
\end{equation*}
and
\begin{equation*}
  \Sigma_{\rms}^{-}=\left\{\left(x_{3},y_{3}\right)\in\poinin[A]{B}\mid x_{3}=\left(1-\epsilon\right)y_{3}^{\frac{e_{\rmB\rmX}}{e_{\rmB\rmY}}}\right\}.
\end{equation*}
These curves lie in the domains \(\Gamma_{\rmX}\) and \(\Gamma_{\rmY}\), respectively.

We define
\begin{equation*}
  \mcl{D}=\left\{\left(X_{3},Y_{3}\right)\in\R^{2}\mid X_{3},Y_{3}<\log h\right\}\subsetneq\negativeR{2}.
\end{equation*}
In logarithmic coordinates, the curves \(\Sigma_{\rms}^{+}\) and \(\Sigma_{\rms}^{-}\) are represented by the affine subspaces
\begin{equation*}
  W_{\rms}^{+}=\left\{\left(X_{3},Y_{3}\right)\in\mcl{D}\mid Y_{3}=\log\left(1-\epsilon\right)+\frac{e_{\rmB\rmY}}{e_{\rmB\rmX}}X_{3}\right\}
\end{equation*}
and
\begin{equation*}
  W_{\rms}^{-}=\left\{\left(X_{3},Y_{3}\right)\in\mcl{D}\mid X_{3}=\log\left(1-\epsilon\right)+\frac{e_{\rmB\rmX}}{e_{\rmB\rmY}}Y_{3}\right\}.
\end{equation*}
These affine subspaces are parallel to linear subspace \(W_{\rms}\)
\begin{equation*}
  W_{\rms}=\left\{\left(X_{3},Y_{3}\right)\mid e_{\rmB\rmY}X_{3}=e_{\rmB\rmX}Y_{3}\right\},
\end{equation*}
which is the curve \(\Sigma_{\rms}\) expressed in logarithmic coordinates, and do not pass through the origin.

Let \(\vartheta\) be a point of \(S\) that is not \(\vartheta_{\rms}\). This point is the projection of any vector in the linear subspace \(W(\vartheta)\) defined by \(X_{4}=-\frac{1+\vartheta}{\vartheta}X_{3}\). However, not all points in \(W(\vartheta)\) lie in the domains of the transition matrix \(M_{\rmX}\) or \(M_{\rmY}\). Some points of \(W(\vartheta)\) lie in the excluded cusp \(\Gamma_{c}\) in logarithmic coordinates, bounded between the lines \(W_{\rms}^{+}\) and \(W_{\rms}^{-}\).

As \(W_{\rms}^{+}\) and \(W_{\rms}^{-}\) are parallel to \(W_{\rms}\), \(W(\vartheta)\) intersects \(W_{\rms}^{+}\) or \(W_{\rms}^{-}\) at some point \(p(\vartheta)\in\negativeR{2}\), given by
\begin{equation*}
  p(\vartheta)=\begin{cases}
    p_{\rmX}(\vartheta)=W(\vartheta)\cap W_{\rms}^{+}=\left(\frac{-\vartheta\log\left(1-\epsilon\right)}{\frac{e_{\rmB\rmY}}{e_{\rmB\rmX}}\vartheta+1+\vartheta},\frac{\left(1+\vartheta\right)\log\left(1-\epsilon\right)}{\frac{e_{\rmB\rmY}}{e_{\rmB\rmX}}\vartheta+1+\vartheta}\right) & \textrm{ if } \vartheta_{\rms}<\vartheta,\\[1em]
    p_{\rmY}(\vartheta)=W(\vartheta)\cap W_{\rms}^{-}=\left(\frac{\vartheta\log\left(1-\epsilon\right)}{\vartheta+\frac{e_{\rmB\rmX}}{e_{\rmB\rmY}}\left(1+\vartheta\right)},\frac{-\left(1+\vartheta\right)\log\left(1-\epsilon\right)}{\vartheta+\frac{e_{\rmB\rmX}}{e_{\rmB\rmY}}\left(1+\vartheta\right)}\right) & \textrm{ if } \vartheta<\vartheta_{\rms}.
  \end{cases}
\end{equation*}

\begin{figure}
  \centering
  \includegraphics[width=0.4\linewidth]{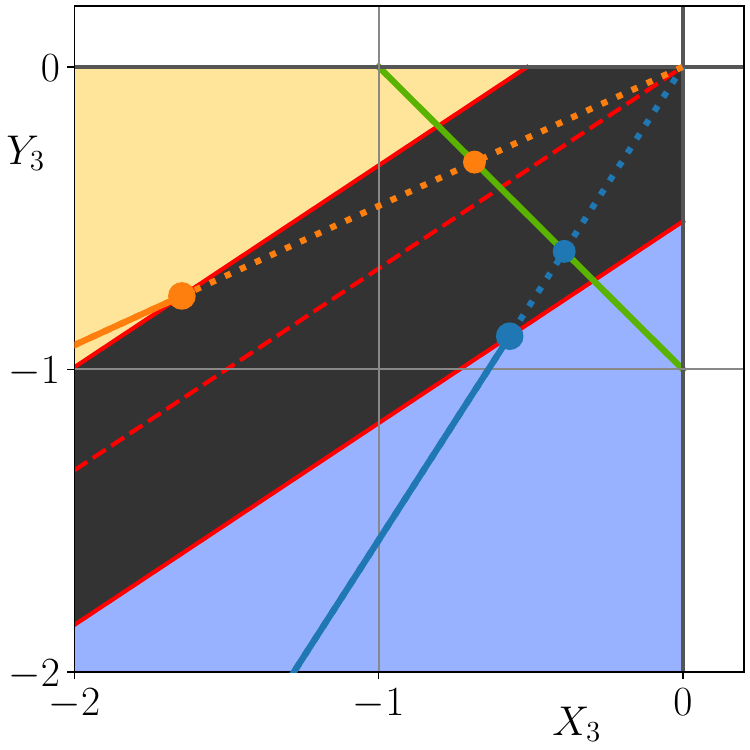}
  \caption{An example of the process of defining the projected map. The solid green line is the set \(S\) and the dashed red line represents the subspace \(W_{\rms}\), which projects to the switching manifold \(\vartheta_{\rms}\). The two solid red lines are \(W_{\rms}^{+}\) and \(W_{\rms}^{-}\). The two points lying on the green line are two examples of points \(\vartheta\in S\) which are not \(\vartheta_{\rms}\). The linear subspace through them is shown as dotted inside the excluded cusp, and solid outside. The projected map is only defined for vectors in the solid line, which represents the half-line \(L(\vartheta)\). All points in these subspaces project to the corresponding point in \(S\). The solid points at the intersection of the affine subspaces \(W_{\rms}^{+}\) and \(W_{\rms}^{-}\) and the linear subspaces \(W(\vartheta)\) are the points \(p(\vartheta)\). In this example, \(e_{\rmB\rmX}=1.2\) and \(e_{\rmB\rmY}=0.8\), and we set \(\epsilon=0.4\) for visual clarity.}
  \label{fig:well_def}
\end{figure}

We set
\begin{equation*}
  \alpha^{*}(\vartheta)=\begin{cases}
    \alpha_{\rmX}^{*}(\vartheta)=\frac{-\log\left(1-\epsilon\right)}{\frac{e_{\rmB\rmY}}{e_{\rmB\rmX}}\vartheta+1+\vartheta} & \textrm{ if } \vartheta_{\rms}<\vartheta,\\[0.75em]
    \alpha_{\rmY}^{*}(\vartheta)=\frac{\log\left(1-\epsilon\right)}{\vartheta+\frac{e_{\rmB\rmX}}{e_{\rmB\rmY}}\left(1+\vartheta\right)} & \textrm{ if } \vartheta<\vartheta_{\rms}.
  \end{cases}
\end{equation*}
All vectors in the half-line \(L(\vartheta)=\left\{\alpha\left(\vartheta,-\left(1+\vartheta\right)\right)\mid\alpha>\alpha^{*}(\vartheta)\right\}\subseteq{\negativeR{2}}\) are in the domain of either \(M_{\rmX}\) or \(M_{\rmY}\) and project onto \(\vartheta\). See \Fref{fig:well_def} for an example of this process.

Therefore, the projected map is well-defined as the projection onto \(S\) of the image under \(M\) of \(\alpha(\vartheta,-1-\vartheta)\) for any \(\alpha>\alpha^{*}(\vartheta)\):
\begin{equation*}
  f\colon\vartheta\mapsto \frac{-1}{(1,1)\cdot M\left(\alpha\left(\vartheta,-\left(1+\vartheta\right)\right)\right)}e_{1}\cdot M\left(\alpha\left(\vartheta,-\left(1+\vartheta\right)\right)\right).
\end{equation*}
In a less precise sense, this construction of the projected map means that, for initial conditions close to the switching manifold \(\vartheta_{\rms}\), the projected map describes the behaviour of corresponding trajectories in only small neighbourhoods of the network, and in the limit \(\vartheta\to\vartheta_{\rms}\), this neighbourhood vanishes. Full details of the above construction can be found in \cite[Appendix~A]{article}.

Lastly, we observe that
\begin{equation*}
  \lim_{X_{3}\to-\infty}\Pi\left(X_{3},\log\left(1-\epsilon\right)+\frac{e_{\rmB\rmY}}{e_{\rmB\rmX}}X_{3}\right)=\vartheta_{\rms};
\end{equation*}
that is, the projection of initial conditions restricted to the line \(W_{\rms}^{+}\) tends to the switching manifold \(\vartheta_{s}\) in the limit as these initial conditions approach the network, and that this is a one-sided limit from above. We similarly observe
\begin{equation*}
  \lim_{X_{3}\to-\infty}\Pi\left(X_{3},\frac{e_{\rmB\rmY}}{e_{\rmB\rmX}}\left(X_{3}-\log\left(1-\epsilon\right)\right)\right)=\vartheta_{\rms};
\end{equation*}
that is, again, the projection of initial conditions restricted to the line \(W_{\rms}^{-}\) tends to the switching manifold \(\vartheta_{s}\) in the limit as these initial conditions approach the network, but is now a one-sided limit from below. That initial conditions on these lines project to the switching manifold in the limit as these points approach the network is in fact true for all points restricted to lines parallel to the line \(W_{\rms}\).

Finally, as the projected maps  \(f_{\rmks}\) and \(f_{\Delta}\) are not defined at \(\vartheta_{\rms}\), we can calculate the limit \(\vartheta\to\vartheta_{\rms}\) by evaluating the limits
\begin{equation*}
  \lim_{X_{3}\to-\infty} M_{\rmX}\left(X_{3},\log\left(1-\epsilon\right)+\frac{e_{\rmB\rmY}}{e_{\rmB\rmX}}X_{3}\right)\equiv \lim_{\vartheta\searrow\vartheta_{\rms}}f\left(\vartheta\right)
\end{equation*}
and
\begin{equation*}
  \lim_{X_{3}\to-\infty} M_{\rmY}\left(X_{3},\frac{e_{\rmB\rmY}}{e_{\rmB\rmX}}\left(X_{3}-\log\left(1-\epsilon\right)\right)\right)\equiv \lim_{\vartheta\nearrow\vartheta_{\rms}}f\left(\vartheta\right).
\end{equation*}

We can understand these limits as evaluating the action on the map at the intersection point \(p\left(\vartheta\right)\) as this point goes to \(-\infty\) along the lines \(W_{\rms}^{+}\) and \(W_{\rms}^{-}\). We highlight now that evaluating the limits \eqref{eqn:ks_lim_right}, \eqref{eqn:ks_lim_left}, and \eqref{eqn:delta_lim} thus requires considering the flow near the two networks in the limit as a trajectory approaches the network.

\section{Completed return maps}\label{sec:ret_map}

We now construct for both networks a \textit{completed return map}---that is, one that captures all trajectories near the network. We follow the standard procedure to analyse the dynamics of trajectories near a heteroclinic network, which is to construct cross-sections transverse to the flow near the network, and construct maps between these sections. We adapt the methodology presented in \cite{kirk_lane_postlethwaite_rucklidge_silber_2010,kirk_postlethwaite_rucklidge_2012} to construct these maps. We compose these maps to derive a completed return map and, in the next section, we analyse this map to explain the (dis)continuity of the projected map of the two networks.

For convenience, we sometimes work in polar coordinates in the \(\left(x_{3},y_{3}\right)\)-plane, giving us coordinates \(\left(x_{1},x_{2},r_{3},\theta_{3}\right)\), where \(r_{3}\) and \(\theta_{3}\) are defined by
\begin{equation*}
  r_{3}^{2}=x_{3}^{2}+y_{3}^{2}
\end{equation*}
and
\begin{equation*}
  \tan\theta_{3}=\frac{y_{3}}{x_{3}}.
\end{equation*}
We also, when simpler, consider the angle \(\varphi_{3}=\frac{\pi}{2}-\theta_{3}\).

\begin{figure}
  \centering
  \begin{tikzpicture}
    \node at (0,0) () {\includegraphics[width=0.6\linewidth]{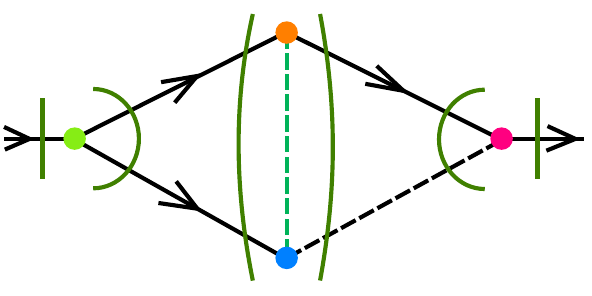}};

    \node at (3.4,-0.2) () {\(A\)};
    \node at (-3.7,-0.2) () {\(B\)};
    \node at (0,-2.2) () {\(X\)};
    \node at (0,2.2) () {\(Y\)};

    \node at (-4.3,1) () {\(\poinin[A]{B}\)};
    \node at (-2.6,1.15) () {\(\poinout[C]{B}\)};

    \node at (-1.35,-1.8) () {\(\poinin[B]{C}\)};
    \node at (1.15,-1.8) () {\(\poinout[A]{C}\)};

    \node at (2.5,1.15) () {\(\poinin[C]{A}\)};
    \node at (4.5,1) () {\(\poinout[B]{A}\)};
  \end{tikzpicture}
  \caption{A diagrammatic representation of both the Kirk--Silber and \(\Delta\)-clique networks, and the sections we define near them, showing how they capture all trajectories sufficiently close to the network. Because certain features differ between the two networks, the curve \(C\) is shown as a dashed line between \(X\) and \(Y\), as is the connection between \(X\) and \(A\). For visual clarity, we have not shown the entire connection \(A\to B\).}
  \label{fig:anal_secs}
\end{figure}
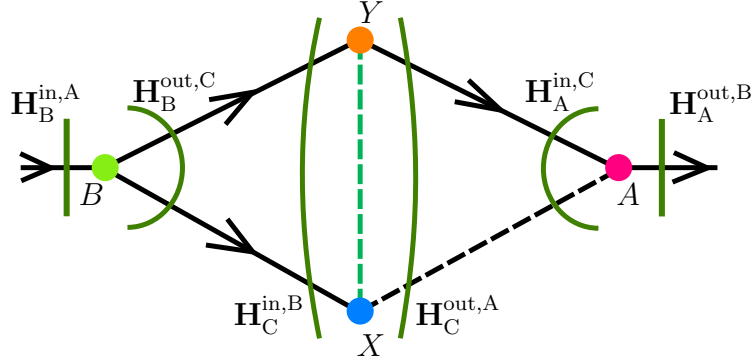

We use the following six cross-sections in our construction of the completed return map. Let \(h\) be a sufficiently small positive constant, so that all six of the following cross-sections are defined in a neighbourhood of approximate linear flow. First, near \(B\), we have
\begin{align*}
  \poinin[A]{B}&=\left\{(x_{1},x_{2},x_{3},y_{3})\mid x_{1}=h,|x_{2}-1|<h,0< x_{3},y_{3}<h\right\},\\
  \poinout[\rmC]{B}&=\left\{(x_{1},x_{2},r_{3},\theta_{3})\mid 0<x_{1}<h,|x_{2}-1|<h,r_{3}=h,0<\theta_{3}<\frac{\pi}{2}\right\}.
\end{align*}
Second, near \(C\),
\begin{align*}
  \poinin[B]{C}&=\left\{(x_{1},x_{2},r_{3},\theta_{3})\mid 0<x_{1}<h,x_{2}=h,|r_{3}-1|<h,0<\theta_{3}<\frac{\pi}{2}\right\},\\
  \poinout[\rmA]{C}&=\left\{(x_{1},x_{2},r_{3},\theta_{3})\mid x_{1}=h,0<x_{2}<h,r_{3}=h,0<\theta_{3}<\frac{\pi}{2}\right\}.
\end{align*}
Last, near \(A\),
\begin{align*}
  \poinin[C]{A}&=\left\{(x_{1},x_{2},r_{3},\theta_{3})\mid |x_{1}-1|<h,0<x_{2}<h,r_{3}=h,0<\theta_{3}<\frac{\pi}{2}\right\},\\
  \poinout[\rmB]{A}&=\left\{(x_{1},x_{2},x_{3},y_{3})\mid |x_{1}-1|<h,x_{2}=h,0< x_{3},y_{3}<h\right\}.
\end{align*}
\Fref{fig:anal_secs} gives a diagrammatic representation of these sections. We note that, by these definitions, no coordinate can be \(0\) for any trajectory that we consider, and therefore we do not consider the dynamics restricted to any proper subspace of the dynamical system.

We construct three completed local maps, which approximate the dynamics in small neighbourhoods of the equilibria \(A\) and \(B\) and the curve \(C\). These maps are
\begin{align*}
  \widehat{\psi}_{\rmB}&\colon\poinin[A]{B}\to\poinout[\rmC]{B},\\
  \widehat{\psi}_{\rmC}&\colon\poinin[B]{C}\to\poinout[\rmA]{C},\\
  \widehat{\psi}_{\rmA}&\colon\poinin[C]{A}\to\poinout[\rmB]{A}.
\end{align*}
We also construct three completed global maps, which approximate the dynamics along the unstable manifolds of the equilibria \(A\) and \(B\), and the curve \(C\):
\begin{align*}
  \widehat{\Psi}_{\rmB\rmC}&\colon\poinout[\rmC]{B}\to\poinin[B]{C},\\
  \widehat{\Psi}_{\rmC\rmA}&\colon\poinout[\rmA]{C}\to\poinin[C]{A},\\
  \widehat{\Psi}_{\rmA\rmB}&\colon\poinout[\rmB]{A}\to\poinin[A]{B}.
\end{align*}
The completed return map \(\widehat{\Phi}_{\rmB}\colon\poinin[A]{B}\to\poinin[A]{B}\) is then the composition
\begin{equation*}
  \widehat{\Phi}_{\rmB}=\widehat{\Psi}_{\rmA\rmB}\widehat{\psi}_{\rmA}\widehat{\Psi}_{\rmC\rmA}\widehat{\psi}_{\rmC}\widehat{\Psi}_{\rmB\rmC}\widehat{\psi}_{\rmB},
\end{equation*}
and we write \(\widehat{\Phi}_{\rmB}\colon\left(x_{3}^{\invar{B}},y_{3}^{\invar{B}}\right)\mapsto\left(x_{3}^{\prime},y_{3}^{\prime}\right)\).

For clarity, we say these maps are \textit{complete} in the sense that they account for all trajectories sufficiently close to the network, and do not require removing from their domain of definition an excluded cusp or its image. We represent completion of a map with a hat above it.

In many cases, because we want our maps to account for all trajectories near the network, we are unable to derive an explicit form for most of the above maps. However, we know enough about the dynamics in relevant subspaces to derive sufficient information of the components of this map for our analysis. By the invariant sphere theorem \cite{field_1996}, we know we can exclude the radial component from each step of our analysis, and so each of the above local and global maps is two-dimensional.

\subsection{The local map \(\widehat{\psi}_{\rmB}\colon\poinin[A]{B}\to\poinout[\rmC]{B}\)}\label{sec:compl:ssec:local_B}

We begin with a point \(\left(x_{3}^{\invar{B}},y_{3}^{\invar{B}}\right)\in\poinin[A]{B}\). Near \(B\), the linearised flow is

\begin{equation*}
  \dot{x}_{1}=-c_{\rmB}x_{1},\quad\dot{x}_{3}=e_{\rmB\rmX}x_{3},\quad\dot{y}_{3}=e_{\rmB\rmY}y_{3}.
\end{equation*}
Integrating these ODEs gives solutions \(x_{1}(t)\), \(x_{3}(t)\), and \(y_{3}(t)\), which can also be used to find expressions for \(r_{3}(t)\) and \(\theta_{3}(t)\). The trajectory crosses \(\poinout[\rmC]{B}\) when \(r_{3}(t)=h\), and so although the residence time \(T_{\rmB}\) cannot be solved for exactly, it satisfies the equation
\begin{equation}\label{eqn:TB}
  \begin{aligned}
    h^{2}
    &=x_{3}^{2}\rme^{2e_{\rmB\rmX}T_{\rmB}}+y_{3}^{2}\rme^{2e_{\rmB\rmY}T_{\rmB}}\\
    &=r_{3}^{2}\cos^{2}\theta_{3}\rme^{2e_{\rmB\rmX}T_{\rmB}}+r_{3}^{2}\sin^{2}\theta_{3}\rme^{2e_{\rmB\rmY}T_{\rmB}},
  \end{aligned}
\end{equation}
and we see that
\begin{equation*}
  \lim_{r_{3}^{\invar{B}}\to 0}T_{\rmB}=\infty.
\end{equation*}

Thus, the local map \(\widehat{\psi}_{\rmB}\colon\poinin[A]{B}\to\poinout[\rmC]{B}\) is defined by
\begin{equation*}
  \widehat{\psi}_{\rmB}\left(x_{3}^{\invar{B}},y_{3}^{\invar{B}}\right)=\left(x_{1}^{\outvar{B}},\theta_{3}^{\outvar{B}}\right),
\end{equation*}
where
\begin{equation*}
  x_{1}^{\outvar{B}}=h\rme^{-c_{\rmB}T_{\rmB}}
\end{equation*}
and
\begin{equation}\label{eqn:theta_B_out}
  \begin{aligned}
    \tan\theta_{3}^{\outvar{B}}
    &=\tan\theta_{3}\rme^{\left(e_{\rmB\rmY}-e_{\rmB\rmX}\right)T_{\rmB}}\\
    &=\frac{y_{3}^{\invar{B}}}{x_{3}^{\invar{B}}}\rme^{\left(e_{\rmB\rmY}-e_{\rmB\rmX}\right)T_{\rmB}}.
  \end{aligned}
\end{equation}

\subsection{The global map \(\widehat{\Psi}_{\rmB\rmC}\colon\poinout[\rmC]{B}\to\poinin[B]{C}\)}\label{sec:compl:ssec:global_BC}

The two-dimensional unstable manifold of \(B\) intersects \(\poinout[\rmC]{B}\) at \(\left(x_{1},r_{3},\theta_{3}\right)=\left(0,h,\theta_{3}\right)\), and intersects \(\poinin[B]{C}\) at \(\left(x_{1},x_{2},\theta_{3}\right)=\left(0,h,\bar{\theta}_{\rmC}(\theta_{3})\right)\), where \(\bar{\theta}_{\rmC}\) is an \(\order{1}\) function of \(\theta_{3}\). Therefore, the effect of the global map \(\widehat{\Psi}_{\rmB\rmC}\) is, to leading order, to rescale \(x_{1}\) by a \(\theta_{3}\)-dependent \(\order{1}\) amount, and to map the outgoing angle to an incoming angle, and so we have
\begin{equation*}
  \widehat{\Psi}_{\rmB\rmC}\left(x_{1}^{\outvar{B}},\theta_{3}^{\outvar{B}}\right)
  =\left(x_{1}^{\invar{C}},\theta_{3}^{\invar{C}}\right)
  =\left(
    E_{\rmC}\left(\theta_{3}^{\outvar{B}}\right)x_{1}^{\outvar{B}},
    \bar{\theta}_{\rmC}\left(\theta_{3}^{\outvar{B}}\right)
  \right).
\end{equation*}
Here, \(E_{\rmC}\) is an \(\order{1}\) function of \(\theta_{3}\). We know little about the function \(\bar{\theta}_{\rmC}\). However, due to the invariance of the \(x_{3}\)- and \(y_{3}\)-axes, we know that \(\bar{\theta}_{\rmC}\left(0\right)=0\) and \(\bar{\theta}_{\rmC}\left(\frac{\pi}{2}\right)=\frac{\pi}{2}\).

\subsection{The local map \(\widehat{\psi}_{\rmC}\colon\poinin[B]{C}\to\poinout[\rmA]{C}\)}\label{sec:compl:ssec:local_C}

The local map \(\widehat{\psi}_{\rmC}\) is the most difficult to construct, and is the only map that differs between the two networks. We need to construct this map in such a way that it captures all trajectories that begin on \(\poinin[\rmA]{\rmB}\), noting that, for any subset of the curve \(C\), there exist initial conditions on \(\poinin[\rmA]{\rmB}\) that give rise to trajectories that pass arbitrarily close to that subset of \(C\). As such, we cannot simply consider the linear flow near each equilibrium that lies in \(C\), as we usually do when constructing return maps near heteroclinic cycles. Here, we follow the methodology in \cite{kirk_lane_postlethwaite_rucklidge_silber_2010,kirk_postlethwaite_rucklidge_2012}.

Recall that, for the Kirk--Silber network, the curve \(C\) is the union of the equilibria \(X\), \(Y\), and \(P\), and the heteroclinic orbits \(P\to X\) and \(P\to Y\), and that, for the \(\Delta\)-clique network, \(C\) is the union of the equilibria \(X\) and \(Y\), and the heteroclinic orbit \(X\to Y\). For both networks, the curve is contained in the \(\left(x_{3},y_{3}\right)\)-plane.

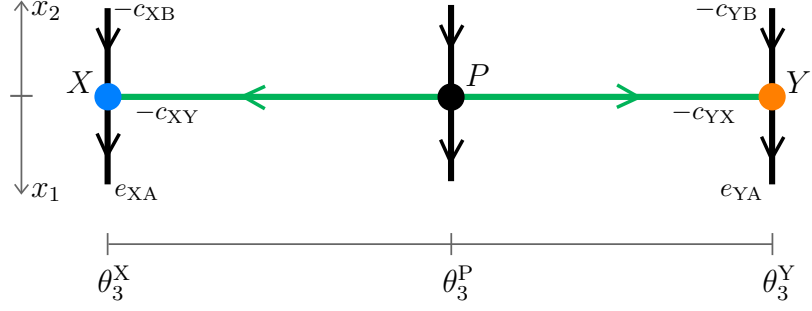
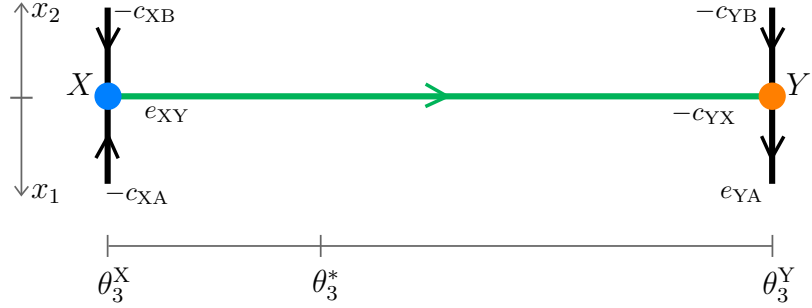
\begin{figure}
  \centering
  \begin{subfigure}[t]{0.7\linewidth}
    \centering
    \begin{tikzpicture}
      \node at (0,0) () {\includegraphics[width=\textwidth]{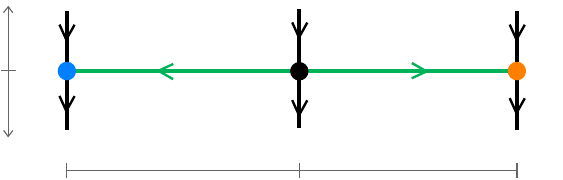}};

      \node at (-5., 1.5) () {\(x_{2}\)};
      \node at (-5.,-0.9) () {\(x_{1}\)};

      \node at (-4.55,0.55) () {\(X\)};
      \node at (0.7,  0.65) () {\(P\)};
      \node at (4.95, 0.55) () {\(Y\)};

      \node at (-4.1,-2.1)  () {\(\theta_{3}^{\rmX}\)};
      \node at ( 0.45,-2.1) () {\(\theta_{3}^{\rmP}\)};
      \node at ( 4.7,-2.1)  () {\(\theta_{3}^{\rmY}\)};

      \node at (-3.7,1.5)  () {\footnotesize\(-c_{\rmX\rmB}\)};
      \node at (-3.4,0.15) () {\footnotesize\(-c_{\rmXY}\)};
      \node at (-3.8,-0.9) () {\footnotesize\(e_{\rmX\rmA}\)};

      \node at (4,1.5)    () {\footnotesize\(-c_{\rmY\rmB}\)};
      \node at (3.7,0.15) () {\footnotesize\(-c_{\rmY\rmX}\)};
      \node at (4.2,-0.9) () {\footnotesize\(e_{\rmY\rmA}\)};
    \end{tikzpicture}
    \vspace{-5mm}
    \caption{Kirk--Silber network.}
    \vspace{10mm}
    \label{fig:C_dynamics:kirk_silber_net}
  \end{subfigure}

  \begin{subfigure}[t]{0.7\linewidth}
    \centering
    \begin{tikzpicture}
      \node at (0,0) () {\includegraphics[width=\textwidth]{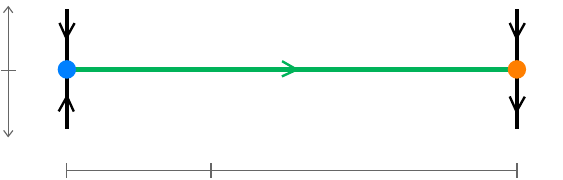}};

      \node at (-5., 1.5) () {\(x_{2}\)};
      \node at (-5.,-0.9) () {\(x_{1}\)};

      \node at (-4.55,0.55) () {\(X\)};
      \node at (4.95, 0.55) () {\(Y\)};

      \node at (-4.1,-2.1) () {\(\theta_{3}^{\rmX}\)};
      \node at (-1.3,-2.1) () {\(\theta_{3}^{*}\)};
      \node at ( 4.7,-2.1) () {\(\theta_{3}^{\rmY}\)};

      \node at (-3.7,1.5)  () {\footnotesize\(-c_{\rmX\rmB}\)};
      \node at (-3.4,0.15) () {\footnotesize\(e_{\rmXY}\)};
      \node at (-3.8,-0.9) () {\footnotesize\(-c_{\rmX\rmA}\)};

      \node at (4,1.5)    () {\footnotesize\(-c_{\rmY\rmB}\)};
      \node at (3.7,0.15) () {\footnotesize\(-c_{\rmY\rmX}\)};
      \node at (4.2,-0.9) () {\footnotesize\(e_{\rmY\rmA}\)};
    \end{tikzpicture}
    \vspace{-5mm}
    \caption{\(\Delta\)-clique network.}
    \label{fig:C_dynamics:delta_clique_net}
  \end{subfigure}%
  \caption{A representation of the dynamics near the curve \(C\) for (a) the Kirk--Silber network, and (b) the \(\Delta\)-clique network, showing the difference in how trajectories pass near the curve \(C\). For both figures, the vertical axis represents the \(x_{2}\) coordinate above the curve \(C\), and the \(x_{1}\) coordinate below the curve \(C\). The horizontal axis represents the angle \(\theta_{3}\), with certain relevant angles labelled. The value of the angle \(\theta_{3}^{*}\) in (\subref{fig:C_dynamics:delta_clique_net}) is indicative only, as its value depends on the values of the parameters in \Tref{tbl:eigs}(b). The value used here is the same as calculated numerically in \Fref{fig:theta_3_C_out_v_theta_3_C_in}(\subref{fig:theta_3_C_out_v_theta_3_C_in:delta_clique_net}).}
  \label{fig:C_dynamics}
\end{figure}

We assumed that the curve \(C\) could be parametrised by the angle \(\theta_{3}\). We also know, by the invariant sphere theorem, that in the \(\left(x_{3},y_{3}\right)\)-plane there is strong contraction onto the curve \(C\). Therefore, the evolution of \(\theta_{3}\) can be described by a \(\theta_{3}\)-dependent function. The local flow near \(C\) is thus given by
\begin{equation*}
  \dot{x}_{1}=g_{\rmA}(\theta_{3})x_{1},\quad\dot{x}_{2}=-g_{\rmB}(\theta_{3})x_{2},\quad\dot{\theta}_{3}=g_{\theta}\left(\theta_{3}\right),
\end{equation*}
where \(g_{\rmA},g_{\rmB},g_{\theta}\colon\left[0,\frac{\pi}{2}\right]\to\R\) are functions of \(\theta_{3}\), and depend on the particular network being considered. A diagrammatic representation of the flow near the curve \(C\) is shown in \Fref{fig:C_dynamics}. For the Kirk--Silber network, we set \(\theta_{3}^{\rmP}\) as the value of \(\theta_{3}\) at the equilibrium \(P\) on the curve \(C\).

We have little precise information about the functions \(g_{\rmA}\), \(g_{\rmB}\), and \(g_{\theta}\). However, in the case of the Kirk--Silber network, we know the following:
\begin{enumerate}[label=(\alph*)]
  \item \(x_{1}\) expands for all values of \(\theta_{3}\), and so \(g_{\rmA}\left(\theta_{3}\right)>0\) for all \(\theta_{3}\in\left[0,\frac{\pi}{2}\right]\). The only known values of \(g_{\rmA}\) are \(g_{\rmA}(0)=e_{\rmX\rmA}\) and \(g_{\rmA}\left(\frac{\pi}{2}\right)=e_{\rmY\rmA}\).
  \item \(x_{2}\) contracts for all values of \(\theta_{3}\), and so \(g_{\rmB}\left(\theta_{3}\right)>0\) for all \(\theta_{3}\in\left[0,\frac{\pi}{2}\right]\). The only known values of \(g_{\rmB}\) are \(g_{\rmB}(0)=c_{\rmX\rmB}\) and \(g_{\rmB}\left(\frac{\pi}{2}\right)=c_{\rmY\rmB}\).
  \item The invariance of the \(x_{3}\)- and \(y_{3}\)-axes implies that \(g_{\theta}(0)=g_{\theta}\left(\frac{\pi}{2}\right)=0\). Moreover, the existence of the equilibrium \(P\) implies that \(g_{\theta}\left(\theta_{3}^{\rmP}\right)=0\). The heteroclinic orbits \(P\to X\) and \(P\to Y\) imply that, respectively, \(g_{\theta}\left(\theta_{3}\right)<0\) for all \(\theta_{3}\in\left(0,\theta_{3}^{\rmP}\right)\), and \(g_{\theta}\left(\theta_{3}\right)>0\) for all \(\theta_{3}\in\left(\theta_{3}^{\rmP},\frac{\pi}{2}\right)\).
  \item The linearisation of the vector field about \(X\) and \(Y\) implies that \(\rmD g_{\theta}(0)=-c_{\rmXY}<0\) and \(\rmD g_{\theta}\left(\frac{\pi}{2}\right)=-c_{\rmY\rmX}<0\). It follows that, for small \(|\theta_{3}|\), \(\dot{\theta}_{3}\approx -c_{\rmXY}\theta_{3}\), and that for small \(\left|\frac{\pi}{2}-\theta_{3}\right|\), \(\dot{\varphi}_{3}\approx-c_{\rmY\rmX}\varphi_{3}\). (Recall that we set \(\varphi_{3}=\frac{\pi}{2}-\theta_{3}\).)
\end{enumerate}
The following, mostly similar, properties are known in the case of the \(\Delta\)-clique network:
\begin{enumerate}[label=(\alph*)]
  \item There exists an angle \(\theta_{3}^{*}\) such that \(x_{1}\) contracts for all \(\theta_{3}\in\left[0,\theta_{3}^{*}\right)\), and expands for all \(\theta_{3}\in\left(\theta_{3}^{*},\frac{\pi}{2}\right]\). (See Figures~\ref{fig:dynamics:D:x1_x3_y3} and \ref{fig:C_dynamics:delta_clique_net}). Therefore, \(g_{\rmA}\left(\theta_{3}\right)<0\) for all \(\theta_{3}\in\left[0,\theta_{3}^{*}\right)\) and \(g_{\rmA}\left(\theta_{3}\right)>0\) for all \(\theta_{3}\in\left(\theta_{3}^{*},\frac{\pi}{2}\right]\). The only known values of \(g_{\rmA}\) are \(g_{\rmA}(0)=-c_{\rmX\rmA}\) and \(g_{\rmA}\left(\frac{\pi}{2}\right)=e_{\rmY\rmA}\).
  \item \(x_{2}\) contracts for all values of \(\theta_{3}\), and so \(g_{\rmB}\left(\theta_{3}\right)>0\) for all \(\theta_{3}\in\left[0,\frac{\pi}{2}\right]\). The only known values of \(g_{\rmB}\) are \(g_{\rmB}(0)=c_{\rmX\rmB}\) and \(g_{\rmB}\left(\frac{\pi}{2}\right)=c_{\rmY\rmB}\).
  \item The invariance of the \(x_{3}\)- and \(y_{3}\)-axes implies that \(g_{\theta}(0)=g_{\theta}\left(\frac{\pi}{2}\right)=0\). Moreover, the heteroclinic orbit \(X\to Y\) implies that \(g_{\theta}\left(\theta_{3}\right)>0\) for all \(\theta_{3}\in\left(0,\frac{\pi}{2}\right)\).
  \item The linearisation of the vector field about \(X\) and \(Y\) implies that \(\rmD g_{\theta}(0)=e_{\rmXY}>0\) and \(\rmD g_{\theta}\left(\frac{\pi}{2}\right)=-c_{\rmY\rmX}<0\). It follows that, for small \(|\theta_{3}|\), \(\dot{\theta}_{3}\approx e_{\rmXY}\theta_{3}\), and that for small \(\left|\frac{\pi}{2}-\theta_{3}\right|\), \(\dot{\varphi}_{3}\approx-c_{\rmY\rmX}\varphi_{3}\).
\end{enumerate}
We make the reasonable simplifying assumption for both networks that the sign of \(g_{\rmB}\) does not change over \(\theta_{3}\in\left[0,\frac{\pi}{2}\right]\). We also assume that \(g_{\rmA}\) does not change sign over \(\theta_{3}\in\left[0,\frac{\pi}{2}\right]\) for the Kirk--Silber network, and only changes sign once over \(\theta_{3}\in\left[0,\frac{\pi}{2}\right]\), at \(\theta_{3}^{*}\), for the \(\Delta\)-clique network. These assumptions hold for the ODEs \eqref{eqn:ODEs}. With the addition of higher order terms to \eqref{eqn:ODEs}, these assumptions could be broken with particular parameter values.

To derive the local map \(\widehat{\psi}_{\rmC}\), we first use separation of variables to solve \(\dot{\theta}_{3}\), and find \(\theta_{3}^{\rmC}(t)\) implicitly:
\begin{equation*}
  t=\int_{\theta_{3}^{\rmC}\left(0\right)}^{\theta_{3}^{\rmC}\left(t\right)}\frac{1}{g_{\theta}\left(\theta\right)}\rmd\theta.
\end{equation*}

Next, the trajectory crosses \(\poinout[\rmA]{\rmC}\) when \(x_{1}(t)=h\), and so although we again cannot solve for the residence time \(T_{\rmC}\) explicitly, it satisfies the equation
\begin{equation}\label{eqn:TC_eq}
  \int_{0}^{T_{\rmC}}g_{\rmA}\left(\theta_{3}\left(\tau\right)\right)\rmd\tau=-\log\left(\frac{x_{1}\left(0\right)}{h}\right)=-\log\left(E_{\rmC}\left(\theta_{3}^{\outvar{B}}\right)\right)+c_{\rmB}T_{\rmB}.
\end{equation}
From this relationship we see that
\begin{equation*}
  \lim_{x_{1}^{\invar{C}}\to 0}T_{\rmC}=\infty.
\end{equation*}

We write \(\theta_{3}^{\invar{C}}\equiv\theta_{3}^{\rmC}\left(0\right)\) and \(\theta_{3}^{\outvar{C}}\equiv\theta_{3}^{\rmC}\left(T_{\rmC}\right)\), and we note that \(\theta_{3}^{\outvar{C}}\) is a function of both \(\theta_{3}^{\invar{C}}\) and \(T_{\rmC}\): \(\theta_{3}^{\outvar{C}}\left(\theta_{3}^{\invar{C}},T_{\rmC}\right)\).

In the case of the \(\Delta\)-clique network, we know that \(g_{\theta}\left(\theta_{3}\right)>0\) for all
\(\theta_{3}\in\left(0,\frac{\pi}{2}\right)\). As such, trajectories that spend a sufficiently long time in a small neighbourhood of the curve \(C\) leave that small neighbourhood close to \(Y\), where \(\theta_{3}=\frac{\pi}{2}\). Therefore, we know that, for a fixed initial condition \(\theta_{3}^{\invar{C}}\),
\begin{equation*}
  \lim_{T_{\rmC}\to \infty}\theta_{3}^{\outvar{C}}\left(\theta_{3}^{\invar{C}},T_{\rmC}\right)=\frac{\pi}{2}.
\end{equation*}
However, since a trajectory strikes \(\poinout[A]{C}\) when \(x_{1}\left(T_{\rmC}\right)=h\), and \(x_{1}\) decays for all \(\theta_{3}<\theta_{3}^{*}\), the value of \(\theta_{3}^{\outvar{C}}\) is bounded below; that is, since \(g_{\rmA}\left(\theta_{3}^{\rmC}\right)<0\) for all \(\theta_{3}<\theta_{3}^{*}\), then, regardless of \(T_{\rmC}\), \(\theta_{3}^{\outvar{C}}>\theta_{3}^{*}\).

In the case of the Kirk--Silber network, trajectories that enter a small neighbourhood of \(C\) along the heteroclinic orbit \(B\to P\) leave near \(P\), where \(\theta_{3}=\theta_{3}^{\rmP}\). All other trajectories that spend a sufficiently long time in a small neighbourhood of \(C\) leave that small neighbourhood close to \(X\) or \(Y\), where \(\theta_{3}=0\) or \(\theta_{3}=\frac{\pi}{2}\), respectively, depending on whether they strike \(\poinin[B]{C}\) between \(X\) and \(P\) or \(P\) and \(Y\), respectively. In particular, since \(g_{\theta}\left(\theta_{3}\right)<0\) for all \(\theta_{3}\in\left(0,\theta_{3}^{\rmP}\right)\), we know that for a fixed initial condition \(\theta_{3}^{\invar{C}}\in\left(0,\theta_{3}^{\rmP}\right)\),
\begin{equation*}
  \lim_{T_{\rmC}\to \infty}\theta_{3}^{\outvar{C}}\left(\theta_{3}^{\invar{C}},T_{\rmC}\right)=0,
\end{equation*}
and
since \(g_{\theta}\left(\theta_{3}\right)>0\) for all \(\theta_{3}\in\left(\theta_{3}^{\rmP},\frac{\pi}{2}\right)\), we know that for a fixed initial condition \(\theta_{3}^{\invar{C}}\in\left(\theta_{3}^{\rmP},\frac{\pi}{2}\right)\),
\begin{equation*}
  \lim_{T_{\rmC}\to \infty}\theta_{3}^{\outvar{C}}\left(\theta_{3}^{\invar{C}},T_{\rmC}\right)=\frac{\pi}{2}.
\end{equation*}
Since \(g_{\theta}\left(\theta_{3}^{\rmP}\right)=0\), then, if \(\theta_{3}^{\invar{C}}=\theta_{3}^{\rmP}\), \(\theta_{3}^{\outvar{C}}=\theta_{3}^{\rmP}\), regardless of \(T_{\rmC}\).

\begin{figure}
  \centering
  \begin{subfigure}[t]{0.48\linewidth}
    \centering
    \includegraphics[width=\textwidth]{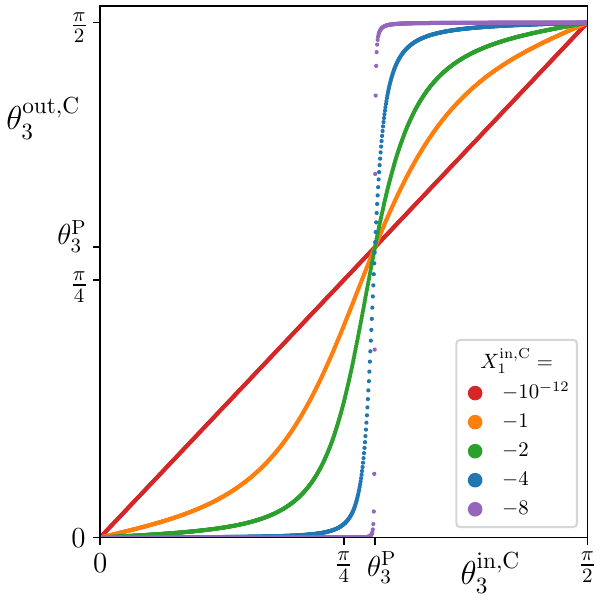}
    \caption{Kirk--Silber network.}
    \label{fig:theta_3_C_out_v_theta_3_C_in:kirk_silber_net}
  \end{subfigure}%
  \hfill%
  \begin{subfigure}[t]{0.48\linewidth}
    \centering
    \includegraphics[width=\textwidth]{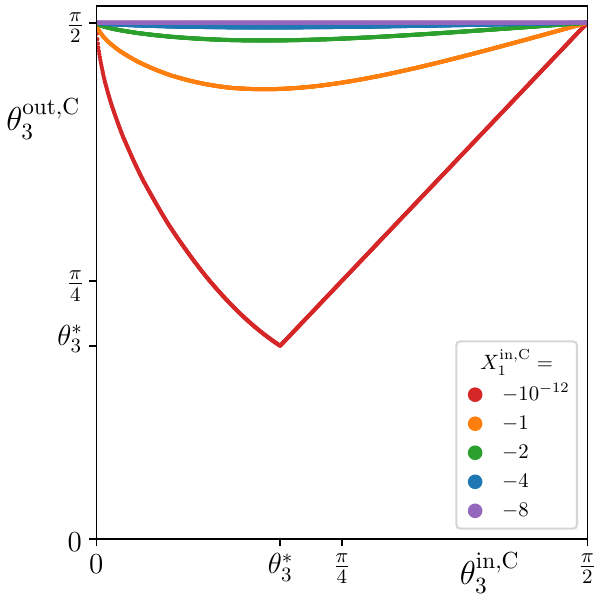}
    \caption{\(\Delta\)-clique network.}
    \label{fig:theta_3_C_out_v_theta_3_C_in:delta_clique_net}
  \end{subfigure}%
  \caption{
    Plots of the outgoing angle \(\theta_{3}^{\outvar{C}}\) as a function of \(\theta_{3}^{\invar{C}}\) for both the Kirk--Silber and \(\Delta\)-clique networks.
    The relationship is computed for various values of \(x_{1}^{\invar{C}}\), given as \(X_{1}^{\invar{C}}\coloneqq\log x_{1}^{\invar{C}}\).
    These plots were computed numerically by integrating ODEs that contain the heteroclinic networks. We see that the function in the case of the Kirk--Silber network is an identity function for small values of \(|X_{1}^{\invar{C}}|\), and converges pointwise to a step function as \(X_{1}\to-\infty\). The function in the case of the \(\Delta\)-clique case is more complicated for small values of \(|X_{1}^{\invar{C}}|\), but converges pointwise to a constant function as \(X_{1}^{\invar{C}}\to-\infty\), which corresponds to \(r_{3}^{\invar{B}}\to0\).}
  \label{fig:theta_3_C_out_v_theta_3_C_in}
\end{figure}

From these relationships, we have our first indication of the origin of the discontinuity of the projected map in the case of the Kirk--Silber network. In particular, considering the outgoing angle \(\theta_{3}^{\outvar{C}}\equiv\theta_{3}^{\rmC}\left(T_{\rmC}\right)\) as a function of the incoming angle \(\theta_{3}^{\invar{C}}\equiv\theta_{3}^{\rmC}(0)\), where we take the \(x_{1}^{\invar{C}}\)-dependent expression \(T_{\rmC}\) as a parameter, we have, in the case of the Kirk--Silber network,
\begin{equation}\label{eqn:ks_Cout_Cin_lim}
  \lim_{x_{1}^{\invar{C}}\to 0}\theta_{3}^{\outvar{C}}\left(\theta_{3}^{\invar{C}}\right)=\begin{cases}
    0,&\ \mathrm{if}\ \theta_{3}^{\invar{C}}<\theta_{3}^{\rmP},\\
    \theta_{3}^{\rmP},&\ \mathrm{if}\ \theta_{3}^{\invar{C}}=\theta_{3}^{\rmP},\\
    \frac{\pi}{2},&\ \mathrm{if}\ \theta_{3}^{\invar{C}}>\theta_{3}^{\rmP},
  \end{cases}
\end{equation}
and, for the \(\Delta\)-clique network,
\begin{equation}\label{eqn:delta_clique_Cout_Cin_lim}
  \lim_{x_{1}^{\invar{C}}\to 0}\theta_{3}^{\outvar{C}}\left(\theta_{3}^{\invar{C}}\right)=\frac{\pi}{2}.
\end{equation}
The function \(\theta_{3}^{\outvar{C}}\left(\theta_{3}^{\invar{C}}\right)\) is a continuous function for both networks, since \(x_{1}^{\invar{C}}\neq 0\) for all trajectories that pass near \(C\). However, in the limit as the trajectory approaches the network, the function is continuous for the \(\Delta\)-clique network---being a constant function---and discontinuous for the Kirk--Silber network, being a step function. In the case of the Kirk--Silber network, a separatrix is formed by the stable manifold of \(P\), and is defined as the surface of trajectories that satisfy \(\theta_{3}^{\outvar{C}}=\theta_{3}^{\invar{C}}\). No such separatrix exists near the \(\Delta\)-clique network.

Although we cannot analytically derive the function \(\theta_{3}^{\outvar{C}}\left(\theta_{3}^{\invar{C}}\right)\), we can, with numerical computation, observe the limits in \eqref{eqn:ks_Cout_Cin_lim} and \eqref{eqn:delta_clique_Cout_Cin_lim}. We show in \Fref{fig:theta_3_C_out_v_theta_3_C_in} plots of these functions for various values of \(X_{1}^{\invar{C}}\coloneqq\log x_{1}^{\invar{C}}\).

Consider first the Kirk--Silber network, shown in \Fref{fig:theta_3_C_out_v_theta_3_C_in}(\subref{fig:theta_3_C_out_v_theta_3_C_in:kirk_silber_net}). For values of \(x_{1}^{\invar{C}}\approx1\), and so \(X_{1}^{\invar{C}}\approx0\), \(T_{C}\) is small by \eqref{eqn:TC_eq}. As such, \(\theta_{3}\left(t\right)\) evolves for only a brief time, and so \(\theta_{3}^{\outvar{C}}\approx\theta_{3}^{\invar{C}}\). We observe in \Fref{fig:theta_3_C_out_v_theta_3_C_in}(\subref{fig:theta_3_C_out_v_theta_3_C_in:kirk_silber_net}) that the trajectories represented by red points form an approximately straight line close to the identity function \(\theta_{3}^{\invar{C}}=\theta_{3}^{\outvar{C}}\). As \(x_{1}^{\invar{C}}\to0\), and \(X_{1}^{\invar{C}}\to-\infty\), the residence time \(T_{C}\) grows, and \(\theta_{3}\left(t\right)\) evolves according to \(g_{\theta}\) for a longer period of time, towards either \(0\) or \(\frac{\pi}{2}\), as can be seen in the orange, green, blue, and purple points which represent successively smaller values of \(x_{1}^{\invar{C}}\) (and larger values of \(|X_{1}^{\invar{C}}|\)). Thus, the value of \(\theta_{3}^{\outvar{C}}\) approaches a step-function, as trajectories that enter a small neighbourhood of the curve \(C\) near \(P\) have a greater period of time to follow the heteroclinic orbits \(P\to X\) or \(P\to Y\), and so leave that small neighbourhood of \(C\) closer to \(X\) or \(Y\).

Now consider the \(\Delta\)-clique network, shown in \Fref{fig:theta_3_C_out_v_theta_3_C_in}(\subref{fig:theta_3_C_out_v_theta_3_C_in:delta_clique_net}). For values of \(x_{1}^{\invar{C}}\approx1\), and so \(X_{1}^{\invar{C}}\approx0\), \(T_{\rmC}\) is small by \eqref{eqn:TC_eq} if \(\theta_{3}^{\invar{C}}>\theta_{3}^{*}\), as, for these values of \(\theta_{3}\), \(x_{1}\) expands. Thus for \(\theta_{3}^{\invar{C}}>\theta_{3}^{*}\), \(\theta_{3}^{\outvar{C}}\approx\theta_{3}^{\invar{C}}\), and the red points with \(\theta_{3}^{\invar{C}}>\theta_{3}^{*}\) form an approximately straight line. However, for \(\theta_{3}^{\invar{C}}<\theta_{3}^{*}\), \(x_{1}\) decays until \(\theta_{3}\left(t\right)\) evolves according to \(g_{\theta}\) to be greater than \(\theta_{3}^{*}\), and \(x_{1}\) begins expanding. For \(\theta_{3}^{\invar{C}}\ll 1\), \(\theta_{3}^{\outvar{C}}\) is close to \(\frac{\pi}{2}\), even for \(x_{1}^{\invar{C}}\approx 1\), as seen in the red points with \(\theta_{3}^{\invar{C}}<\theta_{3}^{*}\). As \(x_{1}^{\invar{C}}\to0\), and \(X_{1}^{\invar{C}}\to-\infty\), the residence time \(T_{C}\) grows, and the value of \(\theta_{3}^{\outvar{C}}\) approaches a constant \(\frac{\pi}{2}\) as \(\theta_{3}\left(t\right)\) evolves according to \(g_{\theta}\), again as seen in the orange, green, blue, and purple points which represent successively smaller values of \(x_{1}^{\invar{C}}\). The value of \(\theta_{3}^{\outvar{C}}\) approaches a constant \(\frac{\pi}{2}\) as trajectories that enter a small neighbourhood of the curve \(C\) at any point have a greater period of time to follow the heteroclinic orbits \(X\to Y\), and so leave that small neighbourhood of \(C\) closer to \(Y\).

In \ref{app:numerics:ssec:theta_in_theta_out}, we describe how these figures were computed.

With assumptions about the functions \(g_{\rmA}\), \(g_{\rmB}\) and \(g_{\theta}\), we could derive exact expressions for components of this local map. In \cite{kirk_postlethwaite_rucklidge_2012}, for example, they assumed that \(g_{\theta}\left(\theta_{3}\right)=-\frac{\lambda}{4}\sin4\theta_{3}\) for the Kirk--Silber network. However, we do not need to make such assumptions here, as the asymptotic behaviour of these components is sufficient for our purposes.

In summary, we write the local map as
\begin{equation*}
  \widehat{\psi}_{\rmC}\left(x_{1}^{\invar{C}},\theta_{3}^{\invar{C}}\right)=\left(x_{2}^{\outvar{C}},\theta_{3}^{\rmC}\left(T_{\rmC}\right)\right),
\end{equation*}
where
\begin{equation*}
  x_{2}^{\outvar{C}}=h\exp\left(-\int_{0}^{T_{\rmC}}g_{\rmB}(\theta_{3}^{\rmC}\left(\tau\right))\rmd\tau\right).
\end{equation*}

\subsection{The global map \(\widehat{\Psi}_{\rmC\rmA}\colon\poinout[\rmA]{C}\to\poinin[C]{A}\)}\label{sec:compl:ssec:global_CA}

The two-dimensional unstable manifold of \(C\) intersects \(\poinout[\rmA]{C}\) at \(\left(x_{1},x_{2},\theta_{3}\right)=\left(h,0,\theta_{3}\right)\) and intersection \(\poinin[C]{A}\) at \(\left(x_{2},r_{3}\theta_{3}\right)=\left(0,h,\bar{\theta}_{\rmA}\left(\theta_{3}\right)\right)\). Again, we then know that the effect of the global map \(\widehat{\Psi}_{\rmC\rmA}\) is, to leading order, to rescale \(x_{2}\) by a \(\theta_{3}\)-dependent \(\order{1}\) amount, and to map the outgoing angle to an incoming angle, and so we have
\begin{equation*}
  \widehat{\Psi}_{\rmC\rmA}\left(x_{2}^{\outvar{C}},\theta_{3}^{\outvar{C}}\right)
  =\left(x_{2}^{\invar{A}},\theta_{3}^{\invar{A}}\right)
  =\left(
    E_{\rmA}\left(\theta_{3}^{\outvar{C}}\right)x_{2}^{\outvar{C}},
    \bar{\theta}_{\rmA}\left(\theta_{3}^{\outvar{C}}\right)
  \right).
\end{equation*}
As for the global map \(\widehat{\Psi}_{\rmB\rmC}\), \(E_{\rmA}\) is an \(\order{1}\) function of \(\theta_{3}\). We also know little about the function \(\bar{\theta}_{\rmA}\), except again that \(\bar{\theta}_{\rmA}\left(0\right)=0\) and \(\bar{\theta}_{\rmA}\left(\frac{\pi}{2}\right)=\frac{\pi}{2}\).

\subsection{The local map \(\widehat{\psi}_{\rmA}\colon\poinin[C]{A}\to\poinout[\rmB]{A}\)}\label{sec:compl:ssec:local_A}

The local map \(\widehat{\psi}_{\rmA}\) is the most straightforward to compute. The flow linearised near \(A\) is
\begin{equation*}
  \dot{x}_{2}=e_{\rmA}x_{2},\quad\dot{x}_{3}=-c_{\rmA\rmX}x_{3},\quad\dot{y}_{3}=-c_{\rmA\rmY}y_{3}.
\end{equation*}
The trajectory crosses \(\poinout[\rmB]{A}\) when \(x_{2}(t)=h\), and therefore the residence time \(T_{\rmA}\) can be found explicitly:
\begin{equation*}
  T_{\rmA}=-\frac{1}{e_{\rmA}}\log\left(\frac{x_{2}}{h}\right).
\end{equation*}
Thus, the local map is
\begin{equation*}
  \widehat{\psi}_{\rmA}\left(x_{2}^{\invar{A}},\theta_{3}^{\invar{A}}\right)
  =\left(x_{3}^{\outvar{A}},y_{3}^{\outvar{A}}\right)
  =\left(
    h\cos\theta_{3}^{\invar{A}}\left(\frac{x_{2}^{\invar{A}}}{h}\right)^{\frac{c_{\rmA\rmX}}{e_{\rmA}}},
    h\sin\theta_{3}^{\invar{A}}\left(\frac{x_{2}^{\invar{A}}}{h}\right)^{\frac{c_{\rmA\rmY}}{e_{\rmA}}}
  \right).
\end{equation*}

\subsection{The global map \(\widehat{\Psi}_{\rmA\rmB}\colon\poinout[\rmB]{A}\to\poinin[A]{B}\)}\label{sec:compl:ssec:global_AB}

The heteroclinic orbit from \(A\) to \(B\) intersects \(\poinout[\rmB]{A}\) at \(\left(x_{2},x_{3},y_{3}\right)=\left(h,0,0\right)\), and intersects \(\poinin[A]{B}\) at \(\left(x_{1},x_{3},y_{3}\right)=\left(h,0,0\right)\). The effect of the global map \(\widehat{\Phi}_{\rmA\rmB}\) is thus, to leading order, to multiply the coordinates \(x_{3}\) and \(y_{3}\) by an \(\order{1}\) amount. The expression for this map is therefore
\begin{equation*}
  \widehat{\Psi}_{\rmA\rmB}\left(x_{3}^{\outvar{A}},y_{3}^{\outvar{A}}\right)
  =\left({x_{3}^{\invar{B}}}^{\prime},{y_{3}^{\invar{B}}}^{\prime}\right)
  =\left(
    E_{\rmB\rmX}x_{3}^{\outvar{A}},
    E_{\rmB\rmY}y_{3}^{\outvar{A}}
  \right).
\end{equation*}

\subsection{Composing the maps}\label{sec:compl:ssec:compose}

Composing these maps together, we can write the completed return map \(\widehat{\Phi}_{\rmB}\colon\poinin[A]{B}\to\poinin[A]{B}\) as
\begin{equation*}
  \widehat{\Phi}_{\rmB}\left(x_{3}^{\invar{B}},y_{3}^{\invar{B}}\right)=(x_{3}^{\prime},y_{3}^{\prime}),
\end{equation*}
where
\begin{equation}\label{eqn:x3_prime}
  x_{3}^{\prime}=E_{\rmB\rmX}h\cos\left(\bar{\theta}_{\rmA}\left(\theta_{3}^{\rmC}\left(T_{\rmC}\right)\right)\right)\left(E_{\rmA}\left(\theta_{3}^{\rmC}\left(T_{\rmC}\right)\right)\exp\left(-\int_{0}^{T_{\rmC}}g_{\rmB}(\theta_{3}^{\rmC}\left(\tau\right))\rmd\tau\right)\right)^{\frac{c_{\rmA\rmX}}{e_{\rmA}}}
\end{equation}
and
\begin{equation}\label{eqn:y3_prime}
  y_{3}^{\prime}=E_{\rmB\rmY}h\sin\left(\bar{\theta}_{\rmA}\left(\theta_{3}^{\rmC}\left(T_{\rmC}\right)\right)\right)\left(E_{\rmA}\left(\theta_{3}^{\rmC}\left(T_{\rmC}\right)\right)\exp\left(-\int_{0}^{T_{\rmC}}g_{\rmB}(\theta_{3}^{\rmC}\left(\tau\right))\rmd\tau\right)\right)^{\frac{c_{\rmA\rmY}}{e_{\rmA}}}.
\end{equation}
The value of \(\theta_{3}^{\rmC}\left(T_{\rmC}\right)\) is found by integrating \(g_{\theta}\) with initial condition
\begin{equation*}
  \theta_{3}^{\rmC}\left(0\right)=\bar{\theta}_{\rmC}\left(\frac{y_{3}^{\invar{B}}}{x_{3}^{\invar{B}}}\rme^{\left(e_{\rmB\rmY}-e_{\rmB\rmX}\right)T_{\rmB}}\right).
\end{equation*}
These expressions for \(x_{3}^{\prime}\) and \(y_{3}^{\prime}\) are valid for both the Kirk--Silber and \(\Delta\)-clique networks. However, each network has a different definition of \(g_{\rmA}\), \(g_{\rmB}\), and \(g_{\theta}\), and of the \(\order{1}\) components of the global maps. Therefore, the actual evaluation of these two expressions differs between the two networks.

This return map \(\widehat{\Phi}_{\rmB}\) is defined on the entirety of the incoming cross-section \(\poinin[\mrm{A}]{\mrm{B}}\), including inside the cusp \(\Gamma_{\rmc}\) excluded from the domain of \(\Phi_{\rmB}\) in \eqref{eqn:ks_ret_map} and \eqref{eqn:delta_ret_map}. Therefore, it is also defined on the switching curve \(\Sigma_{\rms}\) defined by \(x_{3}^{e_{\rmB\rmY}}=y_{3}^{e_{\rmB\rmX}}\), which corresponds to the switching manifold \(\vartheta_{\rms}\) of the projected maps \(f_{\rmks}\) and \(f_{\Delta}\). Moreover, it is a continuous map.

\section{Analysis of the completed return map}\label{sec:anal}

We now analyse the completed return map \(\widehat{\Phi}_{\rmB}\) to determine the continuity of the projected maps \(f_{\rmks}\) and \(f_{\Delta}\) on \(\vartheta_{\rms}\).

From the expressions in \eqref{eqn:x3_prime} and \eqref{eqn:y3_prime}, we calculate the action of the completed return map in logarithmic coordinates:
\begin{equation}\label{eqn:X3_prime}
  \begin{aligned}
    X_{3}^{\prime}
    &\coloneq\log x_{3}^{\prime}\\
    &=\log E_{\rmB\rmX}+\log h+\log\cos\bar{\theta}_{\rmA}\left(\theta_{3}^{\rmC}\left(T_{\rmC}\right)\right)\\
    &\qquad\qquad+\dfrac{c_{\rmA\rmX}}{e_{\rmA}}\left(\log E_{\rmA}\left(\theta_{3}^{\rmC}\left(T_{\rmC}\right)\right)-\int_{0}^{T_{\rmC}}g_{\rmB}(\theta_{3}^{\rmC}\left(\tau\right))\rmd\tau\right)
  \end{aligned}
\end{equation}
and
\begin{equation}\label{eqn:Y3_prime}
  \begin{aligned}
    Y_{3}^{\prime}
    &\coloneq\log y_{3}^{\prime}\\
    &=\log E_{\rmB\rmY}+\log h+\log\sin\bar{\theta}_{\rmA}\left(\theta_{3}^{\rmC}\left(T_{\rmC}\right)\right)\\
    &\qquad\qquad+\dfrac{c_{\rmA\rmY}}{e_{\rmA}}\left(\log E_{\rmA}\left(\theta_{3}^{\rmC}\left(T_{\rmC}\right)\right)-\int_{0}^{T_{\rmC}}g_{\rmB}(\theta_{3}^{\rmC}\left(\tau\right))\rmd\tau\right)
  \end{aligned}
\end{equation}
These expressions allow us to define a map \(\widehat{M}\colon\mcl{D}\to\mcl{D}\):
\begin{equation*}
  \widehat{M}\left(X_{3},Y_{3}\right)=\left(X_{3}^{\prime},Y_{3}^{\prime}\right).
\end{equation*}
The domain of this map is
\begin{equation*}
  \mcl{D}=\left\{\left(X_{3},Y_{3}\right)\in\R^{2}\mid X_{3},Y_{3}<\log h\right\}\subsetneq\negativeR{2}.
\end{equation*}

\subsection{Projection of the completed return map}\label{sec:anal:ssec:proj}

Just as the projected maps \(f_{\rmks}\) and \(f_{\Delta}\) are defined from the transition matrices of the return maps \eqref{eqn:ks_ret_map} and \eqref{eqn:delta_ret_map}, we can define a projected map from the map \(\widehat{M}\). Let \(k>0\) be a positive real number. We define the line \(T^{(k)}\) in the negative quadrant of \(\R^{2}\) by
\begin{equation*}
  T^{(k)}\coloneq\left\{\left(X_{3},Y_{3}\right)\in\mcl{D}\mid X_{3}+Y_{3}=-k\right\}.
\end{equation*}
Note that in the definition given in \eqref{eqn:S}, the set \(S=\left(-1,0\right)\), which is the domain of both \(f_{\rmks}\) and \(f_{\Delta}\), is the straight line from the point \(\left(-1,0\right)\) to the point \(\left(0,-1\right)\). The line \(T^{(k)}\) is instead a subset of \(\mcl{D}\), and so does not extend to the coordinate axes, as all points \(\left(X_{3},Y_{3}\right)\) on the line \(T^{(k)}\) must satisfy \(X_{3}<\log h\) and \(Y_{3}<\log h\). We can identify this line with the open interval \(\left(-k-\log h,\log h\right)\).

\begin{figure}
  \centering
  \includegraphics[width=0.5\linewidth]{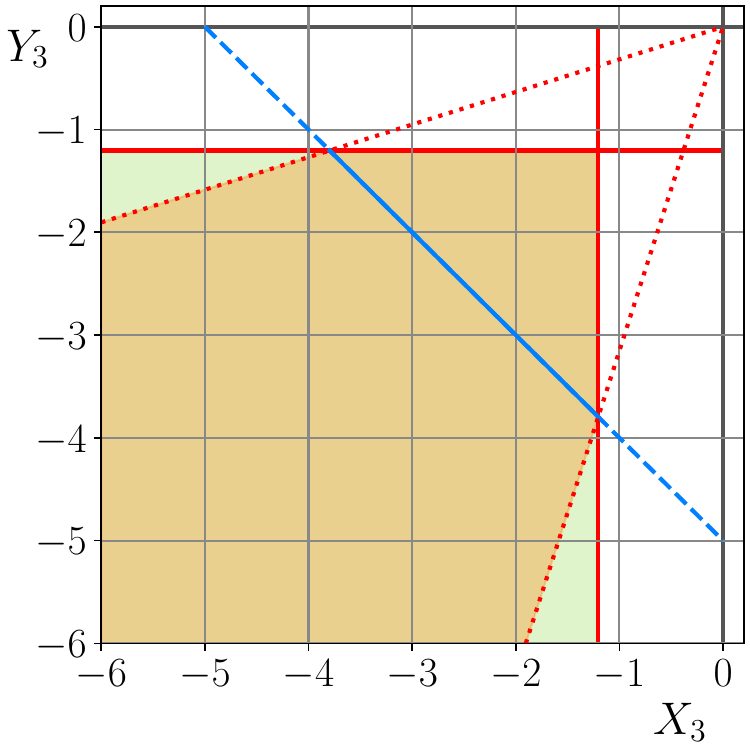}
  \caption{A representation of the domain of the map \(\widehat{M}\) and of the projection \(\Pi^{(k)}\), and the set \(T^{(k)}\). We show a subset of \(\negativeR{2}\). The solid red lines are \(X_{3}=\log h\) and \(Y_{3}=\log h\). The region bounded by these lines, shaded green and orange, is the domain \(\mcl{D}\) of \(\widehat{M}\). The solid blue line is the set \(T^{(k)}\), and the set \(\mcl{D}^{(k)}\) of all points in \(\mcl{D}\) that project onto a point in \(T^{(k)}\) are shaded orange. This region is bounded by the dotted red lines.}
  \label{fig:projection_domain}
\end{figure}

We can rescale this set by \(k\) to define the open interval
\begin{equation*}
  S^{(k)}=\left(-1-\frac{\log h}{k},\frac{\log h}{k}\right).
\end{equation*}
We then have \(\lim_{k\to\infty}S^{(k)}=S\). To generalise the definition of the projection \(\Pi\) in \eqref{eqn:proj}, we define a projection \(\Pi^{(k)}\colon\mcl{D}^{(k)}\to S^{(k)}\) by
\begin{equation*}
  \Pi^{(k)}\colon\left(X_{3},Y_{3}\right)\mapsto\dfrac{-X_{3}}{X_{3}+Y_{3}}.
\end{equation*}
The domain of this projection is 
\begin{equation*}
  \mcl{D}^{(k)}=\left\{\left(X_{3},Y_{3}\right)\in\mcl{D}\biggm\vert Y_{3}<\dfrac{-\log h}{k+\log h}X_{3}\ \textrm{and}\ X_{3}<\dfrac{-\log h}{k+\log h}Y_{3}\right\}
\end{equation*}
This set is the set of all points in \(\negativeR{2}\) that project onto a point in \(T^{(k)}\) that is in the domain of definition of \(\widehat{M}\). We show a visual representation of this domain in \Fref{fig:projection_domain}.

\begin{figure}[!t]
  \begin{subfigure}{0.5\linewidth}
    \includegraphics[width=\linewidth]{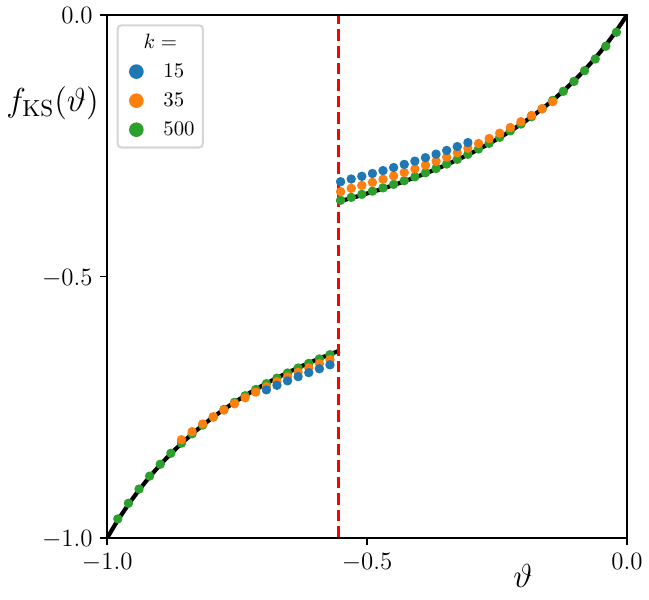}
    \vspace{-3mm}
    \caption{}
    \label{fig:proj_map_approx:ks}
\end{subfigure}
  \begin{subfigure}{0.5\linewidth}
    \includegraphics[width=\linewidth]{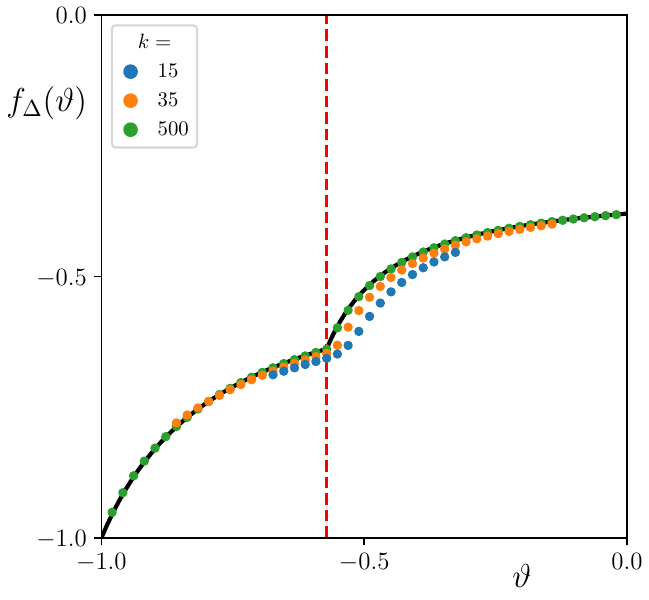}
    \vspace{-3mm}
    \caption{}
    \label{fig:proj_map_approx:delta}
\end{subfigure}
  \caption{Numerical approximations of the projected maps (\subref{fig:proj_map_approx:ks}) \(f_{\mathrm{KS}}^{(k)}\) and (\subref{fig:proj_map_approx:delta}) \(f_{\Delta}^{(k)}\), for various values of \(k\), as indicated. The projected maps \(f_{\mathrm{KS}}\) and \(f_{\Delta}\) are plotted as solid black lines. As \(k\) increases, \(f_{\mathrm{KS}}^{(k)}\) and \(f_{\Delta}^{(k)}\) converge to \(f_{\mathrm{KS}}\) and \(f_{\Delta}\), respectively. Note that, for larger \(k\), the domains of the maps widen. See \Fref{fig:proj_map_approx_sm} for more detailed plots near the switching manifold. See \ref{app:numerics:ssec:projected_k} for an explanation of how these approximations are calculated. For visual clarity, we only plot the components of the projected map within their domain of definition, and we do not shade these domains. The dashed red line corresponds to \(\vartheta_{\rms}\). Parameter values for these Figures~are given in Tables \ref{tbl:parms_fk_ks} and \ref{tbl:parms_fk_delta}.}
  \label{fig:proj_map_approx}
\end{figure}

\begin{figure}[!t]
  \begin{subfigure}{0.5\linewidth}
    \includegraphics[width=\linewidth]{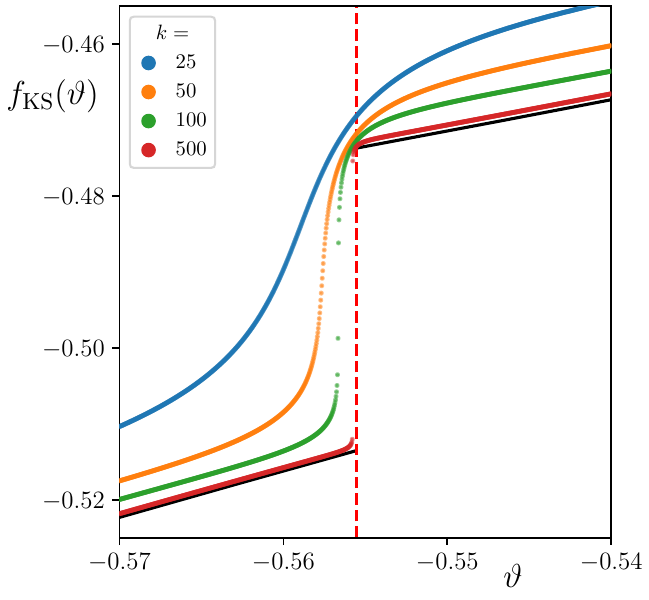}
    \vspace{-3mm}
    \caption{}
    \label{fig:proj_map_approx_sm:ks}
  \end{subfigure}
  \begin{subfigure}{0.5\linewidth}
    \includegraphics[width=\linewidth]{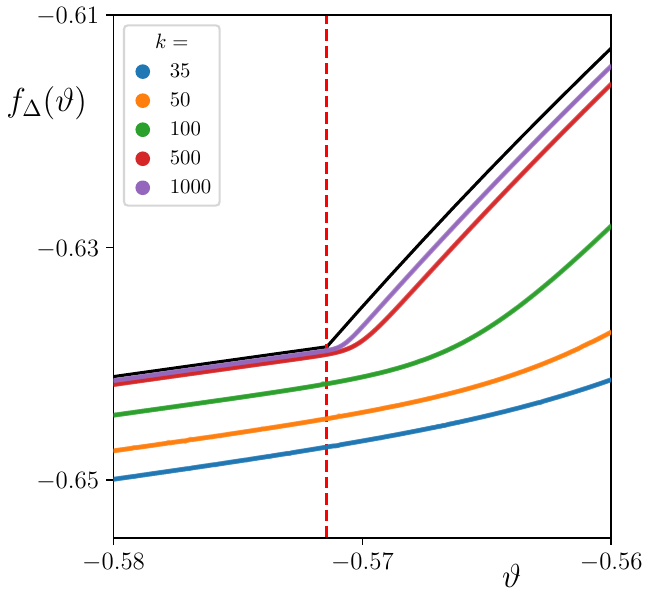}
    \vspace{-3mm}
    \caption{}
    \label{fig:proj_map_approx_sm:delta}
  \end{subfigure}
  \caption{Numerical approximations of the projected maps (\subref{fig:proj_map_approx_sm:ks}) \(f_{\mathrm{KS}}^{(k)}\) and (\subref{fig:proj_map_approx_sm:delta}) \(f_{\Delta}^{(k)}\) near the switching manifold \(\vartheta_{\rms}\) for various values of \(k\), as indicated. The projected maps \(f_{\mathrm{KS}}\) and \(f_{\Delta}\) are plotted as solid black lines. As \(k\) increases, we see that the functions converge to \(f_{\mathrm{KS}}\) and \(f_{\Delta}\). More importantly, however, we see a discontinuity emerge in the case of \(f_{\mathrm{KS}}^{(k)}\). See \ref{app:numerics:ssec:projected_k} for an explanation of how these approximations are calculated. For visual clarity, we only plot the components of the projected map within their domain of definition, and we do not shade these domains. The dashed red line corresponds to \(\vartheta_{\rms}\). Parameter values for these Figures~are given in Tables \ref{tbl:parms_fk_ks} and \ref{tbl:parms_fk_delta}.}
  \label{fig:proj_map_approx_sm}
\end{figure}

With this, we can define a projected map \(f^{(k)}\colon S^{(k)}\to S^{(k)}\):
\begin{equation*}
  f^{(k)}\colon\vartheta\mapsto\Pi^{(k)}\left(\widehat{M}\left(\vartheta,-k-\vartheta\right)\right).
\end{equation*}
The mapping \(f^{(k)}\) is continuous. Since the map \(\widehat{M}\) is not linear, the value of \(f^{(k)}\left(\vartheta\right)\) depends on the value of \(k\); that is, if \(\ell>k\), then for any point \(\vartheta\in S^{(k)}\subsetneq S^{(\ell)}\), \(f^{(\ell)}\left(\vartheta\right)\neq f^{(k)}\left(\vartheta\right)\).

We show numerical approximations of this projected map for both networks in \Fref{fig:proj_map_approx}. We observe that, as \(k\to\infty\), the function \(f^{(k)}\) converges to the projected maps \(f_{\rmks}\) and \(f_{\Delta}\), and we conjecture that this convergence is the case, though we do not prove that here. In \Fref{fig:proj_map_approx:ks}, the map \(f_{\rmks}^{(k)}\) appears to be discontinuous at \(\vartheta_{\rms}\), though only because of the discretisation of the numerical calculations. In \Fref{fig:proj_map_approx_sm}, we show the same calculations in a small neighbourhood of the switching manifold \(\vartheta_{\rms}\). We see that the projected maps \(f^{(k)}\) are continuous at \(\vartheta=\vartheta_{\rms}\) for both the Kirk--Silber and \(\Delta\)-clique networks, though a discontinuity emerges in the limit \(k\to\infty\) for the Kirk--Silber network. Similarly in \Fref{fig:theta_3_C_out_v_theta_3_C_in}(\subref{fig:theta_3_C_out_v_theta_3_C_in:kirk_silber_net}), we also observed a discontinuity emerge in the function \(\theta_{3}^{\outvar{C}}\left(\theta_{3}^{\invar{C}}\right)\) in the limit \(x_{1}^{\invar{C}}\to0\).

\subsection{Calculating the completed return map}\label{sec:anal:ssec:calc}

We now calculate expressions of the completed return map that are valid for trajectories that start near the switching curve \(\Sigma_{\rms}\), and so that do not strike \(\poinin[B]{C}\) near \(X\) or \(Y\), and that are sufficiently close to the heteroclinic network, so that they strike \(\poinout[A]{C}\) near \(X\) or \(Y\). We present full details for the case of the Kirk--Silber network for trajectories that leave \(C\) near \(X\), and provide more abridged arguments for other cases, as the argument follows closely.

\subsubsection{The Kirk--Silber network for trajectories that leave \(C\) near \(X\)}

We begin with the Kirk--Silber network. Consider a point
\begin{equation*}
  \left(x_{3}^{\invar{B}},y_{3}^{\invar{B}}\right)\in\poinin[A]{B}.
\end{equation*}
Applying the basic map defined by the composition \(\widehat{\Psi}_{\rmB\rmC}\widehat{\psi}_{\rmB}\) to this point, we calculate
\begin{equation*}
  \left(\theta_{3}^{\invar{C}},x_{1}^{\invar{C}}\right)=\left(\widehat{\Psi}_{\rmB\rmC}\widehat{\psi}_{\rmB}\right)\left(x_{3}^{\invar{B}},y_{3}^{\invar{B}}\right).
\end{equation*}

Recall that \(\theta_{3}^{\rmP}\) is the value of the angle \(\theta_{3}\) at \(P\). We also define \(\theta_{3}^{\rmX}=0\) as the value of the angle \(\theta_{3}\) at \(X\), and likewise for \(\theta_{3}^{\rmY}=\frac{\pi}{2}\).

To understand the behaviour of trajectories that are close to the switching curve \(\Sigma_{\rms}\) defined by \(x_{3}^{e_{\rmB\rmY}}=y_{3}^{e_{\rmB\rmX}}\), we consider a trajectory beginning at a point \(\left(x_{3}^{\invar{B}},y_{3}^{\invar{B}}\right)\) close to the curve \(\Sigma_{\rms}\). Therefore, this point lies inside the cusp \(\Gamma_{\rmc}\) that is excluded from the domain of the return maps in \eqref{eqn:ks_ret_map} and \eqref{eqn:delta_ret_map}. We assume that this trajectory strikes \(\poinin[B]{C}\) in a small neighbourhood of the equilibrium \(P\), and close to the heteroclinic orbit \(P\to X\). This condition is \(\theta_{3}^{\rmP}-h<\theta_{3}^{\invar{C}}<\theta_{3}^{\rmP}\) for some \(h\ll 1\). For \(\theta_{3}^{\invar{C}}<\theta_{3}^{\rmP}\), we require \(x_{3}^{e_{\rmB\rmY}}>y_{3}^{e_{\rmB\rmX}}\).

We also assume that the trajectory strikes \(\poinout[\mrm{A}]{\mrm{C}}\) in a small neighbourhood of the equilibrium \(X\). This condition is \(\theta_{3}^{\rmX}<\theta_{3}^{\outvar{C}}<\theta_{3}^{\rmX}+h\), and requires \(T_{\rmB}\) to be sufficiently large. We quantify this requirement explicitly in the inequality \eqref{eqn:TB_lower_lim} below.

The flow linearised near \(P\) is given by
\begin{align*}
  \dot{x}_{1}&=g_{\rmA}\left(\theta_{3}^{\rmP}\right)x_{1},\\
  \dot{x}_{2}&=-g_{\rmB}\left(\theta_{3}^{\rmP}\right)x_{2},\\
  \dot{\theta}_{3}&=\lambda_{\rmP}\left(\theta_{3}-\theta_{3}^{\rmP}\right),
\end{align*}
for some real number \(\lambda_{\rmP}>0\).

Our assumption that \(\theta_{3}^{\outvar{C}}<\theta_{3}^{\rmX}+h\) implies that there is some time \(T_{\rmC}^{(\rmP)}\) such that \(\theta_{3}^{\rmC}\left(T_{\rmC}^{(\rmP)}\right)=\theta_{3}^{\rmP}-h\). Solving \(\dot{\theta_{3}}\) with initial condition \(\theta_{3}\left(0\right)=\theta_{3}^{\invar{C}}\), we find
\begin{equation*}
  T_{\rmC}^{(P)}=\frac{-1}{\lambda_{\rmP}}\log\left(\frac{\theta_{3}^{\rmP}-\theta_{3}^{\invar{C}}}{h}\right).
\end{equation*}
We solve both \(\dot{x}_{1}\) and \(\dot{x}_{2}\) with initial conditions
\begin{equation*}
  x_{1}(0)=x_{1}^{\invar{C}}=E_{\rmC}\left(\theta_{3}^{\outvar{B}}\right)he^{-c_{\rmB}T_{\rmB}}
\end{equation*}
and \(x_{2}\left(0\right)=h\), where we recall that \(E_{\rmC}\left(\theta_{3}^{\outvar{B}}\right)\) is the angle at which a trajectory leaving \(\poinout[\rmC]{\rmB}\) then strikes \(\poinin[\rmB]{\rmC}\). We find that
\begin{equation*}
  x_{1}\left(T_{\rmC}^{(\rmP)}\right)=E_{\rmC}\left(\theta_{3}^{out,B}\right)he^{-c_{\rmB}T_{\rmB}}\left(\frac{\theta_{3}^{\rmP}-\theta_{3}^{\invar{C}}}{h}\right)^{-\frac{g_{\rmA}\left(\theta_{3}^{\rmP}\right)}{\lambda_{\rmP}}}
\end{equation*}
and
\begin{equation*}
  x_{2}\left(T_{\rmC}^{(\rmP)}\right)=h\left(\frac{\theta_{3}^{\rmP}-\theta_{3}^{\invar{C}}}{h}\right)^{\frac{g_{\rmB}\left(\theta_{3}^{\rmP}\right)}{\lambda_{\rmP}}}
\end{equation*}

There is then a jump from \(P\) to \(X\) as the trajectory follows the heteroclinic orbit \(P\to X\). In this transition, both \(x_{1}\) and \(x_{2}\) are scaled by an \(h\)-dependent \(\order{1}\) amount, \(E_{\rmP\rmX}^{(1)}\) and \(E_{\rmP\rmX}^{(2)}\), respectively. Therefore, at \(\theta_{3}=\theta_{3}^{\rmX}+h\), the value of \(x_{1}\) is
\begin{equation}\label{eqn:x1_X}
  E_{\rmP\rmX}^{(1)}E_{\rmC}\left(\theta_{3}^{out,B}\right)he^{-c_{\rmB}T_{\rmB}}\left(\frac{\theta_{3}^{\rmP}-\theta_{3}^{\invar{C}}}{h}\right)^{-\frac{g_{\rmA}\left(\theta_{3}^{\rmP}\right)}{\lambda_{\rmP}}}
\end{equation}
and the value of \(x_{2}\) is
\begin{equation*}
  E_{\rmP\rmX}^{(2)}h\left(\frac{\theta_{3}^{\rmP}-\theta_{3}^{\invar{C}}}{h}\right)^{\frac{g_{\rmB}\left(\theta_{3}^{\rmP}\right)}{\lambda_{\rmP}}}.
\end{equation*}
Our assumption that the trajectory strikes \(\poinout[\mrm{A}]{\mrm{C}}\) in a small neighbourhood of \(X\) means that we require the value of \(x_{1}\) in \eqref{eqn:x1_X} to be strictly less than \(h\). Solving this inequality for \(T_{\rmB}\) shows that we require
\begin{equation}\label{eqn:TB_lower_lim}
  T_{\rmB}>\frac{1}{c_{\rmB}}\left(\log\left(E_{\rmP\rmX}^{(1)}E_{\rmC}\left(\theta_{3}^{\outvar{B}}\right)\right)-\frac{g_{\rmA}\left(\theta_{3}^{\rmP}\right)}{\lambda_{\rmP}}\log\left(\frac{\theta_{3}^{\rmP}-\theta_{3}^{\invar{C}}}{h}\right)\right).
\end{equation}
Therefore, our assumption that the trajectory strikes \(\poinin[\mrm{B}]{\mrm{C}}\) near \(P\) and then strikes \(\poinout[\mrm{A}]{\mrm{C}}\) near \(X\) requires \(T_{\rmB}\to\infty\) as \(\theta_{3}^{\invar{C}}\to\theta_{3}^{\rmP}\).

Near \(X\), the linearised flow is
\begin{align*}
  \dot{x}_{1}&=e_{\rmX\rmA}x_{1},\\
  \dot{x}_{2}&=-c_{\rmX\rmB}x_{2},\\
  \dot{\theta}_{3}&=-c_{\rmXY}\theta_{3}.
\end{align*}
We solve for the residence time \(T_{\rmC}^{(\rmX)}\) near \(X\) by solving \(x_{1}\left(T_{\rmC}^{(\rmX)}\right)=h\) with the initial condition of \(x_{1}\) given in \eqref{eqn:x1_X}:
\begin{equation*}
  T_{\rmC}^{(\rmX)}=-\frac{1}{e_{\rmX\rmA}}\left(\log\left(E_{\rmP\rmX}^{(1)}E_{\rmC}\left(\theta_{3}^{\outvar{B}}\right)\right)-\frac{g_{\rmA}\left(\theta_{3}^{\rmP}\right)}{\lambda_{\rmP}}\log\left(\frac{\theta_{3}^{\rmP}-\theta_{3}^{\invar{C}}}{h}\right)-c_{\rmB}T_{\rmB}\right).
\end{equation*}
Therefore, the value of \(x_{2}\) on \(\poinout[\rmA]{C}\) is 
\begin{equation}\label{eqn:KS_X_x2_out_C}
  x_{2}^{\outvar{C}}=E_{\rmP\rmX}^{(2)}h\left(\frac{\theta_{3}^{\rmP}-\theta_{3}^{\invar{C}}}{h}\right)^{\frac{g_{\rmB}\left(\theta_{3}^{\rmP}\right)}{\lambda_{\rmP}}}\left(E_{\rmP\rmX}^{(1)}E_{\rmC}\left(\theta_{3}^{\outvar{B}}\right)e^{-C_{\rmB}T_{\rmB}}\left(\frac{\theta_{3}^{\rmP}-\theta_{3}^{\invar{C}}}{h}\right)^{-\frac{g_{\rmA}\left(\theta_{3}^{\rmP}\right)}{\lambda_{\rmP}}}\right)^{\frac{c_{\rmX\rmB}}{e_{\rmX\rmA}}}
\end{equation}
and the value of \(\theta_{3}^{\outvar{C}}\equiv\theta_{3}^{\rmC}\left(T_{\rmC}\right)\) is
\begin{equation}\label{eqn:KS_X_theta_3_out_C}
  \theta_{3}^{\outvar{C}}=h\left(E_{\rmP\rmX}^{(1)}E_{\rmC}\left(\theta_{3}^{\outvar{B}}\right)e^{-c_{\rmB}T_{\rmB}}\left(\frac{\theta_{3}^{\rmP}-\theta_{3}^{\invar{C}}}{h}\right)^{-\frac{g_{\rmA}\left(\theta_{3}^{\rmP}\right)}{\lambda_{\rmP}}}\right)^{\frac{c_{\rmX\rmY}}{e_{\rmX\rmA}}}.
\end{equation}

\subsubsection{The Kirk--Silber network for trajectories that leave \(C\) near \(Y\)}

A similar analysis can be completed in the case of the Kirk--Silber network when a trajectory strikes \(\poinin[B]{C}\) in a small neighbourhood of the equilibrium \(P\), and close to the heteroclinic orbit \(P\to Y\), and that strikes \(\poinout[A]{C}\) near \(Y\). These conditions are \(\theta_{3}^{\rmP}<\theta_{3}^{\invar{C}}<\theta_{3}^{\rmP}+h\) and \(\theta_{3}^{\rmY}-h<\theta_{3}^{\outvar{C}}<\theta_{3}^{\rmY}\), respectively. In this case, the value of \(x_{2}\) on \(\poinout[\rmA]{C}\) is 
\begin{equation}\label{eqn:KS_Y_x2_out_C}
  x_{2}^{\outvar{C}}=E_{\rmP\rmY}^{(2)}h\left(\frac{\theta_{3}^{\rmP}-\theta_{3}^{\invar{C}}}{h}\right)^{\frac{g_{\rmB}\left(\theta_{3}^{\rmP}\right)}{\lambda_{\rmP}}}\left(E_{\rmP\rmY}^{(1)}E_{\rmC}\left(\theta_{3}^{\outvar{B}}\right)e^{-C_{\rmB}T_{\rmB}}\left(\frac{\theta_{3}^{\rmP}-\theta_{3}^{\invar{C}}}{h}\right)^{-\frac{g_{\rmA}\left(\theta_{3}^{\rmP}\right)}{\lambda_{\rmP}}}\right)^{\frac{c_{\rmY\rmB}}{e_{\rmY\rmA}}}
\end{equation}
and the value of \(\varphi_{3}^{\outvar{C}}\), where we recall \(\varphi_{3}=\frac{\pi}{2}-\theta_{3}\), is
\begin{equation}\label{eqn:KS_Y_theta_3_out_C}
  \theta_{3}^{\outvar{C}}=h\left(E_{\rmP\rmY}^{(1)}E_{\rmC}\left(\theta_{3}^{\outvar{B}}\right)e^{-c_{\rmB}T_{\rmB}}\left(\frac{\theta_{3}^{\rmP}-\theta_{3}^{\invar{C}}}{h}\right)^{-\frac{g_{\rmA}\left(\theta_{3}^{\rmP}\right)}{\lambda_{\rmP}}}\right)^{\frac{c_{\rmY\rmX}}{e_{\rmY\rmA}}}.
\end{equation}

\subsubsection{The \(\Delta\)-clique network}

In the case of the \(\Delta\)-clique network, there is no equilibrium \(P\). Any trajectory that begins in the cusp that is excluded from the cross-section \(\poinin[\mrm{A}]{\mrm{B}}\) strikes \(\poinin[\mrm{B}]{\mrm{C}}\) along the heteroclinic connection \(X\to Y\), and not in a small neighbourhood of \(X\) or \(Y\). Therefore, the initial conditions on \(\poinin[\mrm{A}]{\mrm{C}}\) are \(x_{1}=E_{\rmC}\left(\theta_{3}^{\outvar{B}}\right)he^{-c_{\rmB}T_{\rmB}}\), \(x_{2}=h\), and \(\theta_{3}^{\rmX}+h<\theta_{3}^{\invar{C}}<\theta_{3}^{\rmY}-h\).

The transition of the trajectory along the heteroclinic orbit \(X\to Y\) rescales \(x_{1}\) and \(x_{2}\) by \(h\)-dependent \(\order{1}\) amounts, \(E_{\rmX\rmY}^{(1)}\) and \(E_{\rmX\rmY}^{(2)}\), respectively. Near \(Y\), the linearised flow is
\begin{align*}
  \dot{x}_{1}&=e_{\rmY\rmA}x_{1},\\
  \dot{x}_{2}&=-c_{\rmY\rmB}x_{2},\\
  \dot{\varphi}_{3}&=-c_{\rmY\rmX}\varphi{3}.
\end{align*}
Solving for \(x_{1}\left(T_{\rmC}^{(\rmY)}\right)=h\) and then integrating the above linear system gives \(x_{2}^{\outvar{C}}\)
\begin{equation}\label{eqn:Delta_x2_out_C}
  x_{2}^{\outvar{C}}=E_{\rmX\rmY}^{(2)}h\left(E_{\rmX\rmY}^{(1)}E_{\rmC}\left(\theta_{3}^{\outvar{B}}\right)e^{-c_{\rmB}T_{\rmB}}\right)^{\frac{c_{\rmY\rmB}}{e_{\rmY\rmA}}}
\end{equation}
and the value of \(\theta_{3}^{\outvar{C}}\) is
\begin{equation}\label{eqn:Delta_theta_3_out_C}
  \theta_{3}^{\outvar{C}}=h\left(E_{\rmX\rmY}^{(1)}E_{\rmC}\left(\theta_{3}^{\outvar{B}}\right)e^{-c_{\rmB}T_{\rmB}}\right)^{\frac{c_{\rmY\rmX}}{e_{\rmY\rmA}}}.
\end{equation}

\subsection{Applying the projection}\label{sec:anal:ssec:expl}

We now apply the values for \(\theta_{3}^{\outvar{C}}\) and \(x_{2}^{\outvar{C}}\) to the expressions for \(X_{3}^{\prime}\) and \(Y_{3}^{\prime}\) in \eqref{eqn:X3_prime} and \eqref{eqn:Y3_prime}, respectively.

\subsubsection{The Kirk--Silber network for \(\theta_{3}^{\invar{C}}<\theta_{3}^{\rmP}\)}

For the Kirk--Silber network with \(\theta_{3}^{\invar{C}}<\theta_{3}^{\rmP}\) and \(T_{\rmC}\) sufficiently large, we know that \(\theta_{3}^{\outvar{C}}\equiv\theta_{3}^{\rmC}\left(T_{\rmC}\right)<\theta_{3}^{\rmX}+h\ll 1\); that is, the trajectory leaves a small neighbourhood of \(C\) near \(X\).

We Taylor expand \(\bar{\theta}_{\rmA}\) about the fixed point \(\bar{\theta}_{\rmA}\left(0\right)=0\) to derive the leading order expansion
\begin{equation*}
  \bar{\theta}_{\rmA}\left(\theta_{3}^{\rmC}\left(T_{\rmC}\right)\right)=K_{1}\theta_{3}^{\rmC}\left(T_{\rmC}\right)+\order{\theta_{3}^{\rmC}\left(T_{\rmC}\right)^{2}},
\end{equation*}
for some \(\order{1}\) constant \(K_{1}\). We can then apply a leading order Taylor expansion to the trigonometric terms in \eqref{eqn:X3_prime} and \eqref{eqn:Y3_prime}. The resulting expressions are rather long and cumbersome, and so we do not present them here, but they are the sum of terms involving \(T_{\rmB}\), \(\log\left(\frac{\theta_{3}^{\rmP}-\theta_{3}^{\invar{C}}}{h}\right)\), and the logarithm of constant terms. Since we have assumed that \(T_{\rmB}\gg1\) and \(\theta_{3}^{\rmP}-\theta_{3}^{\invar{C}}\ll h\), we can ignore the logarithm of the constant terms as they do not appear in asymptotically significant expressions in the limit \(T_{\rmB}\to\infty\) or \(\theta_{3}^{\invar{C}}\to\theta_{3}^{\rmP}\). We then substitute the expressions in \eqref{eqn:KS_X_x2_out_C} and \eqref{eqn:KS_X_theta_3_out_C} into \eqref{eqn:X3_prime} and \eqref{eqn:Y3_prime}, and the resulting expressions into the projection \(\Pi^{(k)}\), and rearrange to derive
\begin{equation}\label{eqn:ks_upper_lim_projection}
  \begin{aligned}
    &\qquad\qquad\qquad\qquad\Pi^{(k)}\left(X_{3}^{\prime},Y_{3}^{\prime}\right)=\\
    &\frac{
      -c_{\rmA\rmX}\left(-\frac{e_{\rmX\rmA}}{c_{\rmB}\lambda_{\rmP}}\left(g_{\rmB}(\theta_{3}^{\rmP})-\frac{g_{\rmA}(\theta_{3}^{\rmP})c_{\rmX\rmB}}{e_{\rmX\rmA}}\right)\frac{\log\delta}{T_{\rmB}}+c_{\rmX\rmB}\right)
      }{
        (c_{\rmA\rmX}+c_{\rmA\rmY})\left(-\frac{e_{\rmX\rmA}}{c_{\rmB}\lambda_{\rmP}}\left(g_{\rmB}(\theta_{3}^{\rmP})-\frac{g_{\rmA}(\theta_{3}^{\rmP})c_{\rmX\rmB}}{e_{\rmX\rmA}}\right)\frac{\log\delta}{T_{\rmB}}+c_{\rmX\rmB}\right)+e_{\rmA}\left(\frac{g_{\rmA}(\theta_{3}^{\rmP})c_{\rmX\rmY}}{c_{\rmB}\lambda_{\rmP}}\frac{\log\delta}{T_{\rmB}}+c_{\rmX\rmY}\right)
        }.
  \end{aligned}
\end{equation}
In this expression, we have written
\begin{equation*}
  \delta=\frac{\theta_{3}^{\rmP}-\theta_{3}^{\invar{C}}}{h}.
\end{equation*}

To evaluate the limit in \eqref{eqn:ks_lim_right}, of the projected map \(f_{\rmks}\) as \(\vartheta\searrow\vartheta_{\rms}\), we must follow the definition of the projected map outlined in \Sref{sec:problem:ssec:limit}, and calculate the projection of the action of the transition matrix \(M_{\rmX}\) on points on the line \(Y_{3}=\log \left(1-\epsilon\right)+\frac{e_{\rmB\rmY}}{e_{\rmB\rmX}}X_{3}\). Recall from \eqref{eqn:gamma_sets} in \Sref{sec:problem:ssec:ret_maps} that \(\epsilon\) was a positive constant used in the definition of the excluded cusp. The limit \(\vartheta\searrow\vartheta_{\rms}\) is then the limit \(X_{3}\to-\infty\).

On \(\poinin[A]{B}\), curves of constant angle \(\theta_{3}^{\invar{B}}\) have the form \(x_{3}=ay_{3}\). Therefore, in logarithmic coordinates, curves of constant \(\theta_{3}^{\invar{B}}\) in \(\mcl{D}\) are \(X_{3}=\log a+Y_{3}\). We see from the action of the local map \(\widehat{\psi}_{\rmB}\) in \eqref{eqn:theta_B_out} that curves of constant angle \(\theta_{3}^{\outvar{B}}\) in \(\mcl{D}\) have the form \(X_{3}=\log a+\frac{e_{\rmB\rmX}}{e_{\rmB\rmY}}Y_{3}\).

For any \(0>\log a>\log\left(1-\epsilon\right)\), the line \(X_{3}=\log a+\frac{e_{\rmB\rmX}}{e_{\rmB\rmY}}Y_{3}\) lies inside the excluded cusp. Moreover, the expression in \eqref{eqn:ks_upper_lim_projection} is valid, and can be evaluated at any point on the line \(X_{3}=\log a+\frac{e_{\rmB\rmX}}{e_{\rmB\rmY}}Y_{3}\), whereas the transition matrix \(M_{\rmX}\) cannot. For any point along this line, \(\theta_{3}^{\invar{C}}\), and therefore also \(\delta\), is fixed, but as \(X_{3}\to-\infty\), we see from \eqref{eqn:TB} that \(T_{\rmB}\to\infty\) and thus \(\frac{\log\delta}{T_{\rmB}}\to 0\). Therefore, the expression in \eqref{eqn:ks_upper_lim_projection} simplifies to
\begin{equation*}
  \Pi^{(k)}\left(X_{3}^{\prime},Y_{3}^{\prime}\right)=\frac{-c_{\rmA\rmX}c_{\rmX\rmB}}{(c_{\rmA\rmX}+c_{\rmA\rmY})c_{\rmX\rmB}+e_{\rmA}c_{\rmX\rmY}},
\end{equation*}
which is the value of the limit \(\lim_{\vartheta\searrow\vartheta_{\rms}}f_{\rmks}\left(\vartheta\right)\) in \eqref{eqn:ks_lim_right}. This limit applies for all \(0>\log a>\log\left(1-\epsilon\right)\), and so \(\theta_{3}^{\rmP}-\theta_{3}^{\invar{C}}\) can be chosen to be arbitrarily close to \(0\), and thus the trajectory is arbitrarily  close to the switching manifold.

\subsubsection{The Kirk--Silber network for \(\theta_{3}^{\rmP}<\theta_{3}^{\invar{C}}\)}

In the case of the Kirk--Silber network with \(\theta_{3}^{\rmP}<\theta_{3}^{\invar{C}}<\theta_{3}^{\rmP}+h\), substituting the expressions in \eqref{eqn:KS_Y_x2_out_C} and \eqref{eqn:KS_Y_theta_3_out_C} into \eqref{eqn:X3_prime} and \eqref{eqn:Y3_prime}, disregarding the logarithm of \(\order{1}\) constant terms, which do not affect the dynamics in the limit \(T_{\rmB}\to\infty\), and applying the projection, gives
\begin{equation}\label{eqn:ks_lower_lim_projection}
  \begin{aligned}
    &\qquad\qquad\qquad\qquad\Pi^{(k)}\left(X_{3}^{\prime},Y_{3}^{\prime}\right)=\\
    &\frac{
      -c_{\rmA\rmX}\left(-\frac{e_{\rmY\rmA}}{c_{\rmB}\lambda_{\rmP}}\left(g_{\rmB}(\theta_{3}^{\rmP})-\frac{g_{\rmA}(\theta_{3}^{\rmP})c_{\rmY\rmB}}{e_{\rmY\rmA}}\right)\frac{\log\delta}{T_{\rmB}}+c_{\rmY\rmB}\right)-e_{\rmA}\left(\frac{g_{\rmA}(\theta_{3}^{\rmP})c_{\rmY\rmX}}{c_{\rmB}\lambda_{\rmP}}\frac{\log\delta}{T_{\rmB}}+c_{\rmY\rmX}\right)
      }{
        (c_{\rmA\rmX}+c_{\rmA\rmY})\left(-\frac{e_{\rmY\rmA}}{c_{\rmB}\lambda_{\rmP}}\left(g_{\rmB}(\theta_{3}^{\rmP})-\frac{g_{\rmA}(\theta_{3}^{\rmP})c_{\rmY\rmB}}{e_{\rmY\rmA}}\right)\frac{\log\delta}{T_{\rmB}}+c_{\rmY\rmB}\right)+e_{\rmA}\left(\frac{g_{\rmA}(\theta_{3}^{\rmP})c_{\rmY\rmX}}{c_{\rmB}\lambda_{\rmP}}\frac{\log\delta}{T_{\rmB}}+c_{\rmY\rmX}\right)    
        }.
  \end{aligned}
\end{equation}
We note that the additional term in the numerator compared to \eqref{eqn:ks_upper_lim_projection} is a result of the fact that the leading-order Taylor expansion of \(\log\cos\bar{\theta}_{\rmA}\left(\theta\right)\) around the point \(\frac{\pi}{2}\) is \(\log\theta\), as opposed to \(\log1\) when the Taylor expansion is evaluated around \(0\).

Again, for any \(0<\log a<\log\left(\frac{1}{1-\epsilon}\right)\), the line \(X_{3}=\log a+\frac{e_{\rmB\rmX}}{e_{\rmB\rmY}}Y_{3}\) lies inside the excluded cusp, and the expression in \eqref{eqn:ks_lower_lim_projection} is valid, and the expression can be evaluated at any point on that line. In the limit of \(X_{3}\to-\infty\), we find
\begin{equation*}
  \Pi^{(k)}\left(X_{3}^{\prime},Y_{3}^{\prime}\right)=\frac{-\left(c_{\rmA\rmX}c_{\rmX\rmB}+e_{\rmA}c_{\rmY\rmX}\right)}{(c_{\rmA\rmX}+c_{\rmA\rmY})c_{\rmX\rmB}+e_{\rmA}c_{\rmX\rmY}},
\end{equation*}
which is the value of the limit \(\lim_{\vartheta\nearrow\vartheta_{\rms}}f_{\rmks}\left(\vartheta\right)\) in \eqref{eqn:ks_lim_left}.

\subsubsection{The \(\Delta\)-clique network}

In the case of the \(\Delta\)-clique network, the expressions for \(X_{3}^{\prime}\) and \(Y_{3}^{\prime}\) differ on either side of the switching manifold by only an \(\order{1}\) global constant, since, in both cases, the trajectory follows the connection \(X\to Y\) but does not strike \(\poinin[\mrm{B}]{\mrm{C}}\) near \(X\). Therefore, in both cases, substituting the expressions in \eqref{eqn:Delta_x2_out_C} and \eqref{eqn:Delta_theta_3_out_C} into \eqref{eqn:X3_prime} and \eqref{eqn:Y3_prime}, and the resulting expressions into the projection \(\Pi^{(k)}\), and rearranging gives
\begin{equation*}
  \begin{aligned}
    &\qquad\qquad\qquad\qquad\Pi^{(k)}\left(X_{3}^{\prime},Y_{3}^{\prime}\right)=\\
    &\frac{
      -c_{\rmA\rmX}\left(-\frac{e_{\rmY\rmA}}{c_{\rmB}}\frac{1}{T_{\rmB}}+c_{\rmY\rmB}\right)-e_{\rmA}\left(\frac{c_{\rmY\rmX}}{c_{\rmB}}\frac{1}{T_{\rmB}}+c_{\rmY\rmX}\right)
      }{
        (c_{\rmA\rmX}+c_{\rmA\rmY})\left(-\frac{e_{\rmY\rmA}}{c_{\rmB}}\frac{1}{T_{\rmB}}+c_{\rmY\rmB}\right)+e_{\rmA}\left(\frac{c_{\rmY\rmX}}{c_{\rmB}}\frac{1}{T_{\rmB}}+c_{\rmY\rmX}\right)
    }.
  \end{aligned}
\end{equation*}
The expression is valid, and can be evaluated, at any point on the line \(X_{3}=\log a+\frac{e_{\rmB\rmX}}{e_{\rmB\rmY}}Y_{3}\) for \(0>\log a>\log\left(1-\epsilon\right)\) and \(0<\log a<\log\left(\frac{1}{1-\epsilon}\right)\). Any such value of \(\log a\) lies inside the excluded cusp. In both cases, \(\delta\) is again fixed and \(T_{\rmB}\to\infty\), giving us, in both cases,
\begin{equation*}
  \Pi^{(k)}\left(X_{3}^{\prime},Y_{3}^{\prime}\right)=\frac{-\left(c_{\rmA\rmX}c_{\rmX\rmB}+e_{\rmA}c_{\rmY\rmX}\right)}{(c_{\rmA\rmX}+c_{\rmA\rmY})c_{\rmX\rmB}+e_{\rmA}c_{\rmX\rmY}},
\end{equation*}
which is the value of the limit \(\lim_{\vartheta\to\vartheta_{\rms}}f_{\Delta}\left(\vartheta\right)\) in \eqref{eqn:delta_lim}.

\subsubsection{Summary}

We find that the discontinuity of \(f_{\rmks}\) at \(\vartheta_{\rms}\) is a result of two factors: first, that evaluating the limit \(\vartheta\to\vartheta_{\rms}\) requires us to evaluate the flow in the limit as a trajectory approaches the network; and, second, that there is a separatrix near the network that creates a discontinuity in the flow in the limit as a trajectory approaches the network. In the case of the \(\Delta\)-clique network, however, although the first of these two factors is still required, no discontinuity emerges in the flow in the limit as a trajectory approaches the network because there is no separatrix near the \(\Delta\)-clique network. Therefore, the projected map \(f_{\Delta}\) is continuous at \(\vartheta_{\rms}\).

\section{Discussion}\label{sec:disc}

In \cite{article}, we constructed projected maps to analyse the dynamics of trajectories near three heteroclinic networks in \(\R^{4}\), including the Kirk--Silber and \(\Delta\)-clique networks. In that paper, we observed that the projected map of the Kirk--Silber network was generically discontinuous on its switching manifold, but that the projected map of the \(\Delta\)-clique network was continuous on its switching manifold. In this paper, we have investigated and explained the origin of this discontinuity in terms of the dynamical properties of the flow that contains the respective networks.

The projected maps studied in \cite{article} are derived from the transition matrices of return maps that have the form in \eqref{eqn:gen_ret_map}. In particular, the projected map is derived from the return map to the incoming cross-section \(\poinin[A]{B}\). Because \(B\) is a splitting equilibrium, a subset of \(\poinin[A]{B}\) must be excluded from the domain of definition of \(\Phi_{\rmB}\). Points in this subset give rise to trajectories that do not remain sufficiently close to the heteroclinic network for certain linearisations to remain valid approximations of the flow near the network.

In \Sref{sec:ret_map}, we followed the process adopted in \cite{kirk_lane_postlethwaite_rucklidge_silber_2010,kirk_postlethwaite_rucklidge_2012} and constructed near the two networks a \textit{completed return map} \(\widehat{\Phi}_{\rmB}\colon\poinin[A]{B}\to\poinin[A]{B}\). This return map accounts for all trajectories sufficiently close to the network, including those that begin in the excluded cusp \(\Gamma_{\rmc}\). This map is constructed by linearising the flow near \(A\) and \(B\), and along an invariant curve \(C\) that contains \(X\) and \(Y\). Unlike the typical way of constructing return maps---in which explicit expressions for all components of the map are readily derivable---we know very little about the components of this map. The expressions for the linearised dynamics along \(C\) are particularly opaque. We have complemented the analysis in \cite{kirk_lane_postlethwaite_rucklidge_silber_2010,kirk_postlethwaite_rucklidge_2012} by considering the asymptotic behaviour of components of the completed return map. Crucially, we can deduce that the \(x_{1}^{\invar{C}}\)-dependent function \(\theta_{3}^{\outvar{C}}\left(\theta_{3}^{\invar{C}}\right)\) is discontinuous in the limit \(x_{3}^{\invar{C}}\to0\) for the Kirk--Silber network, but continuous for the \(\Delta\)-clique network.

In \Sref{sec:anal} we considered the logarithmic form of \(\widehat{\Phi}_{\rmB}\) and defined projected maps from these expressions. Due to the nature of the completed return map, the value of these projected maps depended on where the simplex \(S^{(k)}\) was defined. We observed numerically that the projected map \(f^{(k)}\) approaches the original projected maps \(f_{\rmks}\) and \(f_{\Delta}\) in the limit as \(k\to\infty\). We then determined the value of the completed return map for certain trajectories near each network, and applied the projection to these expressions. We were then able to show that the discontinuity of \(f_{\rmks}\) emerges as a combination of two factors. First, there exists a separatrix near the Kirk--Silber network formed by the stable manifold of the equilibrium \(P\), and this separatrix creates a discontinuity in the flow of trajectories near the network in the limit as those trajectories approach the network. Second, the method used to define the projected map requires us to take the limit as trajectories approach the network in order to evaluate the limit \(\vartheta\to\vartheta_{\rms}\). In the case of the \(\Delta\)-clique network, no such separatrix exists, and so the projected map is continuous at the switching manifold. Our methodology also allowed us to demonstrate why the limits in \eqref{eqn:ks_lim_right}, \eqref{eqn:ks_lim_left}, and \eqref{eqn:delta_lim} attain those particular values.

In \cite{article}, we also studied the tournament network, which is the union \(\cyc{X}\cup\cyc{Y}\cup\cyc{XY}\). Although the tournament network also contains a subnetwork whose graph representation is equivalent to the Kirk--Silber network, the nearby dynamics does not feature a structure equivalent to the equilibrium \(P\), and so also does not contain the separatrix formed by its stable manifold. It instead contains two \(\Delta\)-cliques, with trajectories leaving a neighbourhood containing the \(\Delta\)-clique near a common equilibrium. A small extension of the map derived in \Sref{sec:ret_map} and the analysis in \Sref{sec:anal} shows that the projected map of the tournament network is continuous on both of its two switching manifolds.

We proved in \cite{article} that any trajectory asymptotic to the Kirk--Silber, \(\Delta\)-clique, or tournament network is asymptotic to one component cycle, and can switch between two given cycles only once, and only in a particular direction that depends on parameter values. As such, for these three networks composed of only four equilibria, the continuity or discontinuity of the projected map at the switching manifold does not correspond to a qualitative difference in the long-term dynamics observed near these networks.

However, complicated dynamics have been observed near networks with more than four equilibria. For example, in \cite{postlethwaite_rucklidge_2022}, Postlethwaite and Rucklidge observed that near a network composed of five equilibria, regions of parameter space that correspond to stable regular cycling between subcycles of the network form chains with codimension-\(2\) points of zero width. These chains of mode-locking regions have been studied extensively in piecewise-smooth dynamical systems (see \cite[\S~9]{simpson_2016} for an introductory overview), and were originally likened to resembling a string of sausages \cite{yang_hao_1987}. Such strings-of-sausages are typical of piecewise-linear \textit{continuous} maps. In contrast, Podvigina studied in \cite{podvigina_2023} a two-cycle network composed of two \(C_{4}^{-}\) cycles, which are of length \(4\), joined along a common connection, similar to the Kirk--Silber network. They prove that a trajectory may regularly cycle between these two cycles, switching back and forth an infinite number of times. Some of our preliminary numerical results of the stability regions of various cycling patterns near Podvigina's two-cycle network form patterns in parameter space fundamentally different from those observed by Postlethwaite and Rucklidge \cite{postlethwaite_rucklidge_2022}. The projected map of the network studied by Postlethwaite and Rucklidge is continuous on certain (relevant) subsets of the switching manifold, as are certain relevant compositions of maps on certain subsets of the switching manifold. The projected map of Podvigina's two-cycle network is discontinuous on its switching manifold, like the Kirk--Silber network.

Although we leave a detailed study of the nature of regular or irregular cycling near these and other networks as future work, we suspect that the continuity of the projected maps of these networks will affect the structures and patterns of regular and irregular cycling that exist in parameter space. In this paper, we have developed a methodology to explain when the projected map of a network is continuous on its switching manifold, and it will, in future, allow us to decide the continuity of a projected map without constructing it. In turn, we may be able to derive aspects of patterns in parameter space that define the possible behaviour of trajectories from the structure of the heteroclinic network.

\section*{Acknowledgments}

DCGD thanks Vivien Kirk for her guidance, contributions, and feedback provided on a draft of this work as part of his doctoral thesis. DCGD also thanks the School of Mathematics at the University of Leeds, where the initial stages of this work were done, for their hospitality.
This research was funded by the Marsden Fund Council from New Zealand Government funding managed by the Royal Society Te Ap{\=a}rangi (Grant No. 21-UOA-048), and by the Engineering and Physical Sciences Research Council (Grant No. EP/V014439/1).
For the purpose of Open Access, the authors have applied a Creative Commons Attribution (CC BY) license to any Author Accepted Manuscript version arising from this submission.

\section*{Data Availability Statement}

No new data were created or analysed in this study.

\appendix
\section{Numerics}\label{app:numerics}
\setcounter{section}{1}

In this appendix, we describe the numerics used to produce \Fref{fig:theta_3_C_out_v_theta_3_C_in} in \ref{app:numerics:ssec:theta_in_theta_out} and \Fref{fig:proj_map_approx} in \ref{app:numerics:ssec:projected_k}.

For reference, the ODEs \eqref{eqn:ODEs} with the necessary constant terms to ensure the existence of a Kirk--Silber network are
\begin{equation}\label{eqn:ks_ODEs}
  \begin{aligned}
    \dot{x}_{1}&=x_{1}\left(1 - \left\lVert x\right\rVert_{2}^{2} - c_{\rmB}x_{2}^{2} + e_{\rmX\rmA}x_{3}^{2} + e_{\rmY\rmA}y_{3}^{2}\right),\\
    \dot{x}_{2}&=x_{2}\left(1 - \left\lVert x\right\rVert_{2}^{2} + e_{\rmA}x_{1}^{2} - c_{\rmX\rmB}x_{3}^{2} - c_{\rmY\rmB}y_{3}^{2}\right),\\
    \dot{x}_{3}&=x_{3}\left(1 - \left\lVert x\right\rVert_{2}^{2} - c_{\rmA\rmX}x_{1}^{2} + e_{\rmB\rmX}x_{2}^{2} - c_{\rmY\rmX}y_{3}^{2}\right),\\
    \dot{y}_{3}&=y_{3}\left(1 - \left\lVert x\right\rVert_{2}^{2} - c_{\rmA\rmY}x_{1}^{2} + e_{\rmB\rmY}x_{2}^{2} - c_{\rmX\rmY}x_{3}^{2}\right),
  \end{aligned}
\end{equation}
where each constant is a positive real number.

The ODEs \eqref{eqn:ODEs} with the necessary constant terms to ensure the existence of a \(\Delta\)-clique network are
\begin{equation}\label{eqn:delta_clique_ODEs}
  \begin{aligned}
    \dot{x}_{1}&=x_{1}\left(1 - \left\lVert x\right\rVert_{2}^{2} - c_{\rmB}x_{2}^{2} - c_{\rmX\rmA}x_{3}^{2} + e_{\rmY\rmA}y_{3}^{2}\right),\\
    \dot{x}_{2}&=x_{2}\left(1 - \left\lVert x\right\rVert_{2}^{2} + e_{\rmA}x_{1}^{2} - c_{\rmX\rmB}x_{3}^{2} - c_{\rmY\rmB}y_{3}^{2}\right),\\
    \dot{x}_{3}&=x_{3}\left(1 - \left\lVert x\right\rVert_{2}^{2} - c_{\rmA\rmX}x_{1}^{2} + e_{\rmB\rmX}x_{2}^{2} - c_{\rmY\rmX}y_{3}^{2}\right),\\
    \dot{y}_{3}&=y_{3}\left(1 - \left\lVert x\right\rVert_{2}^{2} - c_{\rmA\rmY}x_{1}^{2} + e_{\rmB\rmY}x_{2}^{2} + e_{\rmX\rmY}x_{3}^{2}\right),
  \end{aligned}
\end{equation}
where each constant is a positive real number.

\subsection{Approximation of \(\theta_{3}^{\outvar{C}}\left(\theta_{3}^{\invar{C}}\right)\)}\label{app:numerics:ssec:theta_in_theta_out}

To derive a numerical approximation of \(\theta_{3}^{\outvar{C}}\left(\theta_{3}^{\invar{C}}\right)\), we begin with a point \(x\left(0\right)\in\poinin[\mrm{B}]{\mrm{C}}\), which implies that \(x_{2}\left(0\right)=h\) and \(x_{1}\left(0\right)<h\). As such, \(x_{1}\left(0\right),x_{2}\left(0\right)\ll 1\). To minimise round-off error when computing \(\dot{x}_{1}\) and \(\dot{x}_{2}\), we use the change of coordinates \(X_{1}\coloneq\log x_{1}\) and \(X_{2}\coloneq\log x_{2}\), and compute \(\dot{X}_{1}\) and \(\dot{X}_{2}\). Since \(x_{1}\left(t\right),x_{2}\left(t\right)\ll 1\) while the trajectory is in a small neighbourhood of \(C\), we can ignore them from the right-hand side of the ODEs. Thus, we have the following system of ODEs for the Kirk--Silber network
\begin{equation}\label{eqn:ks_ODEs_near_C}
  \begin{aligned}
    \dot{X}_{1}&=1 - \left(x_{3}^{2}+y_{3}^{2}\right) + e_{\rmX\rmA}x_{3}^{2} + e_{\rmY\rmA}y_{3}^{2},\\
    \dot{X}_{2}&=1 - \left(x_{3}^{2}+y_{3}^{2}\right) - c_{\rmX\rmB}x_{3}^{2} - c_{\rmY\rmB}y_{3}^{2},\\
    \dot{x}_{3}&=x_{3}\left(1 - \left(x_{3}^{2}+y_{3}^{2}\right) - c_{\rmY\rmX}y_{3}^{2}\right),\\
    \dot{y}_{3}&=y_{3}\left(1 - \left(x_{3}^{2}+y_{3}^{2}\right) - c_{\rmX\rmY}x_{3}^{2}\right),
  \end{aligned}
\end{equation}
and the following system of ODEs for the \(\Delta\)-clique network
\begin{equation}\label{eqn:delta_clique_ODEs_near_C}
  \begin{aligned}
    \dot{X}_{1}&=1 - \left(x_{3}^{2}+y_{3}^{2}\right) - c_{\rmX\rmA}x_{3}^{2} + e_{\rmY\rmA}y_{3}^{2},\\
    \dot{X}_{2}&=1 - \left(x_{3}^{2}+y_{3}^{2}\right) - c_{\rmX\rmB}x_{3}^{2} - c_{\rmY\rmB}y_{3}^{2},\\
    \dot{x}_{3}&=x_{3}\left(1 - \left(x_{3}^{2}+y_{3}^{2}\right) - c_{\rmY\rmX}y_{3}^{2}\right),\\
    \dot{y}_{3}&=y_{3}\left(1 - \left(x_{3}^{2}+y_{3}^{2}\right) + e_{\rmX\rmY}x_{3}^{2}\right).
  \end{aligned}
\end{equation}
For an initial condition \(\left(X_{1}\left(0\right),\log h,x_{3}\left(0\right),y_{3}\left(0\right)\right)\in\poinin[\mrm{B}]{\mrm{C}}\), with \(X_{1}\left(0\right)<0\), these ODEs can be integrated until \(X_{1}=0\), at which time \(\left(X_{1}\left(t\right),X_{2}\left(t\right),x_{3}\left(t\right),y_{3}\left(t\right)\right)\in\poinout[\mrm{A}]{\mrm{C}}\).

To produce \Fref{fig:theta_3_C_out_v_theta_3_C_in}, we integrated \eqref{eqn:ks_ODEs_near_C} and \eqref{eqn:delta_clique_ODEs_near_C} with a range of initial conditions. In particular, we chose first \(X_{1}\left(0\right)\in\left\{-10^{-12},-1,-2,-4,-8\right\}\). For each value of \(X_{1}\left(0\right)\), we then sampled \(1000\) values of the angle \(\theta_{3}\) from the open interval \(\left(\frac{\pi}{2}\cdot\frac{1}{1001},\frac{\pi}{2}-\frac{\pi}{2}\cdot\frac{1}{1001}\right)\), in particular, setting \(\theta_{3}^{(\ell)}=\ell\frac{\pi}{2}\cdot\frac{1}{1001}\) for all positive integers \(1\leq\ell\leq 1000\). (Note we deliberately exclude \(\theta_{3}\left(0\right)=0,\frac{\pi}{2}\), as we trivially know the angle on \(\poinout[\mrm{A}]{\mrm{C}}\) in these cases.) We then used as initial conditions \(x_{3}\left(0\right)=\cos\theta_{3}^{(\ell)}\) and \(y_{3}\left(0\right)=\sin\theta_{3}^{(\ell)}\).

We used the eighth-order Dormand--Prince method \cite[\S~II.5]{harrier_norsett_wanner_2008} of the Python scientific computing library SciPy \cite{scipy} to integrate the ODEs \eqref{eqn:ks_ODEs_near_C} and \eqref{eqn:delta_clique_ODEs_near_C} with these initial conditions until the solver detected that \(X_{1}\left(t\right)=0\) at \(t=T_{C}\). We then calculated \(\theta_{3}\left(T_{C}\right)\coloneq\arctan\frac{y_{3}\left(T_{C}\right)}{x_{3}\left(T_{C}\right)}\). The results, plotting \(\theta_{3}\left(T_{C}\right)\) against the initial angle \(\theta_{3}^{(\ell)}\), are shown in \Fref{fig:theta_3_C_out_v_theta_3_C_in}.

\subsection{Approximations of \(f_{\mathrm{KS}}^{(k)}\) and \(f_{\Delta}^{(k)}\)}\label{app:numerics:ssec:projected_k}

For a given \(k>1\), we begin with a subset of \(50\) equispaced points in \(S^{(k)}\) defined by
\begin{equation*}
  \left\{-X_{3}^{(\ell)}=\frac{k}{51}\ell\ \bigg|\ 1\leq\ell\leq 50\right\}.
\end{equation*}
The point \(\left(X_{3}^{(\ell)},Y_{3}^{(\ell)}\right)\in S^{(k)}\), where \(Y_{3}^{(\ell)}=-k-X_{3}^{(\ell)}\), lies on \(\poinin[\mrm{A}]{\rmB}\) expressed in logarithmic coordinates. Therefore, as an initial condition, we know for that point that \(x_{1}\left(0\right)=h\) and \(|x_{2}\left(0\right)-h|=0\). For \Fref{fig:proj_map_approx}, we used \(h=10^{-2}\), and, for simplicity, we chose \(x_{2}\left(0\right)=1\).

To minimise round-off error, we integrate the ODEs \eqref{eqn:ks_ODEs} and \eqref{eqn:delta_clique_ODEs} in logarithmic coordinates, with the initial condition \(X\left(0\right)=\left(\log 10^{-2},0,X_{3}^{(\ell)},Y_{3}^{(\ell)}\right)\). Again, we used the eighth-order Dormand--Prince method of SciPy until the solver detected that the trajectory intersected \(\poinin[\mrm{A}]{\rmB}\) again. This intersection point is a numerical approximation of \(\left(X_{3}^{\prime},Y_{3}^{\prime}\right)\). The projection \(\Pi^{(k)}\) of this point was then calculated. This projection provided a numerical approximation of \(f_{\mathrm{KS}}^{(k)}\) or \(f_{\Delta}^{(k)}\) at the point \(\vartheta=X_{3}^{{(\ell)}}\). We last rescaled the set \(S^{(k)}\) to unit length. These numerical approximations are presented in Figures~\ref{fig:proj_map_approx} and \ref{fig:proj_map_approx_sm}.

We note that a point \(\left(X_{3}^{(\ell)},Y_{3}^{(\ell)}\right)\) on the line \(X_{3}+Y_{3}=-k\) may satisfy either \(\exp X_{3}^{(\ell)}>\log h\) or \(\exp Y_{3}^{(\ell)}>\log h\), in which case the point is not actually in \(\mcl{D}\). In this case, we exclude the point from consideration as it does not lie in the set \(T^{(k)}\), and for this reason the domain of \(f^{(k)}\) lengthens as \(k\) increases.

\section{Parameter values}\label{app:parms}

\Fref{fig:proj_maps:kirk_silber_net} is calculated with the projected map in \eqref{eqn:ks_proj_map} with parameter values:

\begin{table}[H]
  \centering
  \caption{Parameter values for \Fref{fig:proj_maps:delta_clique_net}}
  \label{tbl:parms_proj_map_ks}
  \begin{tabular}{| c c || c c || c c || c c |}
    \hline
    \multicolumn{2}{|c|}{\(A\)} & \multicolumn{2}{|c|}{\(B\)} & \multicolumn{2}{|c|}{\(X\)} & \multicolumn{2}{|c|}{\(Y\)} \\
    \hline
    \hline
    \(e_{\rmA}\)     & \(0.8\) & \(c_{\rmB}\)     & \(0.3\) & \(e_{\rmX\rmA}\) & \(0.8\) & \(e_{\rmY\rmA}\) & \(0.8\) \\
    \(c_{\rmA\rmX}\) & \(2.9\) & \(e_{\rmB\rmX}\) & \(1.2\) & \(c_{\rmX\rmB}\) & \(2.9\) & \(c_{\rmY\rmB}\) & \(2.3\) \\
    \(c_{\rmA\rmY}\) & \(2\)   & \(e_{\rmB\rmY}\) & \(0.8\) & \(c_{\rmX\rmY}\) & \(5\)   & \(c_{\rmY\rmX}\) & \(3\) \\
    \hline
  \end{tabular}
\end{table}

\Fref{fig:proj_maps:delta_clique_net} is calculated with the projected map in \eqref{eqn:delta_proj_map} with parameter values:

\begin{table}[H]
  \centering
  \caption{Parameter values for \Fref{fig:proj_maps:delta_clique_net}}
  \label{tbl:parms_proj_map_delta}
  \begin{tabular}{| c c || c c || c c || c c |}
    \hline
    \multicolumn{2}{|c|}{\(A\)} & \multicolumn{2}{|c|}{\(B\)} & \multicolumn{2}{|c|}{\(X\)} & \multicolumn{2}{|c|}{\(Y\)} \\
    \hline
    \hline
    \(e_{\rmA}\)     & \(7.48\)  & \(c_{\rmB}\)     & \(4.74\) & \(c_{\rmX\rmA}\) & \(0.25\)  & \(e_{\rmY\rmA}\) & \(3.58\) \\
    \(c_{\rmA\rmX}\) & \(0.69\)  & \(e_{\rmB\rmX}\) & \(1.2\)  & \(c_{\rmX\rmB}\) & \(10.21\) & \(c_{\rmY\rmB}\) & \(2.36\) \\
    \(c_{\rmA\rmY}\) & \(10.77\) & \(e_{\rmB\rmY}\) & \(0.8\)  & \(e_{\rmX\rmY}\) & \(0.57\)  & \(c_{\rmY\rmX}\) & \(8.68\) \\
    \hline
  \end{tabular}
\end{table}

\Fref{fig:theta_3_C_out_v_theta_3_C_in:kirk_silber_net} is calculated as outlined in \ref{app:numerics:ssec:theta_in_theta_out}, with parameter values:

\begin{table}[H]
  \centering
  \caption{Parameter values for \Fref{fig:theta_3_C_out_v_theta_3_C_in:kirk_silber_net}}
  \label{tbl:parms_CinCout_ks}
  \begin{tabular}{| c c || c c |}
    \hline
    \multicolumn{2}{|c|}{\(X\)} & \multicolumn{2}{|c|}{\(Y\)} \\
    \hline
    \hline
    \(e_{\rmX\rmA}\) & \(0.8\) & \(e_{\rmY\rmA}\) & \(0.9\) \\
    \(c_{\rmX\rmB}\) & \(1.1\) & \(c_{\rmY\rmB}\) & \(1.2\) \\
    \(c_{\rmX\rmY}\) & \(1.2\) & \(c_{\rmY\rmX}\) & \(0.8\) \\
    \hline
  \end{tabular}
\end{table}

\Fref{fig:theta_3_C_out_v_theta_3_C_in:delta_clique_net} is calculated as outlined in \ref{app:numerics:ssec:theta_in_theta_out}, with parameter values:

\begin{table}[H]
  \centering
  \caption{Parameter values for \Fref{fig:theta_3_C_out_v_theta_3_C_in:delta_clique_net}}
  \label{tbl:parms_CinCout_delta}
  \begin{tabular}{| c c || c c |}
    \hline
    \multicolumn{2}{|c|}{\(X\)} & \multicolumn{2}{|c|}{\(Y\)} \\
    \hline
    \hline
    \(c_{\rmX\rmA}\) & \(0.4\) & \(e_{\rmY\rmA}\) & \(0.9\) \\
    \(c_{\rmX\rmB}\) & \(2.5\) & \(c_{\rmY\rmB}\) & \(1.1\) \\
    \(e_{\rmX\rmY}\) & \(0.9\) & \(c_{\rmY\rmX}\) & \(0.2\) \\
    \hline
  \end{tabular}
\end{table}

Figures~\ref{fig:proj_map_approx:ks} and \ref{fig:proj_map_approx_sm:ks} are calculated with the completed projected map \(f^{(k)}_{\mrm{KS}}\) and the following parameter values:

\begin{table}[H]
  \centering
  \caption{Parameter values for Figures~\ref{fig:proj_map_approx:ks} and \ref{fig:proj_map_approx_sm:ks}}
  \label{tbl:parms_fk_ks}
  \begin{tabular}{| c c || c c || c c || c c |}
    \hline
    \multicolumn{2}{|c|}{\(A\)} & \multicolumn{2}{|c|}{\(B\)} & \multicolumn{2}{|c|}{\(X\)} & \multicolumn{2}{|c|}{\(Y\)} \\
    \hline
    \hline
    \(e_{\rmA}\)     & \(0.8\) & \(c_{\rmB}\)     & \(0.8\) & \(e_{\rmX\rmA}\) & \(0.8\) & \(e_{\rmY\rmA}\) & \(0.9\) \\
    \(c_{\rmA\rmX}\) & \(1.2\) & \(e_{\rmB\rmX}\) & \(1\)   & \(c_{\rmX\rmB}\) & \(1.2\) & \(c_{\rmY\rmB}\) & \(1.2\) \\
    \(c_{\rmA\rmY}\) & \(1.2\) & \(e_{\rmB\rmY}\) & \(0.8\) & \(c_{\rmX\rmY}\) & \(0.2\) & \(c_{\rmY\rmX}\) & \(0.1\) \\
    \hline
  \end{tabular}
\end{table}

Figures~\ref{fig:proj_map_approx:delta} and \ref{fig:proj_map_approx_sm:delta} are calculated with the completed projected map \(f^{(k)}_{\Delta}\) and the following parameter values:

\begin{table}[H]
  \centering
  \caption{Parameter values for Figures~\ref{fig:proj_map_approx:delta} and \ref{fig:proj_map_approx_sm:delta}}
  \label{tbl:parms_fk_delta}
  \begin{tabular}{| c c || c c || c c || c c |}
    \hline
    \multicolumn{2}{|c|}{\(A\)} & \multicolumn{2}{|c|}{\(B\)} & \multicolumn{2}{|c|}{\(X\)} & \multicolumn{2}{|c|}{\(Y\)} \\
    \hline
    \hline
    \(e_{\rmA}\)     & \(1.2\) & \(c_{\rmB}\)     & \(1.3\) & \(c_{\rmX\rmA}\) & \(0.8\) & \(e_{\rmY\rmA}\) & \(0.9\) \\
    \(c_{\rmA\rmX}\) & \(0.3\) & \(e_{\rmB\rmX}\) & \(1.2\) & \(c_{\rmX\rmB}\) & \(2.3\) & \(c_{\rmY\rmB}\) & \(1\)   \\
    \(c_{\rmA\rmY}\) & \(1.8\) & \(e_{\rmB\rmY}\) & \(0.9\) & \(e_{\rmX\rmY}\) & \(0.5\) & \(c_{\rmY\rmX}\) & \(2.4\) \\
    \hline
  \end{tabular}
\end{table}

\section*{References}
\bibliographystyle{unsrt}
\bibliography{main}
\end{document}